# DSm VECTOR SPACES OF REFINED LABELS

W. B. Vasantha Kandasamy
Florentin Smarandache

**2011**

# CONTENTS

# PREFACE

The study of DSm linear algebra of refined labels have been done by Florentin Smarandache, Jean Dezert, and Xinde Li.

In this book the authors introduce the notion of DSm vector spaces of refined labels. The reader is requested to refer the paper as the basic concepts are taken from that paper [35].

This book has six chapters. The first one is introductory in nature just giving only the needed concepts to make this book a self contained one. Chapter two introduces the notion of refined plane of labels, the three dimensional space of refined labels DSm vector spaces. Clearly any n-dimensional space of refined labels can be easily studied as a matter of routine.

This chapter defines the notion of DSm vector space of refined labels and it contains around 7 definitions, 47 examples and 33 theorems. The new notion of different types of special DSm vector spaces are described with 82 examples in chapter three. DSm semivector spaces of ordinary labels and refined labels are introduced and studied in chapter four. Chapter five suggests some applications of these new structures. Over 125 problems are given in chapter six; some of which are simple and some of them are at research level.

We thank Dr. K.Kandasamy for proof reading and being extremely supportive.

W.B.VASANTHA KANDASAMY
FLORENTIN SMARANDACHE

**Chapter One**

# INTRODUCTION

In this chapter we just recall the new notions of refined labels and linear algebra of refined labels. We mainly do this to make the book a self contained one. For these concepts are used to build new types of linear algebras. For more about these concepts please refer [7, 34-5]. Let $L_1, L_2, \ldots, L_m$ be labels, where m ≥ 1 is an integer. The set of labels are extended by using $L_0$ to be the minimal or minimum label and $L_{m+1}$ to be the maximum label.

We say the labels are equidistant if the qualitative distant between any two consecutive labels is the same, we get an exact qualitative result and a qualitative basic belief assignment (bba) is considered normalized, if the sum of all its qualitative masses is equal to $L_{max} = L_{m+1}$. If the labels are not equidistant, we still can use all qualitative operators defined in the Field and the Linear Algebra of Refined Labels (FLARL), but the qualitative result is approximate and a qualitative bba is considered quasi-normalized if the sum of all its masses is equal to $L_{max}$. We consider a relation of order defined on these labels which can be "smaller", "less in equality" "lower" etc; $L_1 < L_2 < \ldots < L_m$. Connecting them to the classical interval [0, 1] we have so,

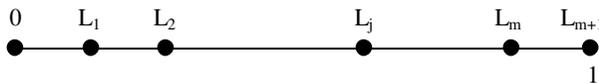

$0 \equiv L_1 < L_2 < \ldots < L_i < \ldots < L_m < L_{m+1} \equiv 1$, and $L_i = \dfrac{i}{m+1}$

for $i \in \{0, 1, 2, \ldots, m, m+1\}$.

Now we proceed onto define the notion of ordinary labels [7, 34-5]. The set of labels $\tilde{L} \triangleq \{L_0, L_1, L_2, \ldots, L_i, \ldots, L_m, L_{m+1}\}$ whose indexes are positive integers between 0 and m+1 is called the set of 1-Tuple labels. We call a set of labels to be equidistant labels, if the geometric distance between any two consecutive $L_i$ and $L_{i+1}$ is the same, that is $L_{i+1} - L_i$ = constant for any i, $1 \leq i \leq m + 1$. A set of labels is said to be of non-equidistant if the distances between consecutive labels are not the same that is for $i \neq j$; $L_{j+1} - L_j \neq L_{i+1} - L_i$. We see the set of 1-Tuple labels is isomorphic with the numerical set $\{\dfrac{i}{m+1}$; i = 0, 1, \ldots, m+1\}$ through the isomorphism $f_{\tilde{L}}^{(L_i)} = \dfrac{i}{m+1}$; $0 \leq i \leq m+1$. Now we proceed onto recall the definition of refined labels [7, 34-5]. We just theoretically extend the set of labels $\tilde{L}$ on the left and right sides of the intervals [0, 1] towards $-\infty$ and respectively $+\infty$.

Thus $L_z \triangleq \left\{ \dfrac{j}{m+1} \,\middle|\, j \in Z \right\}$ where Z is the set of positive and negative integers including zero.

Thus $L_z = \{\ldots L_{-j}, \ldots, L_{-1}, 0, L_1, \ldots, L_j, \ldots\} = \{L_j \mid j \in Z\}$; ie., the set of extended Labels with positive and negative indices. Similarly we define $L_Q \triangleq \{L_q \mid q \in Q\}$ as the set of labels whose indices are rational or fractions. $L_Q$ is isomorphic $f_Q(L_q) = \frac{q}{m+1}$; $q \in Q$. Even more generally they define $L_R \triangleq \left\{ \frac{r}{m+1} \middle| r \in R \right\}$ where R is the set of real numbers. $L_R$ is isomorphic with R through the isomorphism $f_R(L_r) = \frac{r}{m+1}$ for any $r \in R$ [7, 34-5].

We will just recall the proof / definitions which makes {$L_R$, +, ×} a field called the DSm field of refined labels. For the first time the notion of decimal or refined labels that is labels with index as decimal is defined [7, 34-5]. For example $L_{3/2}$ which is $L_{1.5}$ means a label in the middle of the label interval [$L_1$, $L_2$]. They have theoretically introduced negative labels $L_{-i}$ which is equal to $-L_i$ that occur in qualitative calculations.

Further ($L_R$, +, ×, .) where '.' means scalar product is a commutative linear algebra over the field of real numbers R with unit element, and for which each non-null element is invertible with respect to multiplication of labels.

This is called DSm field and Linear Algebra of Refined labels (FLARL for short) [7, 34-5].

We just recall the definition of qualitative operators on FLARL. We consider a, b, c in R and labels $L_a = \frac{a}{m+1}$, $L_b = \frac{b}{m+1}$ and $L_c = \frac{c}{m+1}$.

$$L_a + L_b = \frac{a}{m+1} + \frac{b}{m+1} = \frac{a+b}{m+1}$$
$$= L_{a+b}.$$

Likewise $L_a - L_b = L_{a-b}$ as $\dfrac{a}{m+1} - \dfrac{b}{m+1} = \dfrac{a-b}{m+1}$.

Consider $L_a \cdot L_b = L_{(ab)/m+1}$

Since $\dfrac{a}{m+1} \cdot \dfrac{b}{m+1} = \dfrac{(ab)/m+1}{m+1}$.

The notion of scalar multiplication for any $\alpha, \beta$ in R is as follows:

$\alpha L_a = L_a \cdot \alpha = L_{\alpha a}$;

since $\alpha \cdot \dfrac{a}{m+1} = \dfrac{a \cdot \alpha}{m+1}$.

If $\alpha = -1$ then we get $L_a(-1) = -L_a = L_{-a}$.

Also $\dfrac{L_a}{\beta} = L_a \div \beta = 1/\beta \cdot L_a = L_{a/\beta}$ ($\beta \neq 0$), $\beta \in R$.

Now we proceed onto define vector division or division of labels.

$L_a \div L_b = L_{(a/b)(m+1)}$.

$\dfrac{a}{m+1} \div \dfrac{b}{m+1} = \dfrac{(a/b)m+1}{m+1}$.

Now we define scalar power.

$(L_a)^p = L_{a^p/(m+1)^{p-1}}$ since

$$\left(\frac{a}{m+1}\right)^p = \frac{a^p/(m+1)^{p-1}}{m+1} \text{ for all } p \in R.$$

We can also define scalar root

$$a \sqrt[k]{L_\alpha} = (L_a)^{1/k} = L_{a^{1/k}/(m+1)^{1/k-1}}$$

which results from replacing p = 1/k in the power formula for all k integer greater than or equal to two. Thus ($L_R$, +, ×) is a field and is isomorphic with set of real number (R, +, ×) is defined as the DSm field of refined reals [7, 34-5]. The field isomorphism being defined by $f_R : L_R \to R$ where $f_R (L_R) = \frac{r}{m+1}$ satisfies the axioms.

$f_R (L_a + L_b) = f_R (L_a) + f_R (L_b)$ since

$f_R (L_a + L_b) = f_R (L_{a+b}) = \frac{a+b}{m+1}$ and

$f_R (L_a) + f_R (L_b) = \frac{a}{m+1} + \frac{b}{m+1} = \frac{a+b}{m+1}.$

$f_R (L_a \times L_b) = f_R (L_a) \cdot f_R (L_b)$

since

$f_R (L_a \times L_b) = f_R (L_{(ab)/(m+1)})$

$= \frac{ab}{m+1}$ and $f_R (L_a) \cdot f_R (L_b) = \frac{a}{m+1} \cdot \frac{b}{m+1} = \frac{ab}{(m+1)^2}.$

$(L_R, +, .)$ is a vector space of refined labels over the field of real numbers R since $(L_R, +)$ is a commutative group and the scalar multiplication which is an external operation.

Consider $1.L_a = L_{1.a} = L_a$.

If $\alpha, \beta \in R$ that $(\alpha.\beta) L_a = \alpha (\beta L_a)$ since both left and right sides are equal we see $(\alpha .\beta) L_a = L_{\alpha\beta a}$.

Further $\alpha (L_a + L_b) = \alpha L_a + \alpha L_b$

since $\alpha (L_a + L_b) = \alpha L_{a+b} = L_{\alpha(a+b)}$

$= L_{\alpha a} + L_{\alpha b} = L_{\alpha a} + L_{\alpha b} = \alpha L_a + \alpha L_b$.

Consider $(\alpha + \beta) = \alpha L_a + \beta L_b$,

since $(\alpha+\beta).L_a = L_{(\alpha+\beta)a} = L_{\alpha a+\beta a}$

$= L_{\alpha a} + L_{\beta \alpha} = \alpha L_a + \beta L_a$.

$(L_R, +, \times, .)$ is a linear algebra of Refined Labels over the field R of real numbers, called DSm Linear Algebra of Refined Labels. (DSm – LARL for short), which is commutative with identity element which is $L_{m+1}$ for vector multiplication and whose non null elements (labels) are invertible with respect to vector multiplication. This occurs since $(L_R, +, .)$ is a vector space $(L_R, \times)$ is a commutative group, the set of scalars R is well known as a field.

Clearly vector multiplication is associative. For consider $L_a \times (L_b \times L_c) = (L_a \times L_b) \times L_c$. To prove associativity, we know $L_a \times (L_b \times L_c) = L_a \times L_{bc/m+1} = L_{a.b.c/(m+1)^2}$ while $(L_a \times L_b) \times L_c = (L_{ab/m+1}) L_c = L_{a.b.c/(m+1)^2}$.

Hence the claim.

The vector multiplication is distributive with respect to addition.

Consider $L_a \times (L_b + L_c)$; to show $L_a \times (L_b + L_c) = L_a \times L_b + L_a \times L_c$ since $L_a \times (L_b + L_c) = L_a \times L_{b+c} = L_{a(b+c)/m+1}$ ----- (I)

Consider $L_a \times L_b + L_a \times L_c = L_{ab/m+1} + L_{ac/m+1} = L_{(ab+ac)/m+1}$

$= L_{a(b+c)/m+1}$.

Consequently we have $(L_a + L_b) \times L_c = L_a \times L_c + L_b \times L_c$.

Consider $(L_a + L_b) \times L_c = (L_{a+b}) \times L_c$
$= L_{(a+b)c/m+1} = L_{ac+bc/m+1}$.
$= L_{ac/m+1} + L_{bc/m+1}$
$= L_a \times L_c + L_b \times L_c$.

Finally we show $\alpha (L_a \times L_b) = (\alpha L_a) \times L_b = L_a \times (\alpha L_b)$ for all $\alpha \in R$.

Consider $\alpha (L_a \times L_b) = \alpha \cdot L_{ab/m+1}$
$= L_{\alpha(ab)/m+1} = L_{(\alpha b)b/m+1}$
$= L_{a(\alpha b)/m+1}$ as $L_{\alpha ab/m+1}$
$= L_{(\alpha a)b/m+1}$     $= L_{\alpha a} \times L_b$
$= \alpha L_a \times L_b$
$= L_a \times \alpha L_b$
$= L_a \times L_{\alpha b}$ (The argument is a matter of routine).

The unitary element for vector multiplication is $L_{m+1}$.

For all $a \in R$; $L_a \times L_{m+1} = L_{m+1} \times L_a = L_{a(m+1)/m+1} = L_a$

All $L_a \neq L_0$ are invertible with respect to vector multiplication and the inverse of $L_a$ is $(L_a)^{-1}$.

Consider $(L_a)^{-1} = L_{(m+1)^2/a} = \dfrac{1}{L_a}$

Since $L_a \times L_{a^{-1}} = L_a \times L_{(m+1)^2/a}$

$= L_{(a(m+1)^2/a)/m+1} = L_{m+1}$.

Hence the DSm linear algebra is a Division Algebra. DSm Linear Algebra is also a Lie Algebra since we can define a law

$(L_a, L_b) = [L_a, L_b] = L_a \times L_b - L_b \times L_a = L_0$ such that

$[L_a, L_a] = L_0$ and the Jacobi identity is satisfied.

$[L_a [L_b, L_c]] + [L_b [L_c, L_a]] + [L_c [L_a, L_b]] = L_0.$

Actually $(L_R, +, \times, .)$ is a field and therefore in particular a ring and any ring with the law: $[x, y] = xy - yx$ is a Lie algebra.

We can extend the field isomorphism $f_R$ to a linear algebra isomorphism by defining $f_R : R. L_R \to R.R.$ with $f_R(\alpha. L_{r_1}) = \alpha f_R(L_{r_1})$ since $f_R(\alpha L_{r_1}) = f_R(L_{(\alpha r_1)}) = \alpha r_1 / m+1 = \dfrac{\alpha r_1}{m+1} = \alpha.f_R(L_{r_1})$. Since $(R, +, .)$ is a trivial linear algebra over the field of reals $R$ and because $(L_R, +, .)$ is isomorphic with it through the above $f_R$ linear algebra isomorphism; it results that $(L_R, +, .)$ is also a linear algebra which is associative and commutative [7, 34-5].

We proceed onto recall more new operators like scalar-vector (mixed) addition, scalar-vector (mixed) subtraction,

scalar-vector (mixed) division, vector power and vector root on $(L_R, +, \times)$. For all $L_a \in L_R$ there exists a unique $r = \dfrac{a}{m+1}$ such that $L_a = r$ and inversely (reciprocally) for every $r \in R$ there exists a unique $L_a$ in $L_R$, $L_a \equiv L_{r(m+1)}$ such that $r = L_a$.

Let $\alpha \in R$ and $L_a \in L_R$.

$L_a + \alpha = \alpha + L_a = L_{a+\alpha(m+1)}$ since $L_a + \alpha = L_a + \dfrac{\alpha(m+1)}{m+1}$

$= L_a + L_{\alpha(m+1)} = L_{a+\alpha\ (m+1)}$.

This is defined as the scalar vector addition [7, 34-5].

On similar lines one can define scalar vector mixed subtraction

$L_a - \alpha = L_a - \alpha\ (m+1)$

since $L_a - \alpha = L_a - \dfrac{\alpha(m+1)}{m+1}$

$= L_a - L_{\alpha(m+1)} = L_{a-\alpha\ (m+1)}$.

$\alpha - L_a = L_{\alpha(m+1)-a}$.

Since $\alpha - L_a = \dfrac{\alpha(m+1)}{m+1} - L_a$.

$= L_{\alpha(m+1)} - L_a = L_{\alpha\ (m+1)-a}$. [7, 34-5]

Now we proceed onto recall the notion of scalar-vector mixed division.

$$L_a \div \alpha = \frac{L_a}{\alpha} = \frac{1}{\alpha} \cdot L_a = L_{\frac{a}{\alpha}} \quad \text{for } \alpha \neq 0 \in R, \text{ which is}$$

equivalent to the scalar multiplication $\left(\dfrac{1}{\alpha}\right) L_a$ where $\dfrac{1}{\alpha} \in R$.

$$\alpha \div L_\alpha = L_{a(m+1)^2/\alpha} \text{ since}$$

$$\alpha \div L_\alpha = \frac{\alpha(m+1)}{(m+1)} \div L_a = L_{\alpha(m+1)} \div L_a$$

$$= L_{(\alpha(m+1)/a) \cdot m+1}$$

$$= L_{\alpha(m+1)^2/a}.$$

Now we proceed onto define vector power in $(L_R, +, \times)$.

For $L_a$, $L_b$ in $L_R$; $(L_a)^{L_b} = L_{a^{b/m+1}/(m+1)^{b-m-1/m+1}}$,

since $(L_a)^{L_b} = (L_a)^{b/m+1} = L_{a^{b/m+1}/(m+1)^{\frac{b}{m+1}-1}}$

$= L_{a^{b/m+1}/(m+1)^{b-m-1/m+1}}$.

Now we recall the notion of vector roots

$\sqrt[L_b]{L_a} = L_{a^{m+1/b}/(m+1)^{m-b+1}}$,

since $\sqrt[L_b]{L_a} = (L_a)^{1/L_b} = (L_a)^{1/b/m+1}$

$= (L_a)^{m+1/b} = L_{a^{m+1/b}/(m+1)^{(m+1/b)-1}}$

$= L_{a^{m+1/b}/(m+1)^{m-b+1/b}}$.

For more refer [7, 34-5].

**Chapter Two**

# DSM VECTOR SPACES

We in this chapter define the new notion of Real Label Plane or Label Real Plane and Three dimensional real Label space. Also we define label vector spaces of different types and label linear algebras of different types and operations on them.

**DEFINITION 2.1:** *Let $L_{R \times R} = \{L_R \times L_R\} = \{(L_a, L_b) = \left(\dfrac{a}{m+1}, \dfrac{b}{m+1}\right) \mid a, b \in R; m \geq 2\}$; $L_{R \times R}$ is defined as the refined plane of labels. In fact any element of $L_{R \times R} = L_R \times L_R$ is an ordered pair set of labels.*

The following properties are direct and hence left as an exercise for the reader to prove.

**THEOREM 2.1:** *$L_{R \times R} = L_R \times L_R$ is isomorphic with real plane $R \times R$.*

**THEOREM 2.2:** $L_{R\times R} = L_R \times L_R$ is a commutative ring with unit under the operation + and ×.

**THEOREM 2.3:** $L_{R\times R} = L_R \times L_R$ has zero divisors.

Now define a map $\eta : L_R \times L_R = L_{R\times R} \to R \times R$ by

$$\eta(L_a, L_b) = \left(\frac{a}{m+1}, \frac{b}{m+1}\right).$$

Clearly $\eta$ is a ring isomorphism.

Now we define the notion of the space of three dimensional labels or three dimensional labels.

**DEFINITION 2.2:** Let $L_{R\times R\times R} = L_R \times L_R \times L_R = \{(L_a, L_b, L_c) \mid a, b, c \in R; L_a, L_b, L_c \in L_R\}$, we define $L_R \times L_R \times L_R$ as the three dimensional space of refined labels.

Thus $L_R \times L_R \times L_R = \{(L_a, L_b, L_c) = \left(\frac{a}{m+1}, \frac{b}{m+1}, \frac{c}{m+1}\right) \mid$ a, b, c $\in$ R $\}$.

We leave the proof of the following theorems to the reader.

**THEOREM 2.4:** $L_{R\times R\times R} = L_R \times L_R \times L_R$ is isomorphic with the three dimensional real space.

**THEOREM 2.5:** $L_{R\times R\times R}$ is an abelian group under addition.

**THEOREM 2.6:** $L_{R\times R\times R} = L_R \times L_R \times L_R$ is a monoid under multiplication which is commutative with $(L_{m+1}, L_{m+1}, L_{m+1})$ as the identity element under multiplication.

**THEOREM 2.7:** $L_{R \times R \times R} = L_R \times L_R \times L_R$ *is a commutative ring with unit.*

Now as in case of the real plane we can in case of refined plane of labels define distance between two pairs of labels.

For instance P $(L_a, L_b)$ and Q$(L_c, L_d)$ are a given pair of labels in the refined plane.
To find $PQ^2$,

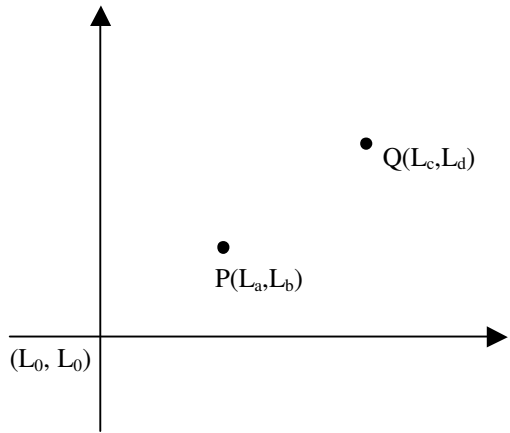

$$\begin{aligned} PQ^2 &= (L_a - L_c)^2 + (L_b - L_d)^2 \\ &= (L_{a-c})^2 + (L_{b-d})^2 \\ &= L_t^2 + L_s^2 \text{ (where } a - c = t \text{ and } b - d = s) \\ &= L_{t^2/m+1} + L_{s^2/m+1} \\ &= L_{t^2+s^2/m+1}. \end{aligned}$$

Thus PQ = $\sqrt{L_{t^2+s^2/m+1}}$. Refined Label plane enables one to plot the labels on them just like the real plane. The analogous result is also true in case of refined label three dimensional space.

We can thus define a n-dimensional refined label space $V = \underbrace{L_R \times L_R \times ... \times L_R}_{n-\text{times}} = \{(L_{a_1}, L_{a_2}, ..., L_{a_n}) / L_{a_i} \in L_R; 1 \leq i \leq n\}$.

Now we proceed onto define refined label vector spaces over the real field $F = R$ the field of reals.

It is pertinent to mention here that $L_{m+1}$ acts as the identity element of $L_R$ with respect to multiplication.

Further $L_a \times L_{m+1} = L_{m+1} \times L_a$

$= L_a$ (as $L_a = \dfrac{a}{m+1}$

$L_{m+1} = \dfrac{m+1}{m+1}$, $L_a \times L_{m+1} = \dfrac{a}{m+1} \times \dfrac{m+1}{m+1} = \dfrac{a}{m+1}$).

Thus we see the linear algebra of refined labels $(L_R, +, \times, .)$ is of infinite dimension over R.

Now we can define by taking $L_R \times L_R = \{(L_a, L_b) \mid L_a, L_b$ belongs to $L_R\}$ to be again a linear algebra of refined labels over R. Here
$(L_a, L_b)(L_c, L_d) = (L_a L_c, L_b L_d)$
$= (L_{ac/m+1}, L_{bd/m+1})$.

Further for $\alpha \in R$ and $(L_a, L_b) \in L_R \times L_R$ we have a

$(L_a, L_b) = (\alpha L_a, \alpha L_b)$
$= (L_{\alpha a}, L_{\alpha b})$ since $\dfrac{\alpha.a}{m+1} = \dfrac{\alpha a}{m+1}$ and $\dfrac{\alpha.b}{m+1} = \dfrac{\alpha b}{m+1}$.

It is easily verified $L_R \times L_R$ is again a linear algebra of refined labels or we will be calling $L_R \times L_R$ as the DSm ring of refined labels. Likewise $L_R \times L_R \times L_R$ is again the DSm ring of refined labels. We will also call $L_R \times L_R$ as the DSm row matrix ring of refined labels. But we will more generally call them as DSm ring (or DSm commutative ring) of refined labels. So ($L_R \times L_R$, +, ×, .) where '.' means scalar product will be known as the DSm ring and linear algebra of refined labels. Thus we can have ($L_R \times L_R \times \ldots \times L_R$, +, ×, .) algebra of refined labels.

Thus $V = \{(L_R \times L_R \times \ldots \times L_R) = \{(L_{a_1}, L_{a_2}, \ldots, L_{a_n})/L_{a_i}\} \in L_R; 1 \leq i \leq n\}$, +, ×, .} is DSm ring and Linear algebra of refined labels.

Now we proceed onto define the notion of DSm ring of matrices.

Consider $T = \left\{ \begin{bmatrix} L_a & L_b \\ L_c & L_d \end{bmatrix} \middle| L_a, L_b, L_c, L_d \in L_R \right\}$ be the collection of all 2 × 2 labels from the DSm field of refined labels. T is a group under matrix addition. Infact T is an abelian group under matrix addition.

For take A and B in T where

$A = \begin{bmatrix} L_a & L_b \\ L_c & L_d \end{bmatrix}$ and $B = \begin{bmatrix} L_t & L_r \\ L_s & L_k \end{bmatrix}$.

Now $A + B = \begin{bmatrix} L_a & L_b \\ L_c & L_d \end{bmatrix} + \begin{bmatrix} L_t & L_r \\ L_s & L_k \end{bmatrix}$

$= \begin{bmatrix} L_a + L_t & L_b + L_r \\ L_c + L_s & L_d + L_k \end{bmatrix} = \begin{bmatrix} L_{a+t} & L_{b+r} \\ L_{c+s} & L_{d+k} \end{bmatrix}$.

Clearly A + B is in T.

Further $L_0 = 0$ acts as the additive identity in $L_R$. So in

$$T(0) = \begin{bmatrix} L_0 & L_0 \\ L_0 & L_0 \end{bmatrix} = \begin{bmatrix} 0 & 0 \\ 0 & 0 \end{bmatrix}$$ acts as the additive identity.

Clearly it can be seen $A+B = B + A$ in T.

Now we proceed onto show how $A \times B$ in T is defined, A and B are given above

$$A \times B = \begin{bmatrix} L_a & L_b \\ L_c & L_d \end{bmatrix} \times \begin{bmatrix} L_t & L_r \\ L_s & L_k \end{bmatrix}$$

$$= \begin{bmatrix} L_a L_t + L_b L_s & L_a L_r + L_b L_k \\ L_c L_t + L_s L_d & L_c L_r + L_d L_k \end{bmatrix} = \begin{bmatrix} L_{\frac{at+bs}{m+1}} & L_{\frac{ar+bk}{m+1}} \\ L_{\frac{ct+sd}{m+1}} & L_{\frac{cr+dk}{m+1}} \end{bmatrix}$$

is in T; as basically the elements $L_a$, $L_b$, $L_c$, ..., $L_s$, $L_k$ are from $L_R$.

Now $L_R$ is only a semigroup of refined labels. For we if

$$A = \begin{bmatrix} 0 & 0 \\ 0 & L_a \end{bmatrix} \text{ and } B = \begin{bmatrix} L_b & 0 \\ 0 & 0 \end{bmatrix}$$

in T then

$$AB = \begin{bmatrix} 0 & 0 \\ 0 & 0 \end{bmatrix} \in T. \text{ Also } BA = \begin{bmatrix} 0 & 0 \\ 0 & 0 \end{bmatrix} \in T,$$

$$\left( \text{since } \begin{bmatrix} 0 & 0 \\ 0 & 0 \end{bmatrix} = \begin{bmatrix} L_0 & L_0 \\ L_0 & L_0 \end{bmatrix} \right).$$

Since every element in T need not necessarily be invertible we see T is only a DSm semigroup of matrix refined labels. Also T has zero divisors hence T is not a group only a semigroup. It is a matter of routine to check that T satisfies the distributive laws, hence T is defined as the DSm ring of refined labels.

Consider W = {all n × n matrices of labels with labels from $L_R$}; W is easily verified to be a DSm ring of refined labels.

Now using these DSm ring of matrices of labels we can define the notion of DSm Linear algebra of refined labels.

**DEFINITION 2.3:** *Let V = {All n × n matrices of labels from $L_R$}; V is defined as the a DSm linear algebra of refined labels over the reals R.*

We will illustrate this situation by some examples.

***Example 2.1:*** Let $V = \left\{ \begin{pmatrix} L_a & L_b & L_c \\ L_d & L_e & L_f \\ L_g & L_h & L_n \end{pmatrix} \right.$ where $L_a$, $L_b$, $L_c$, $L_d$, $L_e$, …, $L_n$ are in $L_R$} be the collection of all 3 × 3 label matrices.

V is a DSm linear algebra of matrix refined labels.

Now we proceed onto define the notion of DSm vector space of refined labels.

For this we have to consider an additive abelian group of refined labels from $L_R$ on which multiplication cannot be defined. In such cases we have only DSm vector space of refined labels and not a DSm linear algebra of refined labels which are not DSm linear algebra of refined labels over the reals R. Consider $M = \left\{ \begin{bmatrix} L_a \\ L_b \end{bmatrix} \middle| L_a, L_b \in L_R \right\}$, M under matrix addition is an abelian group with $\begin{bmatrix} 0 \\ 0 \end{bmatrix} = \begin{bmatrix} L_0 \\ L_0 \end{bmatrix}$ as the additive identity. Clearly for $\begin{bmatrix} L_a \\ L_b \end{bmatrix} \in M$ we cannot define any form of product in

M. Now consider M; M is a vector space over the reals R so M is a DSm vector space of refined labels. We can have infinite number of DSm vector space of refined labels.

**DEFINITION 2.4:** *Let*

$$V = \left\{ \begin{bmatrix} L_{a_1} \\ \vdots \\ L_{a_n} \end{bmatrix} \middle| L_{a_i} \in L_R ; 1 \le i \le n \right\}$$

*be the collection of n × 1 row matrices with entries from $L_R$. V is a vector space over the reals R for*

*(1) V is an abelian group under addition.*

*(2) For $\alpha \in R$ and $s = \begin{bmatrix} L_{a_1} \\ L_{a_2} \\ \vdots \\ L_{a_n} \end{bmatrix}$ in V we have*

$$\alpha s = \alpha \begin{bmatrix} L_{a_1} \\ L_{a_2} \\ \vdots \\ L_{a_n} \end{bmatrix} = \begin{bmatrix} \alpha L_{a_1} \\ \alpha L_{a_2} \\ \vdots \\ \alpha L_{a_n} \end{bmatrix} = \begin{bmatrix} L_{\alpha a_1} \\ L_{\alpha a_2} \\ \vdots \\ L_{\alpha a_n} \end{bmatrix}.$$

*Thus V is defined as the DSm vector space of refined labels over the reals R.*

We will illustrate this by some examples.

***Example 2.2:*** Let us consider

$$W = \left\{ \begin{bmatrix} L_{a_1} \\ L_{a_2} \\ L_{a_3} \\ L_{a_4} \\ L_{a_5} \end{bmatrix} \middle| L_{a_i} \in L_R; 1 \le i \le 5 \right\};$$

W is an abelian group under addition. For any $\alpha \in R$ and

$$T = \begin{bmatrix} L_{a_1} \\ L_{a_2} \\ L_{a_3} \\ L_{a_4} \\ L_{a_5} \end{bmatrix} \text{ in W we have } \alpha T = \begin{bmatrix} L_{\alpha a_1} \\ L_{\alpha a_2} \\ L_{\alpha a_3} \\ L_{\alpha a_4} \\ L_{\alpha a_5} \end{bmatrix} \text{ is in W. Thus W is a}$$

DSm vector space of refined labels over R.

Likewise we can give examples of several DSm vector spaces of refined labels which are not DSm linear algebras of refined labels over R.

Consider $V = \{ \left( L_{a_{ij}} \right) | 1 \le i \le n$ and $1 \le j \le m$; where $L_{a_{ij}} \in L_R\}$, V is only a group under matrix label addition. We call V the collection of $n \times m$ matrix refined labels. V is not a group under multiplication for on V, multiplication cannot be defined. For every $B = \left( L_{b_{ij}} \right)_{n \times m}$ where $1 \le i \le n$ and $1 \le j \le m$ with $L_{b_{ij}} \in L_R$ and $\alpha \in R$ we can define $\alpha B = (\alpha L_{b_{ij}}) = \left( L_{\alpha b_{ij}} \right)$ and $\alpha B$ is in V and $\left( L_{\alpha b_{ij}} \right) \in L_R$. Thus V is a DSm vector space of $n \times m$ matrix of refined labels over the reals R.

We can give some examples of these structures.

*Example 2.3:* Let

$$V = \left\{ \begin{bmatrix} L_a & L_b & L_c & L_j \\ L_d & L_e & L_f & L_k \\ L_g & L_h & L_i & L_n \end{bmatrix} \right.$$

where $L_a, L_b, L_c, L_j, \ldots, L_i, L_n$ are in $L_R$} be the collection of 3 × 4 refined labels of matrices. V is a DSm vector space of refined label matrices over the field, R of reals.

*Example 2.4:* Let

$$P = \left\{ \begin{bmatrix} L_a & L_b \\ L_c & L_d \\ L_f & L_g \\ L_n & L_t \\ L_p & L_q \end{bmatrix} \right.$$

where $L_a, L_b, L_c, L_d, \ldots, L_n, L_t, L_p$ and $L_q$ are in $L_R$}, P is a DSm vector space of refined labels over the field of reals R. Clearly these are not DSm linear algebras they are only DSm vector spaces. Now having seen examples of DSm vector spaces which are not DSm linear algebras over the real field R. But obviously DSm linear algebras are always DSm vector spaces over R.

*Example 2.5:* Let $P = \left\{ \left( L_{a_{ij}} \right)_{5 \times 4} \right.$ refined labels with $L_{a_{ij}} \in L_R$; 1 ≤ i ≤ 5 and 1 ≤ j ≤ 4} be a group of refined labels under addition. P is a DSm vector space refined labels over the reals R. Clearly P is not a DSm linear algebra of refined labels over R. Now we have built DSm matrix of refined labels which are groups under addition and some of them are not even semigroups under multiplication.

We now define DSm polynomials of refined labels.

Let $L_R$ be the DSm field of refined labels. Let x be an indeterminate or variable. Consider $p_l(x) = \sum_{i=0}^{\infty} L_{a_i} x^i$ where $L_{a_i} \in L_R$; we call $p_l(x)$ the polynomial refined labels with coefficients which are refined labels from $L_R$.

We can define addition of two polynomial refined labels.

We denote by $L_R[x] = \left\{ \sum_{i=0}^{\infty} L_{a_i} x^i \,\middle|\, L_{a_i} \in L_R \right\}$, x an indeterminate. $L_R[x]$ is defined as the polynomials with refined label coefficients.

**THEOREM 2.8:** *$L_R[x]$ is a DSm ring of polynomials.*

We define $L_R[x]$ as the polynomial ring with refined label coefficients.

For if $p_l(x) = \sum_{i=0} L_{a_i} x^i$ and $g_l(x) = \sum L_{a_i} x^i$ are in $L_R[x]$ then $p_l(x) + g_l(x)$

$= \sum L_{a_i} x^i + \sum L_{b_i} x^i$

$= \left( L_{a_0} + L_{b_0} \right) + \left( L_{a_1} + L_{b_1} \right) x + \ldots + \left( L_{a_i} + L_{b_i} \right) x^i + \ldots$

$= L_{a_0+b_0} + L_{a_1+b_1} x + \ldots + L_{a_i+b_i} x^i + \ldots$ is in $L_R[x]$.

Thus $0 = 0 + 0x + \ldots + 0x^i + \ldots = L_0 + L_0 x + \ldots + L_0 x^i + \ldots$ as $L_0 = 0$ is the zero polynomial in $L_R[x]$. We can easily verify $L_R[x]$ is an abelian group under addition known as the DSm polynomial group of refined labels.

Further $L_R[x]$ is closed with respect to polynomial multiplication. So $L_R[x]$ can only be a semigroup under multiplication.

It is easily verified that multiplication distributes over addition, hence $(L_R[x], +, \times)$ is a ring defined as the DSm polynomial ring with refined label coefficients. Clearly $L_R[x]$ is a commutative ring and $L_R[x]$ has no zero divisors so we call define $L_R[x]$ as the DSm integral domain of refined labels.

Further $L_R \subseteq L_R[x]$ so $L_R[x]$ is a Smarandache DSm ring.

Now using $L_R[x]$ we can construct linear algebra which we define as DSm linear algebra of polynomials with refined labels.

We just recall this situation in which follows.
$$L_R[x] = \left\{ \sum_{i=0}^{\infty} L_{a_i} x^i \,\middle|\, L_{a_i} \in L_R \right\}$$
is a DSm ring defined as the DSm polynomial ring or integral domain of refined labels.

Now consider $L_R[x]$ the set of polynomials with coefficients from $L_R$. $L_R[x]$ is an additive abelian group and hence $L_R[x]$ is a DSm vector space of refined labels over R or infact DSm is a linear algebra of refined labels over R.

We will illustrate how the product is made. Consider
$$p_1(x) = L_{a_1} + L_{a_2} x^2 + L_{a_3} x^4 + L_{a_4} x^7 + L_{a_5} x^9$$
and
$$q_1(x) = L_{b_1} + L_{b_2} x^5 + L_{b_3} x^{10}$$
in $L_R[x]$.
$$p_1(x) \cdot q_1(x) = L_{a_1} + L_{a_2} x^2 + L_{a_3} x^4 + L_{a_4} x^7 + L_{a_5} x^9$$
$$\times L_{b_1} + L_{b_2} x^5 + L_{b_3} x^{10}$$

$$= L_{a_1b_1/m+1} + L_{a_2b_1/m+1}x^2 + L_{a_3b_1/m+1}x^4 + L_{a_4b_1/m+1}x^7 + L_{a_5b_1/m+1}x^9$$
$$+ L_{a_1b_2/m+1}x^5 + L_{a_2b_2/m+1}x^7 + L_{a_3b_2/m+1}x^9 + L_{a_4b_2/m+1}x^{12} + L_{a_5b_2/m+1}x^{14}$$
$$+ L_{a_1b_3/m+1}x^{10} + L_{a_2b_3/m+1}x^{12} + L_{a_3b_3/m+1}x^{14} + L_{a_4b_3/m+1}x^{17} + L_{a_5b_3/m+1}x^{19}$$

$$= L_{a_1b_1/m+1} + L_{a_2b_1/m+1}x^2 + L_{a_3b_1/m+1}x^4 + L_{a_1b_2/m+1}x^5 + (L_{a_4b_1/m+1}$$
$$+ L_{a_2b_2/m+1})x^7 + (L_{a_5b_1/m+1} + L_{a_3b_2/m+1})x^9 + L_{a_1b_3/m+1}x^{10} + (L_{a_4b_2/m+1}$$
$$+ L_{a_2b_3/m+1})x^{12} + (L_{a_5b_2/m+1} + L_{a_3b_3/m+1})x^{14} + L_{a_4b_3/m+1}x^{17}$$
$$+ L_{a_5b_3/m+1}x^{19} \in L_R[x].$$

Of course addition is also simple using the fact

$$\sum L_{a_i}x^i + \sum L_{b_i}x^i = \sum (L_{a_i} + L_{b_i})x^i$$
$$= \sum L_{a_i+b_i}x^i.$$

Now $\sum L_{a_i}x^i \times \sum L_{b_j}x^j = \sum L_{a_ib_j}x^{i+j}$ and so on.

Thus $L_R[x]$ is a DSm linear algebra over the reals R.

Now having seen examples of DSm linear algebra and vector spaces we proceed onto study properties in them.

**DEFINITION 2.5:** *Let V be a DSm linear algebra of refined labels over the real field R. Let $W \subseteq V$; (W a proper subset of V) if W itself is a DSm linear algebra over the real field R then we define W to be a DSm linear subalgebra of refined labels of V over the real field R.*

We will illustrate this situation by some examples.

***Example 2.6:*** Let $V = \{L_R \times L_R \times L_R = (L_a, L_b, L_c) \mid L_a, L_b, L_c \in L_R, a, b, c \in R\}$ be a DSm linear algebra refined labels over the real field R. Consider $W = \{(L_a, 0, L_b) \mid L_a, L_b \in L_R\} \subseteq V$. W is a DSm linear subalgebra of refined labels over the field R of V.

***Example 2.7:*** Let

$$V = \left\{ \begin{bmatrix} L_{a_1} \\ L_{a_2} \\ L_{a_3} \\ L_{a_4} \end{bmatrix} \middle| L_{a_1}, L_{a_2}, L_{a_3}, L_{a_4} \in L_R \right\}$$

be a DSm vector space of refined labels over the field R.
Consider

$$W = \left\{ \begin{bmatrix} 0 \\ L_a \\ L_b \\ 0 \end{bmatrix} \middle| L_a, L_b \in L_R \right\} \subseteq V;$$

it is easily verified W is a DSm vector subspace of refined labels over the field R of V.

***Example 2.8:*** Let

$$V = \left\{ \begin{bmatrix} L_a & L_b & L_c \\ L_d & L_f & L_h \\ L_d & L_e & L_n \end{bmatrix} \right.$$

where $L_a, L_b, L_c, L_d, L_f, L_h, L_d, L_e, L_n$ are in $L_R\}$ be a DSm linear algebra of refined labels over the real field R.
Consider

$$T = \left\{ \begin{bmatrix} L_a & 0 & 0 \\ 0 & L_b & 0 \\ 0 & 0 & L_c \end{bmatrix} \middle| L_a, L_b, L_c \text{ are in } L_R \right\} \subseteq V;$$

T is clearly a DSm linear subalgebra refined labels of V over the real field R.

*Example 2.9*: Let

$$M = \left\{ \begin{bmatrix} L_a & L_b & L_c \\ L_{a_1} & L_{b_1} & L_{c_1} \\ L_{a_2} & L_{b_2} & L_{c_2} \\ L_{a_3} & L_{b_3} & L_{c_3} \\ L_{a_4} & L_{b_4} & L_{c_4} \end{bmatrix} \right.$$

where $L_a$, $L_b$, $L_c$, $L_{a_i}, L_{b_i}, L_{c_i}$ are in $L_R$; $1 \leq i \leq 4$} be a DSm vector space of refined labels over the real field R.

Consider

$$N = \left\{ \begin{bmatrix} L_a & L_b & L_c \\ 0 & 0 & 0 \\ L_{a_1} & L_{a_2} & L_{a_3} \\ 0 & 0 & 0 \\ 0 & 0 & 0 \end{bmatrix} \middle| L_a, L_b, L_c, L_{a_i} \in L_R; 1 \leq i \leq 3 \right\}$$

be a DSm vector subspace of refined labels of M over the real field R.

*Example 2.10:* Let

$$P = \left\{ \begin{pmatrix} L_{a_1} & \cdots & L_{a_8} \\ L_{b_1} & \cdots & L_{b_8} \end{pmatrix} \middle| L_{a_i}, L_{b_i} \in L_R; 1 \leq i \leq 8 \right\}$$

be a DSm vector space of refined labels over the real field R.

Take

$$M = \left\{ \begin{pmatrix} L_{a_1} & 0 & L_{a_2} & 0 & L_{a_3} & 0 & 0 & L_{a_4} \\ 0 & L_{b_1} & 0 & L_{b_2} & 0 & L_{b_3} & 0 & 0 \end{pmatrix} \middle| \begin{array}{l} L_{a_i}, L_{b_j} \in L_R; \\ 1 \leq i \leq 4 \, \& \, 1 \leq j \leq 3 \end{array} \right\}$$

$\subseteq$ P is a DSm vector subspace of refined labels of P over the real field R.

***Example 2.11***: Let $X = \left\{ \begin{bmatrix} L_a & L_b \\ L_c & L_d \end{bmatrix} \middle| L_a, L_b, L_c, L_d \in L_R \right\}$ be a DSm linear algebra of refined labels over the real field R.

Consider $P = \left\{ \begin{bmatrix} L_a & L_b \\ 0 & 0 \end{bmatrix} \middle| L_a, L_b \in L_R \right\} \subseteq X$, P is a DSm linear subalgebra of refined labels over the real field R of X.

***Example 2.12:*** Let $Y = \left\{ \begin{pmatrix} L_a & L_b & L_c \\ L_d & L_e & L_f \\ L_g & L_h & L_n \end{pmatrix} \right.$ where $L_a$, $L_b$, $L_c$, $L_d$, $L_e$, $L_f$, $L_g$, $L_h$, $L_n \in L_R \}$ be a DSm linear algebra of refined labels over the real field R.

Consider $M = \left\{ \begin{pmatrix} L_a & L_b & L_c \\ 0 & L_d & 0 \\ L_e & 0 & L_p \end{pmatrix} \middle| L_a, L_b, L_c, L_d, L_e, L_p \in L_R \right\}$

$\subseteq$ Y, we see M is not multiplicatively closed for if A, B are in M then AB $\notin$ M. Thus M $\subseteq$ Y is only a DSm vector subspace of refined labels of Y over the reals R. We call such DSm vector subspaces of refined labels of a DSm linear algebra of refined labels to be a pseudo DSm vector subspace of refined labels of Y over the real field R.

Now we want to define yet another concept. Suppose B be a DSm linear algebra of refined labels over the real field R. We see Q $\subseteq$ R is a subfield of R. Let V $\subseteq$ B; if V is a DSm linear algebra of refined labels over the field Q then we define V to be a subfield DSm linear subalgebra of refined labels of B over the subfield Q contained in R.

Interested reader can construct examples of these structures.

We now proceed onto define direct sum of DSm subspaces of refined labels and just sum of DSm subspaces of refined labels of a DSm vector space of refined labels.

**DEFINITION 2.6:** *Let V be a DSm vector space of refined labels over the real field R. Let $W_1, W_2, \ldots, W_n$ be n DSm vector subspaces of refined labels of V over the real field R. If $V = W_1 + \ldots + W_n$ and $W_i \cap W_j = (0) = (L_0)$; if $i \neq j$; $1 \leq i, j \leq n$ then we define $W_1, W_2, \ldots, W_n$ to be the DSm direct sum of vector subspaces refined labels of V over the field F.*

*If $V = W_1 + \ldots + W_n$ with $W_i \cap W_j \neq (0) \neq (L_0)$ for some $i \neq j$, $1 \leq i, j \leq n$ then we define V to be the sum of the vector subspaces of refined labels of V over the field F.*

We will illustrate this situation by some examples.

***Example 2.13:*** Let $V = \{(L_{a_1}, L_{a_2}, \ldots, L_{a_9}) \mid L_{a_i} \in L_R; 1 \leq i \leq 9\}$ be a DSm linear algebra of refined labels over the real field R.
Consider
$$W_1 = \{(L_{a_1}, 0, L_{a_2}, 0, \ldots, 0) \mid L_{a_1}, L_{a_2} \in L_R\},$$
$$W_2 = \{(0,0,0,0,L_{a_3}, L_{a_4}, 0, 0, L_{a_5}) \mid L_{a_3}, L_{a_4}, L_{a_5} \in L_R\},$$
$$W_3 = \{(0, L_{a_6}, 0,0,0,0, L_{a_7}, L_{a_8}, 0) \mid L_{a_6}, L_{a_7}, L_{a_8} \in L_R\}$$
and
$$W_4 = \{(0,0,0,L_{a_9},0,0,0,0,0) \mid L_{a_9} \in L_R\}$$
be four DSm linear subalgebras of V of refined labels over the field R.

Clearly $V = W_1 + W_2 + W_3 + W_4$ and $W_i \cap W_j = (0, 0, \ldots, 0)$ if $i \neq j$; $1 \leq i, j \leq 4$.

Thus V is the direct sum of DSm linear subalgebras of refined labels over the field R.

Now it may so happen that $W_1, \ldots, W_n$ may not be DSm linear subalgebras but only a DSm vector subspaces of refined labels over the field R in which case we write V as a pseudo direct sum of DSm vector subspaces of refined labels of the DSm linear algebra of refined labels over the real field R.

*Example 2.14:* Let

$$V = \left\{ \begin{pmatrix} L_{a_1} & L_{a_3} & L_{a_5} & L_{a_7} \\ L_{a_2} & L_{a_4} & L_{a_6} & L_{a_8} \end{pmatrix} \middle| L_{a_i} \in L_R; 1 \leq i \leq 8 \right\}$$

be a DSm vector space of refined labels over the real field R.

Consider

$$W_1 = \left\{ \begin{pmatrix} 0 & 0 & L_{a_2} & 0 \\ L_{a_1} & 0 & 0 & L_{a_3} \end{pmatrix} \middle| L_{a_1}, L_{a_2}, L_{a_3} \in L_R \right\},$$

$$W_2 = \left\{ \begin{pmatrix} L_{a_4} & 0 & 0 & L_{a_5} \\ 0 & 0 & 0 & 0 \end{pmatrix} \middle| L_{a_4}, L_{a_5} \in L_R \right\},$$

$$W_3 = \left\{ \begin{pmatrix} 0 & L_{a_6} & 0 & 0 \\ 0 & L_{a_7} & 0 & 0 \end{pmatrix} \middle| L_{a_6}, L_{a_7} \in L_R \right\}$$

and $W_4 = \left\{ \begin{pmatrix} 0 & 0 & 0 & 0 \\ 0 & 0 & L_{a_8} & 0 \end{pmatrix} \middle| L_{a_8} \in L_R \right\}$

be DSm vector subspaces of refined labels of V over the real field R.

Clearly $W_i \cap W_j = \begin{pmatrix} 0 & 0 & 0 & 0 \\ 0 & 0 & 0 & 0 \end{pmatrix}$; if $i \neq j$; $1 \leq i, j \leq 4$.

Further $V = W_1 + W_2 + W_3 + W_4$. Thus V is a direct sum of DSm vector subspaces of refined labels of V over the field R.

*Example 2.15:* Let

$$V = \left\{ \begin{bmatrix} L_{a_1} & L_{a_2} & L_{a_3} \\ L_{a_4} & L_{a_5} & L_{a_6} \\ L_{a_7} & L_{a_8} & L_{a_9} \end{bmatrix} \middle| L_{a_i} \in L_R; 1 \leq i \leq 9 \right\}$$

be a DSm linear algebra of refined labels over the field R.

Consider

$$W_1 = \left\{ \begin{bmatrix} L_{a_1} & 0 & L_{a_2} \\ L_{a_3} & 0 & L_{a_4} \\ L_{a_5} & L_{a_6} & 0 \end{bmatrix} \,\middle|\, L_{a_i} \in L_R; 1 \leq i \leq 6 \right\},$$

$$W_2 = \left\{ \begin{bmatrix} 0 & L_{a_7} & 0 \\ 0 & 0 & 0 \\ 0 & 0 & 0 \end{bmatrix} \,\middle|\, L_{a_7} \in L_R \right\}$$

and

$$W_3 = \left\{ \begin{bmatrix} 0 & 0 & 0 \\ 0 & L_{a_8} & 0 \\ 0 & 0 & L_{a_9} \end{bmatrix} \,\middle|\, L_{a_8}, L_{a_9} \in L_R \right\}$$

be DSm vector subspaces of refined labels of V over the field R.

Clearly

$$W_i \cap W_j = \begin{pmatrix} 0 & 0 & 0 \\ 0 & 0 & 0 \\ 0 & 0 & 0 \end{pmatrix} \text{ if } i \neq j; \, 1 \leq i, j \leq 3.$$

Further $V = W_1 + W_2 + W_3$, hence V is the pseudo direct sum of DSm vector subspace of refined labels of V, the DSm linear algebra of refined labels over the real field R.

*Example 2.16:* Let

$$V = \left\{ \sum_{i=0}^{\infty} L_{a_i} x^i \,\middle|\, L_{a_i} \in L_R \right\}$$

be a DSm linear algebra of refined labels over the real field R.

Consider $\quad W_1 = \left\{ \sum_{i=0}^{2n} L_{a_i} x^i \,\middle|\, L_{a_i} \in L_R \right\} \subseteq V,$

is a DSm vector pseudo subspace of V of refined labels over the field R.

$$W_2 = \left\{ \sum_{i=2n+1}^{120n} L_{a_i} x^i \,\middle|\, L_{a_i} \in L_R \right\} \subseteq V$$

is again a pseudo DSm vector subspace of V over the field R.

However we cannot by this way define $W_i$ and write $V = \sum_i W_i$. But however if we take

$$V_1 = \left\{ \sum_{i=0}^{120} L_{a_i} x^i \;\middle|\; L_{a_i} \in L_R \right\}, \quad V_2 = \left\{ \sum_{i=121}^{500} L_{a_i} x^i \right\}$$

and 
$$V_3 = \left\{ \sum_{i=501}^{\infty} L_{a_i} x^i \;\middle|\; L_{a_i} \in L_R \right\}$$

as pseudo DSm vector subspaces of V of refined labels over the field R then we see $V_i \cap V_j = 0$ if $i \ne j$, $1 \le i, j \le 3$ and $V = V_1 + V_2 + V_3$ is a direct sum of pseudo DSm vector subspace of V of refined labels over the field R. Now having see examples of direct sum now we proceed onto give examples of sum of DSm linear subalgebras and DSm subvector subspaces over a field R.

*Example 2.17:* Let

$$V = \left\{ \begin{bmatrix} L_{a_1} & L_{a_4} \\ L_{a_2} & L_{a_5} \\ L_{a_3} & L_{a_6} \end{bmatrix} \;\middle|\; L_{a_i} \in L_R ; 1 \le i \le 6 \right\}$$

be a DSm vector space of refined labels over the field $F = R$.

Consider

$$W_1 = \left\{ \begin{bmatrix} L_{a_1} & L_{a_4} \\ 0 & L_{a_3} \\ 0 & 0 \end{bmatrix} \;\middle|\; L_{a_i} \in L_R ; 1 \le i \le 3 \right\} \subseteq V,$$

$$W_2 = \left\{ \begin{bmatrix} L_{a_4} & 0 \\ L_{a_5} & 0 \\ 0 & 0 \end{bmatrix} \;\middle|\; L_{a_4}, L_{a_5} \in L_R \right\} \subseteq V,$$

$$W_3 = \left\{ \begin{bmatrix} 0 & L_{a_1} \\ 0 & 0 \\ L_{a_2} & L_{a_3} \end{bmatrix} \;\middle|\; L_{a_i} \in L_R ; 1 \le i \le 3 \right\} \subseteq V,$$

$$W_4 = \left\{ \begin{bmatrix} L_{a_1} & L_{a_2} \\ 0 & 0 \\ L_{a_3} & L_{a_4} \end{bmatrix} \middle| L_{a_i} \in L_R; 1 \le i \le 4 \right\} \subseteq V$$

and $$W_5 = \left\{ \begin{bmatrix} 0 & 0 \\ L_{a_1} & L_{a_2} \\ L_{a_3} & 0 \end{bmatrix} \middle| L_{a_i} \in L_R; 1 \le i \le 3 \right\} \subseteq V;$$

be a DSm vector subspaces of V of refined labels over the field R. Clearly $V_i \cap V_j \ne \begin{bmatrix} 0 & 0 \\ 0 & 0 \\ 0 & 0 \end{bmatrix}$ if $i \ne j$; $1 \le i, j \le 5$ however $V = W_1 + W_2 + W_3 + W_4 + W_5$. Thus V is the sum of DSm vector subspace of refined labels of V over the field R.

*Example 2.18*: Let

$$V = \left\{ \begin{bmatrix} L_{a_1} & L_{a_2} & L_{a_3} & L_{a_4} \\ L_{b_1} & L_{b_2} & L_{b_3} & L_{b_4} \\ L_{c_1} & L_{c_2} & L_{c_3} & L_{c_4} \end{bmatrix} \middle| L_{a_i}, L_{b_i}, L_{c_i} \in L_R; 1 \le i \le 4 \right\}$$

be a DSm vector space of refined labels over the field R.
Consider

$$W_1 = \left\{ \begin{bmatrix} L_{a_1} & 0 & 0 & 0 \\ L_{a_2} & 0 & 0 & 0 \\ L_{a_3} & 0 & 0 & 0 \end{bmatrix} \middle| L_{a_i} \in L_R; 1 \le i \le 3 \right\}$$

$$W_2 = \left\{ \begin{bmatrix} 0 & L_{b_1} & 0 & 0 \\ 0 & L_{b_2} & 0 & 0 \\ 0 & L_{b_3} & L_{b_4} & 0 \end{bmatrix} \middle| L_{b_i} \in L_R; 1 \le i \le 4 \right\}$$

and $W_3 = \left\{ \begin{bmatrix} 0 & 0 & 0 & L_{a_1} \\ 0 & 0 & 0 & L_{a_2} \\ 0 & 0 & 0 & L_{a_3} \end{bmatrix} \middle| L_{a_i} \in L_R; 1 \leq i \leq 3 \right\}$

be DSm vector subspaces of refined labels of V.

Clearly $W_i \cap W_j = \begin{pmatrix} 0 & 0 & 0 & 0 \\ 0 & 0 & 0 & 0 \\ 0 & 0 & 0 & 0 \end{pmatrix}$ if $i \neq j$ but $V \neq W_1 + W_2 +$ $W_3$ as some elements are missing to be full of V; we see $W_1 + W_2 + W_3 \subseteq V$ (the containment being proper).

*Example 2.19:* Let

$$V = \left\{ \begin{bmatrix} L_{a_1} & L_{a_2} & L_{a_3} \\ L_{a_4} & L_{a_5} & L_{a_6} \\ L_{a_7} & L_{a_8} & L_{a_9} \\ L_{a_{10}} & L_{a_{11}} & L_{a_{12}} \\ L_{a_{13}} & L_{a_{14}} & L_{a_{15}} \\ L_{a_{16}} & L_{a_{17}} & L_{a_{18}} \\ L_{a_{19}} & L_{a_{20}} & L_{a_{21}} \end{bmatrix} \middle| L_{a_i} \in L_R; 1 \leq i \leq 21 \right\}$$

be a DSm vector space of refined labels over the field R.

Consider

$$W_1 = \left\{ \begin{bmatrix} L_{a_1} & 0 & L_{a_7} \\ L_{a_2} & 0 & 0 \\ L_{a_3} & 0 & 0 \\ L_{a_4} & L_{a_8} & 0 \\ L_{a_5} & 0 & 0 \\ L_{a_6} & 0 & L_{a_9} \\ 0 & 0 & L_{a_{10}} \end{bmatrix} \middle| L_{a_i} \in L_R; 1 \leq i \leq 10 \right\},$$

$$W_2 = \left\{ \begin{bmatrix} L_{a_1} & L_{a_3} & L_{a_4} \\ L_{a_2} & 0 & 0 \\ L_{a_5} & 0 & 0 \\ 0 & L_{a_7} & 0 \\ 0 & L_{a_8} & 0 \\ 0 & L_{a_9} & L_{a_{11}} \\ 0 & L_{a_{10}} & 0 \end{bmatrix} \middle| L_{a_i} \in L_R; 1 \leq i \leq 11 \right\}$$

$$W_3 = \left\{ \begin{bmatrix} 0 & L_{a_1} & 0 \\ L_{a_2} & 0 & 0 \\ L_{a_3} & 0 & 0 \\ L_{a_5} & L_{a_7} & 0 \\ 0 & 0 & 0 \\ 0 & 0 & L_{a_9} \\ L_{a_8} & 0 & 0 \end{bmatrix} \middle| L_{a_i} \in L_R; 1 \leq i \leq 9 \right\}$$

and

$$W_4 = \left\{ \begin{bmatrix} L_{a_1} & 0 & 0 \\ 0 & 0 & 0 \\ 0 & 0 & 0 \\ L_{a_2} & L_{a_3} & L_{a_4} \\ 0 & 0 & 0 \\ L_{a_6} & L_{a_7} & L_{a_5} \\ 0 & L_{a_8} & L_{a_9} \end{bmatrix} \middle| L_{a_i} \in L_R; 1 \leq i \leq 9 \right\}$$

be DSm vector subspaces of the refined labels of V over the field R.

$$W_1 \cap W_2 = \begin{bmatrix} L_{a_1} & 0 & L_{a_2} \\ L_{a_3} & 0 & 0 \\ L_{a_4} & 0 & 0 \\ 0 & L_{a_5} & 0 \\ 0 & 0 & 0 \\ 0 & 0 & L_{a_6} \\ 0 & 0 & 0 \end{bmatrix}, W_1 \cap W_3 = \begin{bmatrix} 0 & 0 & 0 \\ L_{a_1} & 0 & 0 \\ L_{a_2} & 0 & 0 \\ L_{a_5} & L_{a_6} & 0 \\ 0 & 0 & 0 \\ 0 & 0 & L_{a_6} \\ 0 & 0 & 0 \end{bmatrix},$$

$$W_1 \cap W_4 = \begin{bmatrix} L_{a_1} & 0 & 0 \\ 0 & 0 & 0 \\ 0 & 0 & 0 \\ L_{a_4} & L_{a_2} & 0 \\ 0 & 0 & 0 \\ L_{a_3} & 0 & L_{a_6} \\ 0 & 0 & L_{a_5} \end{bmatrix}, W_2 \cap W_3 = \begin{bmatrix} 0 & L_{a_1} & 0 \\ L_{a_2} & 0 & 0 \\ L_{a_3} & 0 & 0 \\ 0 & L_{a_4} & 0 \\ 0 & 0 & 0 \\ 0 & 0 & L_{a_5} \\ 0 & 0 & 0 \end{bmatrix},$$

$$W_2 \cap W_4 = \begin{bmatrix} L_{a_1} & 0 & 0 \\ 0 & 0 & 0 \\ 0 & 0 & 0 \\ 0 & L_{a_2} & 0 \\ 0 & 0 & 0 \\ 0 & L_{a_3} & L_{a_4} \\ 0 & L_{a_5} & 0 \end{bmatrix} \text{ and } W_3 \cap W_4 = \begin{bmatrix} 0 & 0 & 0 \\ 0 & 0 & 0 \\ 0 & 0 & 0 \\ L_{a_1} & L_{a_2} & 0 \\ 0 & 0 & 0 \\ 0 & 0 & L_{a_5} \\ 0 & 0 & 0 \end{bmatrix}.$$

Thus $W_i \cap W_j \neq (0)$ for $1 \leq i, j \leq 4$ with $i \neq j$. However $V \neq W_1 + W_2 + W_3 + W_4$ as only $W_1 + W_2 + W_3 + W_4 \subset V$. Certain elements are left out. Now having seen examples of direct sum of DSm vector subspaces of refined labels and the set of DSm vector subspaces which does not give direct sum or sum we now proceed onto discuss about DSm linear transformation of refined labels. Let V and W be any two DSm vector spaces of refined labels over the field of reals R.

If T : V → W is such that T (a+b) = T (a) + T(b) and T(ca) = cT (a) for all a, b ∈ V and c ∈ R we define T to be a DSm linear transformation of refined labels of V to W.

We will illustrate this situation by some examples.

***Example 2.20:*** Let V = { $(L_{a_1}, L_{a_2}, L_{a_3})$ | $L_{a_i}$ ∈ $L_R$, 1 ≤ i ≤ 3} and

$$W = \left\{ \begin{bmatrix} L_{a_1} & L_{a_2} \\ L_{a_4} & L_{a_3} \\ L_{a_5} & L_{a_6} \end{bmatrix} \middle| L_{a_i} \in L_R; 1 \leq i \leq 6 \right\}$$

be two DSm vector spaces of refined labels over the field of reals R. Define η : V → W by

$$\eta((L_a, L_b, L_c)) = \begin{bmatrix} L_a & L_a \\ L_b & L_b \\ L_c & L_c \end{bmatrix}$$

for all $(L_a, L_b, L_c)$ in V. It is easily verified η: V → W is a DSm vector space linear transformation of refined label over the field R.

***Example 2.21:*** Let

$$V = \left\{ \begin{bmatrix} L_{a_1} & L_{a_2} & L_{a_3} & L_{a_4} \\ L_{a_5} & L_{a_6} & L_{a_7} & L_{a_8} \end{bmatrix} \middle| L_{a_i} \in L_R; 1 \leq i \leq 8 \right\}$$

be a DSm vector space of refined labels over the field R.

$$W = \left\{ \begin{bmatrix} L_{a_1} & L_{a_2} & L_{a_3} & L_{a_4} \\ L_{a_5} & L_{a_6} & L_{a_7} & L_{a_8} \\ L_{a_9} & L_{a_{10}} & L_{a_{11}} & L_{a_{12}} \end{bmatrix} \middle| L_{a_i} \in L_R; 1 \leq i \leq 12 \right\}$$

be a DSm vector space of refined labels over the field R.

Define T : V → W by

$$T \begin{bmatrix} L_{a_1} & L_{a_2} & L_{a_3} & L_{a_4} \\ L_{a_5} & L_{a_6} & L_{a_7} & L_{a_8} \end{bmatrix} = \begin{bmatrix} L_{a_1} & 0 & L_{a_4} & 0 \\ L_{a_2} & L_{a_7} & L_{a_8} & L_{a_5} \\ L_{a_3} & 0 & 0 & L_{a_6} \end{bmatrix}.$$

Clearly T : V → W is a DSm vector linear transformation of refined label spaces.

Consider P : W → V defined by

$$P \left( \begin{bmatrix} L_{a_1} & L_{a_2} & L_{a_3} & L_{a_4} \\ L_{a_5} & L_{a_6} & L_{a_7} & L_{a_8} \\ L_{a_9} & L_{a_{10}} & L_{a_{11}} & L_{a_{12}} \end{bmatrix} \right) = \begin{bmatrix} L_{a_1} & L_{a_2} & L_{a_3} & L_{a_4} \\ L_{a_5} & L_{a_6} & L_{a_7} & L_{a_8} \end{bmatrix}.$$

P is a DSm linear transformation of W to V.

If V = W then as in case of usual vector spaces define the DSm linear transformation as DSm linear operator of refined label spaces.

*Example 2.22*: Let

$$M = \left\{ \begin{bmatrix} L_{a_1} & L_{a_2} \\ L_{a_3} & L_{a_4} \end{bmatrix} \middle| L_{a_i} \in L_R; 1 \le i \le 4 \right\}$$

and

$$P = \left\{ \begin{bmatrix} L_{a_1} & L_{a_2} & L_{a_3} \\ L_{a_4} & L_{a_5} & L_{a_6} \\ L_{a_7} & L_{a_8} & L_{a_9} \end{bmatrix} \middle| L_{a_i} \in L_R; 1 \le i \le 9 \right\}$$

be two DSm linear algebras of refined labels over the real field R.

Define T : M → P by

$$T\left(\begin{bmatrix} L_{a_1} & L_{a_2} \\ L_{a_3} & L_{a_4} \end{bmatrix}\right) = \begin{bmatrix} L_{a_1} & L_{a_2} & 0 \\ L_{a_3} & L_{a_4} & 0 \\ 0 & 0 & 0 \end{bmatrix},$$

T is a DSm linear transformation of linear algebra of refined labels over R.

Suppose S : P → M defined by

$$S\left(\begin{bmatrix} L_{a_1} & L_{a_2} & L_{a_3} \\ L_{a_4} & L_{a_5} & L_{a_6} \\ L_{a_7} & L_{a_8} & L_{a_9} \end{bmatrix}\right) = \begin{bmatrix} L_{a_1} & L_{a_2} \\ L_{a_3} & L_{a_4} \end{bmatrix}$$

S is also a DSm linear transformation of refined labels of linear algebras.

*Example 2.23:* Let

$$V = \left\{ \begin{bmatrix} L_{a_1} & L_{a_2} \\ L_{a_3} & L_{a_4} \end{bmatrix} \middle| L_{a_i} \in L_R; 1 \leq i \leq 4 \right\}$$

be a DSm linear algebra of refined labels over F.

Define T : V → V by

$$T\left(\begin{bmatrix} L_{a_1} & L_{a_2} \\ L_{a_3} & L_{a_4} \end{bmatrix}\right) = \begin{bmatrix} L_{a_1} & L_{a_2} \\ 0 & L_{a_3} \end{bmatrix}.$$

It is easily verified T is DSm linear operator on V.

*Example 2.24*: Let

$$W = \left\{ \sum_{i=0}^{\infty} L_{a_i} x^i \middle| L_{a_i} \in L_R \right\}$$

be a DSm linear algebra of refined labels over the real field R. Define T : W → W by

$$T\left(\sum_{i=0}^{\infty} L_{a_i} x^i\right) = \sum_{i=0}^{\infty} L_{a_i} x^{2i}.$$

It is easily verified T is a DSm linear operator on W.

Now having seen DSm linear operator on DSm linear algebra we proceed onto define concepts like DSm projections. Let V be a DSm vector space of refined labels over reals. Suppose $W_1, W_2, \ldots, W_t$ be t, DSm subspaces of refined labels of V such that

$V = W_1 + \ldots + W_t$ is a direct sum.

Define $T_i$ an DSm linear operator on V such that

$$T_i(v) = \begin{cases} \omega_i \in W_i \text{ if } v \in V \\ 0 \text{ if } v \notin W_i \text{ and } v \in V. \end{cases}$$

$T_i$ is also a DSm linear operator of a special type.

Let V be a DSm vector space of refined labels over the real field R. Let E be a DSm linear operator on V where E is a projection such that $E^2 = E$.

Then we have the following interesting properties. Suppose E is a DSm projection. Let R be the DSm range space of E and let N be the DSm null space of E.

Then the refined label of vector β is in R the DSm range space of V if and only if Eβ = β.

If β = Eα then Eβ = Eα² = Eα = β. Conversely if β = Eβ then of course β is in the range of E.

V = R + N, DSm direct sum of DSm subspaces. The unique expression for α as a sum of vectors in R and N is
$$\alpha = E\alpha + (\alpha - E\alpha).$$

Before we proceed to make conclusions of R, N and E we give an illustrative example.

***Example 2.25***: Let $V = (L_R \times L_R \times L_R) = \{(L_a, L_b, L_c) / L_a, L_b, L_c$ are in $L_R\}$ be a DSm vector space of refined labels over the real field R.

Let E be a DSm linear operator on V such that $E : V \to V$ is such that $E(v) = E((L_a, L_b, L_c)) = (L_a, 0, L_c)$ for every $v \in V$.

It is easily verified E is a DSm linear operator on V. The DSm null space of $E = \{(v \in V \mid E(v) = (0, 0, 0)\} = \{(0, L_a, 0) / L_a \in L_R\} = N$. It is easily verified N is a DSm vector subspace of refined labels of V.

Now consider the DSm range space of the DSm linear operator on V. The DSm range space $W = \{(L_a, 0, L_b) / L_a, L_b \in L_R\}$.

It is easily verified that $V = W + N$ and $W \cap N = (0, 0, 0)$, that is V is the DSm direct sum of DSm vector subspaces of V.

Further if $v \in V$ then we see $v = (L_a, 0, L_b) \in W$ now $E_v = E(L_a, 0, L_b) = (L_a, 0, L_b)$. If $w = Ev$, then $Ew = Ev^2 = Ev = w$. Conversely if $v = Ev$ then v is in W the DSm range space of the DSm linear operator of refined labels E. Thus $V = W + N$.
Any $v \in V$ is such that $Ev + (v - Ev)$ for if $v = (L_a, L_b, L_c)$ then $E(v) = E((L_a, L_b, L_c)) = (L_a, 0, L_c)$. Thus $Ev + (v - Ev) = (L_a, 0, L_c) + (L_a, L_b, L_c) - (L_a, 0, L_c) = (L_a, L_b, L_c) = V$. This representation is unique.

Thus we can conclude as in case of vector spaces from the above results that if W and N DSm subspaces of refined labels of V such that $V = W + N$ (as direct sum) then there is one and only one DSm projection operator E which has the DSm range space of refined labels to be W and DSm null space of E of refined labels to be N.

This DSm operator E is defined as the projection on W along N.

Before we discuss about the properties of DSm projection on a DSm vector space of refined labels we proceed onto recall some analogous results on DSm vector space of refined labels.

Suppose $W_1, W_2, \ldots, W_k$ are k DSm subspaces of refined labels of V, we say $W_1, \ldots, W_k$ are independent if
$$\alpha_1 + \ldots + \alpha_k = 0; \alpha_i \in W_i$$

implies $\alpha_i = 0$, $1 \leq i \leq k$.

We will give some illustrative examples of them for this concept can be developed as a matter of routine.

***Example 2.26:*** Let

$$V = \left\{ \begin{pmatrix} L_{a_1} & L_{a_2} & L_{a_3} \\ L_{a_4} & L_{a_5} & L_{a_6} \\ L_{a_7} & L_{a_8} & L_{a_9} \\ L_{a_{10}} & L_{a_{11}} & L_{a_{12}} \end{pmatrix} \right.$$

where $L_{a_i} \in L_R$; $1 \leq i \leq 12\}$ be a DSm vector space of refined labels over the reals R.

Consider

$$W_1 = \left\{ \begin{pmatrix} L_{a_1} & L_{a_2} & L_{a_3} \\ 0 & 0 & 0 \\ 0 & 0 & 0 \\ 0 & 0 & 0 \end{pmatrix} \middle| L_{a_i} \in L_R; 1 \leq i \leq 3 \right\} \subseteq V$$

is a DSm vector subspace of refined labels over R.

$$W_2 = \left\{ \begin{pmatrix} 0 & 0 & 0 \\ L_{a_1} & L_{a_2} & 0 \\ L_{a_3} & L_{a_4} & 0 \\ 0 & 0 & 0 \end{pmatrix} \middle| L_{a_i} \in L_R; 1 \leq i \leq 4 \right\} \subseteq V$$

is a DSm vector subspace of refined labels of V over the real field R.

Let

$$W_3 = \left\{ \begin{pmatrix} 0 & 0 & 0 \\ 0 & 0 & 0 \\ 0 & 0 & 0 \\ L_{a_1} & L_{a_2} & L_{a_3} \end{pmatrix} \middle| L_{a_i} \in L_R; 1 \leq i \leq 3 \right\} \subseteq V$$

be a DSm vector subspace of refined labels of V over the real field R and

$$W_4 = \left\{ \begin{pmatrix} 0 & 0 & 0 \\ 0 & 0 & L_{a_1} \\ 0 & 0 & L_{a_2} \\ 0 & 0 & 0 \end{pmatrix} \middle| L_{a_i} \in L_R; 1 \le i \le 2 \right\} \subseteq V$$

be a DSm vector subspace of V of refined labels over R.

We see $W_1$, $W_2$, $W_3$ and $W_4$ are independent DSm vector subspace of refined labels over R.

For we see $\alpha_1 + \alpha_2 + \alpha_3 + \alpha_4 = \begin{pmatrix} 0 & 0 & 0 \\ 0 & 0 & 0 \\ 0 & 0 & 0 \\ 0 & 0 & 0 \end{pmatrix}$ for $\alpha_i \in W_i$;

$1 \le i \le 4$ is possible only if $\alpha_i = \begin{pmatrix} 0 & 0 & 0 \\ 0 & 0 & 0 \\ 0 & 0 & 0 \end{pmatrix}$ for each $\alpha_i$.

Further $W_i \cap W_j = (0)$; $i \ne j$.

*Example 2.27:* Let

$$V = \left\{ \begin{pmatrix} L_{a_1} & L_{a_2} & L_{a_3} \\ L_{a_4} & L_{a_5} & L_{a_6} \\ L_{a_7} & L_{a_8} & L_{a_9} \end{pmatrix} \middle| L_i \in L_R; 1 \le i \le 9 \right\}$$

be a DSm vector space of refined labels over the real field R.
Choose

$$W_1 = \left\{ \begin{pmatrix} L_1 & 0 & 0 \\ 0 & L_2 & 0 \\ 0 & 0 & L_3 \end{pmatrix} \middle| L_i \in L_R; 1 \le i \le 3 \right\} \subseteq V$$

be a subvector space of refined labels over the real field R.

$$W_2 = \left\{ \begin{pmatrix} 0 & L_2 & L_3 \\ 0 & 0 & L_1 \\ 0 & 0 & 0 \end{pmatrix} \middle| L_i \in L_R; 1 \le i \le 3 \right\} \subseteq V$$

be a DSm vector subspace of V of refined labels over R.

$$W_3 = \left\{ \begin{pmatrix} 0 & 0 & 0 \\ L_1 & 0 & 0 \\ L_2 & L_3 & 0 \end{pmatrix} \middle| L_i \in L_R; 1 \le i \le 3 \right\} \subseteq V$$

be a DSm vector subspace of refined labels of V over R. Clearly $W_1$ $W_2$ $W_3$ are independent for

$$\begin{pmatrix} L_1 & 0 & 0 \\ 0 & L_2 & 0 \\ 0 & 0 & L_3 \end{pmatrix} + \begin{pmatrix} 0 & L_4 & L_5 \\ 0 & 0 & L_6 \\ 0 & 0 & 0 \end{pmatrix} + \begin{pmatrix} 0 & 0 & 0 \\ L_7 & 0 & 0 \\ L_8 & L_9 & 0 \end{pmatrix} = \begin{pmatrix} 0 & 0 & 0 \\ 0 & 0 & 0 \\ 0 & 0 & 0 \end{pmatrix}$$

is not possible unless each of the terms are zero.

$$\text{Further } W_1 \cap W_2 = \begin{pmatrix} 0 & 0 & 0 \\ 0 & 0 & 0 \\ 0 & 0 & 0 \end{pmatrix} = W_2 \cap W_3 = W_1 \cap W_3.$$

Each vector $\alpha = \alpha_1 + \alpha_2 + \alpha_3$ in V can be expressed as a sum and if $\alpha = \beta_1 + \beta_2 + \beta_3$ where $\alpha_i, \beta_i \in W_i$, i = 1, 2, 3 then $\alpha - \alpha = (\alpha_1 - \beta_1) + (\alpha_2 - \beta_2) + (\alpha_3 - \beta_3) = 0$, hence $\alpha_i - \beta_i = 0$; i = 1, 2, 3

Thus $W_1, W_2, W_3$ are independent,

Now the following result is a matter of routine and the reader is expected to prove.

*Lemma 2.1:* Let V be a DSm vector space of refined labels over the field of reals R. Let $W = W_1 + \ldots + W_k$ be the DSm subspace of V spanned by the DSm subspaces of refined labels $W_1, W_2, \ldots, W_k$ of V.

The following are equivalent:

(a) $W_1, W_2, \ldots, W_k$ are independent
(b) For each j, $2 \leq j \leq k$, we have $W_j \cap (W_1 + \ldots + W_{j-1}) = \{0\}$.

By examples 2.26 and 2.27 we see that the above lemma is true.

It is still interesting to study the structure of linear transformation of two DSm vector spaces of refined labels V and W over the real field R. Let V be a DSm vector space of refined labels over the real field R. If T is a DSm linear operator on V and if W is a DSm subspace of V of refined labels, we say W is invariant under T if for each vector a in W the vector Ta is in W; that is T (W) is contained in W.

We will illustrate this situation by some examples.

*Example 2.28:* Let

$$V = \left\{ \begin{pmatrix} L_{a_1} & L_{a_2} \\ L_{a_3} & L_{a_4} \end{pmatrix} \middle| L_{a_i} \in L_R; 1 \leq i \leq 4 \right\}$$

be a DSm vector space of refined labels over the real field R. Let T be a DSm linear operator on V.

Let

$$W = \left\{ \begin{pmatrix} L_{a_1} & 0 \\ 0 & L_{a_2} \end{pmatrix} \middle| L_{a_1}, L_{a_2} \in L_R \right\} \subseteq V.$$

W is a DSm vector subspace of V of refined labels over R.

Now

$$T \left( \begin{pmatrix} L_{a_1} & L_{a_2} \\ L_{a_3} & L_{a_4} \end{pmatrix} \right) = \begin{pmatrix} L_{a_1} & 0 \\ 0 & L_{a_2} \end{pmatrix}.$$

We see $T(W) \subseteq W$; hence W the DSm vector subspace of refined labels is invariant under T.

These DSm subspaces for which T (W) ⊆ W will be known as the invariant DSm subspaces of refined labels invariant under T.

Further from [34-5] it is given for every r ∈ R we have r = $L_a$. If we assume $L_a$ ∈ $L_R$ then we have a unique r ∈ R such that $L_a$ = r.

Thus we see ($L_R$, +, ×, .) is a DSm vector space of refined labels of dimension one over R. Thus all the DSm spaces discussed by us are of finite dimensional expect the DSm polynomial vector space where the coefficients of the polynomials are refined labels.

Now we can as in case of usual vector spaces define the notion of basis and dimension, we call them as DSm basis and DSm dimension, we will proceed onto give examples of them.

***Example 2.29:*** Let V = ($L_R$ × $L_R$ × $L_R$ × $L_R$) = {($L_a$, $L_b$, $L_c$, $L_d$) | $L_a$, $L_b$, $L_c$, $L_d$ ∈ $L_R$} be a DSm linear algebra of refined labels over the real field R.

Consider B = {($L_{m+1}$, 0, 0, 0), (0, $L_{m+1}$, 0, 0), (0, 0, $L_{m+1}$, 0), (0, 0, 0, $L_{m+1}$)} ⊆ V; B is a DSm basis of V and the DSm dimension of V over R is four.

As in case of usual vector spaces we see that a set with more than four non zero elements are linearly dependent. We have several DSm vector spaces of any desired dimension.

***Example 2.30:*** Let

$$V = \left\{ \begin{pmatrix} L_{a_1} & L_{a_2} & L_{a_3} \\ L_{a_4} & L_{a_5} & L_{a_6} \\ L_{a_7} & L_{a_8} & L_{a_9} \end{pmatrix} \middle| L_i \in L_R; 1 \le i \le 9 \right\}$$

be a DSm linear algebra of refined labels over the reals. Consider

$$B = \left\{ \begin{pmatrix} L_a & 0 & 0 \\ 0 & 0 & 0 \\ 0 & 0 & 0 \end{pmatrix}, \begin{pmatrix} 0 & L_b & 0 \\ 0 & 0 & 0 \\ 0 & 0 & 0 \end{pmatrix}, \begin{pmatrix} 0 & 0 & L_c \\ 0 & 0 & 0 \\ 0 & 0 & 0 \end{pmatrix}, \right.$$

$$\begin{pmatrix} 0 & 0 & 0 \\ 0 & L_d & 0 \\ 0 & 0 & 0 \end{pmatrix}, \begin{pmatrix} 0 & 0 & 0 \\ L_e & 0 & 0 \\ 0 & 0 & 0 \end{pmatrix}, \begin{pmatrix} 0 & 0 & 0 \\ 0 & 0 & L_f \\ 0 & 0 & 0 \end{pmatrix},$$

$$\left. \begin{pmatrix} 0 & 0 & 0 \\ 0 & 0 & 0 \\ L_n & 0 & 0 \end{pmatrix}, \begin{pmatrix} 0 & 0 & 0 \\ 0 & 0 & 0 \\ 0 & L_t & 0 \end{pmatrix}, \begin{pmatrix} 0 & 0 & 0 \\ 0 & 0 & 0 \\ 0 & 0 & L_s \end{pmatrix} \right\} \subseteq V;$$

is a subset of V, which is easily verified to be a DSm linearly independent subset of V.

This forms a DSm basis of V over R (none of $L_a$, $L_b$, ..., $L_s$ are 0). Consider

$$W = \left\{ \begin{pmatrix} L_{a_1} & L_{a_2} & 0 \\ 0 & 0 & 0 \\ 0 & 0 & 0 \end{pmatrix}, \begin{pmatrix} 0 & 0 & L_{a_3} \\ 0 & 0 & 0 \\ 0 & 0 & 0 \end{pmatrix}, \right.$$

$$\begin{pmatrix} 0 & 0 & 0 \\ 0 & L_{a_4} & L_{a_5} \\ 0 & 0 & 0 \end{pmatrix}, \begin{pmatrix} 0 & 0 & 0 \\ L_{a_6} & 0 & 0 \\ L_{a_7} & 0 & 0 \end{pmatrix},$$

$$\left. \begin{pmatrix} 0 & 0 & 0 \\ 0 & 0 & 0 \\ 0 & L_{a_8} & L_{a_9} \end{pmatrix} \middle| L_{a_i} \in L_R \text{ and none of } L_{a_i} = 0 \text{ and } 1 \leq i \leq 9 \right\} \subseteq V,$$

W is a DSm linearly independent subset of V but W is not a DSm basis of V over the field R.

***Example 2.31:*** Let

$$V = \left\{ \sum_{i=0}^{\infty} L_{a_i} x^i \middle| L_{a_i} \in R \right\}$$

be a DSm linear algebra of refined labels over the field R. Clearly V is an infinite dimensional DSm linear algebra over R.

We see P = $\{L_{m+1}x^3, L_{m+1}x^7, L_{m+1}x^8, L_{m+1}x^{120}\} \subseteq$ V is a DSm linearly independent subset of V however P is not a DSm basis of V.

Consider B = $\{L_a x^2 + L_b x^3, L_c x^5 + L_t x^9, L_a x^2, L_c x^5\} \subseteq$ V is a DSm linearly dependent subset of V. Now we have seen examples of DSm basis, DSm dimension and DSm linearly independent subset and DSm linearly dependent subset of a DSm linear algebra (vector space) of refined labels over the real field R.

Now we proceed onto define DSm characteristic values or characteristic DSm values or refined label values of a DSm linear algebra (or vector space) of refined labels over R.

We call them as refined label values as we deal with refined labels DSm linear algebra (or vector spaces).

**DEFINITION 2.7:** *Let V be a DSm vector space of refined labels over the real field F. Let T be a DSm linear operator on V. A DSm characteristic value or characteristic DSm value or refined label value associated with T (or of T) is a scalar c in R such that there is a non zero label vector $\alpha$ in V with $T\alpha = c\alpha$. If c is a DSm characteristic value of T or characteristic DSm value or characteristic refined label value of T then*

*(a) for any $\alpha$ in V such that $T\alpha = c\alpha$ is called the DSm characteristic DSm vector of T or characteristic refined*

*label vector of T associated with the characteristic refined label value c of R.*

(b) *the collection of all DSm characteristic vectors $\alpha$ such that $T\alpha = c\alpha$ is called the characteristic DSm space of refined labels associated with c.*

We can also as in case of usual vector spaces call the DSm characteristic values (or refined label values) as DSm characteristic roots (or refined label roots), DSm latent roots (or refined latent label roots), eigen label values (or DSm eigen values), DSm proper values (proper refined label values) or DSm spectral values (refined spectral label values).

We call them in this book as refined characteristic label values or just refined label values.

Let T be a DSm linear operator on a DSm vector space V of refined values over the field R. If for any scalar c in R, the set of label vectors $\alpha$ such that $T\alpha = c\alpha$ is a DSm subspace of refined labels of V.

It is the DSm null space of the DSm linear transformation (T – cI). We call c the DSm characteristic value of T if this subspace is different from the zero subspace i.e. if (T – cI) fails to be one to one. If the underlying DSm space of refined labels V is finite dimensional, (T – cI) fails to be one to one only when its determinant value is different from zero.

We assume all properties enjoyed by the usual vector spaces are true in case of DSm vector spaces of refined labels. However we will derive the analogue.

We can use the theorem on finite dimensional vector spaces which says every n-dimensional vector space over the field F is isomorphic to the space $F^n$.

Thus we can say if V is a DSm vector space of refined labels over the reals R of dimension n then
$$V \cong R^n = \underbrace{R \times R \times ... \times R}_{n-\text{times}}.$$

For if $V = (L_R \times L_R)$ a DSm vector space of dimension two over R then $V \cong R \times R$.

Likewise if

$$V = \left\{ \begin{pmatrix} L_{a_1} & L_{a_2} & L_{a_3} \\ L_{a_4} & L_{a_5} & L_{a_6} \end{pmatrix} \middle| L_{a_i} \in L_R ; 1 \leq i \leq 6 \right\}$$

be a DSm vector space of refined labels over the real field R. $V \cong R \times R \times R \times R \times R \times R = R^6$.

For if

$$B = \left\{ \begin{pmatrix} L_{m+1} & 0 & 0 \\ 0 & 0 & 0 \end{pmatrix}, \begin{pmatrix} 0 & L_{m+1} & 0 \\ 0 & 0 & 0 \end{pmatrix}, \begin{pmatrix} 0 & 0 & L_{m+1} \\ 0 & 0 & 0 \end{pmatrix}, \right.$$

$$\left. \begin{pmatrix} 0 & 0 & 0 \\ L_{m+1} & 0 & 0 \end{pmatrix}, \begin{pmatrix} 0 & 0 & 0 \\ 0 & L_{m+1} & 0 \end{pmatrix}, \begin{pmatrix} 0 & 0 & 0 \\ 0 & 0 & L_{m+1} \end{pmatrix} \right\} \subseteq V$$

be a basis of V over the reals R.

Let

$$A = \begin{pmatrix} L_{a_1} & L_{a_2} & L_{a_3} \\ L_{a_4} & L_{a_5} & L_{a_6} \end{pmatrix} \in V.$$

Now A in V can be written using the basis B as follows:

$$A = r_1 \begin{pmatrix} L_{m+1} & 0 & 0 \\ 0 & 0 & 0 \end{pmatrix} + r_2 \begin{pmatrix} 0 & L_{m+1} & 0 \\ 0 & 0 & 0 \end{pmatrix} + r_3 \begin{pmatrix} 0 & 0 & L_{m+1} \\ 0 & 0 & 0 \end{pmatrix}$$

$$+ r_4 \begin{pmatrix} 0 & 0 & 0 \\ L_{m+1} & 0 & 0 \end{pmatrix} + r_5 \begin{pmatrix} 0 & 0 & 0 \\ 0 & L_{m+1} & 0 \end{pmatrix} + r_6 \begin{pmatrix} 0 & 0 & 0 \\ 0 & 0 & L_{m+1} \end{pmatrix}$$

$$= \begin{pmatrix} L_{r_1 m+1} & L_{r_2 m+1} & L_{r_3 m+1} \\ L_{r_4 m+1} & L_{r_5 m+1} & L_{r_6 m+1} \end{pmatrix}$$

$$= \begin{pmatrix} L_{a_1} & L_{a_2} & L_{a_3} \\ L_{a_4} & L_{a_5} & L_{a_6} \end{pmatrix}$$

as $r_i = \dfrac{a_i}{m+1}$; i =1, 2, 3, 4, 5, 6.

Thus we can say

$$V \cong \underbrace{R \times ... \times R}_{6-\text{times}}.$$

Likewise one can show for any DSm n dimensional DSm vector space of refined labels over R is isomorphic to

$$\underbrace{R \times R \times ... \times R}_{n-\text{times}} = R^n.$$

We can as in case of usual vector spaces derive properties about DSm vector space of refined labels over the reals R with appropriate modifications. The following theorems can be proved as a matter of routine.

**THEOREM 2.9:** *Let V be a DSm vector space refined labels over the reals. Intersection of any collection of DSm subspaces of V is a DSm subspace of refined labels of V over the reals.*

**THEOREM 2.10:** *Let V be a DSm vector space of refined labels. The DSm subspace spanned by a non empty subset S of the DSm vector space V is the set of all linear combinations of vectors in S.*

We just give an hint for the proof. If we assume W to be a DSm vector subspace of V spanned by the set S then each linear combination a = $x_1 \alpha_1 + x_2 \alpha_2 + ... + x_n \alpha_m$; $\alpha_1, \alpha_2, ..., \alpha_m \in$ S (set of vectors in S); $x_i \in$ R; $1 \leq i \leq n$ is in W.

Hence the claim [34-5, 37].

**THEOREM 2.11:** *Let V be a DSm vector space of refined labels over R. If V is spanned as a DSm vector space by $\beta_1, \beta_2, ..., \beta_m$. Then any independent set of vectors in V is finite and contains no more than m elements.*

**THEOREM 2.12:** *Let V be a DSm vector space of refined labels of finite dimension over the real field R, then any two basis of V have the same number of elements.*

**THEOREM 2.13:** *Let V be a DSm vector space of refined labels of finite dimension over the field R. If n = dimension of V. Then (a) any subset of V which contains more than n elements is linearly dependent (b) No subset of V which contains less than n vectors can span V.*

**THEOREM 2.14 :** *Let V be a DSm vector space of refined labels over the field R. Suppose β is a vector in V which is not in the DSm subspace spanned by S then the set obtained by adjoining β to S is linearly independent.*

**THEOREM 2.15:** *Let W be a DSm subspace of refined labels of the DSm finite dimensional vector space V of refined labels over R, then every linearly independent subset of W is finite and is part of a (finite) basis for W.*

**COROLLARY 2.1:** *If in a finite dimensional DSm vector space of refined labels V over the field R. W is a proper DSm vector subspace of V over R then dim W < dim V.*

**COROLLARY 2.2**: *In any finite dimensional DSm vector space V of refined labels every non empty linearly independent set of vectors is part of a basis.*

Now we will leave the proof of all these theorems and corollaries to the reader. However we will give some illustrative examples of them.

*Example 2.32:* Let

$$V = \left\{ \begin{pmatrix} L_{a_1} & L_{a_2} \\ L_{a_3} & L_{a_4} \\ L_{a_5} & L_{a_6} \end{pmatrix} \middle| L_{a_i} \in L_R; 1 \leq i \leq 6 \right\}$$

be a DSm vector space of refined labels over the real field R.

Consider

$$W_1 = \left\{ \begin{pmatrix} L_{a_1} & L_{a_2} \\ 0 & 0 \\ L_{a_3} & 0 \end{pmatrix} \middle| L_{a_i} \in L_R; 1 \le i \le 3 \right\},$$

$$W_2 = \left\{ \begin{pmatrix} 0 & L_{a_1} \\ 0 & L_{a_2} \\ L_{a_3} & 0 \end{pmatrix} \middle| L_{a_i} \in L_R; 1 \le i \le 3 \right\}$$

and

$$W_3 = \left\{ \begin{pmatrix} L_{a_1} & 0 \\ L_{a_2} & L_{a_3} \\ L_{a_4} & 0 \end{pmatrix} \middle| L_{a_i} \in L_R; 1 \le i \le 3 \right\}$$

be three DSm vector subspaces of refined labels over the field R. We see

$$\bigcap_{i=1}^{3} W_i = \left\{ \begin{pmatrix} 0 & 0 \\ 0 & 0 \\ L_a & 0 \end{pmatrix} \ne \begin{pmatrix} 0 & 0 \\ 0 & 0 \\ 0 & 0 \end{pmatrix} \middle| L_{a_i} \in L_R \right\} \subseteq V$$

is a DSm vector subspace of V over R of refined labels.

***Example 2.33***: Let $V = \{(L_{a_1}, L_{a_2}, ..., L_{a_{15}}) \mid L_{a_i} \in L_R\}$ be a DSm vector space of refined labels over R. Consider the set $S = \{(L_{a_1}, 0, 0, 0, 0, L_{a_6}, 0, 0, 0, 0, L_{a_{12}}), (L_{a_1}, L_{a_2}, L_{a_3}, L_{a_4}, 0, ..., 0), (0, 0, 0, 0, 0, L_{a_1}, L_{a_2}, L_{a_3}, 0, ..., 0, L_{a_{11}}, 0), (0, L_{a_1}, 0, 0, 0, 0, 0, 0, 0, 0, L_{a_{10}}, 0, 0)\} \subseteq V$ be a proper subset of V. W be the DSm subspace of refined labels spanned by S. Now $W = \{(L_{a_1}, L_{a_2}, L_{a_3}, L_{a_4}, 0, L_{a_6}, L_{a_7}, L_{a_8}, 0, L_{a_{10}}, L_{a_{11}}, L_{a_{12}}) \mid L_{a_i} \in L_R; i = 1, 2, 3, 4, 6, 7, 8, 10, 11,$ and $12\} \subseteq V$ is a DSm vector subspace of refined labels spanned by S.

We see every element in W is a linear combination of set of vectors from S.

*Example 2.34:* Let

$$V = \left\{ \begin{pmatrix} L_{a_1} & L_{a_2} \\ L_{a_3} & L_{a_4} \end{pmatrix} \middle| L_{a_i} \in L_R; 1 \le i \le 4 \right\}$$

be a DSm vector space of refined labels over the field R. Clearly DSm dimension of V is four.

Consider a set

$$S = \left\{ \begin{pmatrix} L_{a_1} & L_{a_2} \\ 0 & 0 \end{pmatrix}, \begin{pmatrix} L_{a_1} & L_{a_2} \\ L_{a_3} & 0 \end{pmatrix}, \begin{pmatrix} L_{a_1} & 0 \\ 0 & 0 \end{pmatrix}, \right.$$

$$\left. \begin{pmatrix} 0 & 0 \\ L_{a_1} & L_{a_2} \end{pmatrix}, \begin{pmatrix} L_{a_1} & L_{a_2} \\ 0 & L_{a_3} \end{pmatrix} \right\} \subseteq V.$$

Clearly S is a subset of V with cardinality five and we see this set is a linearly dependent subset of V.

*Example 2.35:* Let

$$V = \left\{ \begin{pmatrix} L_{a_1} & L_{a_2} \\ L_{a_3} & L_{a_4} \\ L_{a_5} & L_{a_6} \\ L_{a_7} & L_{a_8} \\ L_{a_9} & L_{a_{10}} \\ L_{a_{11}} & L_{a_{12}} \end{pmatrix} \middle| L_{a_i} \in L_R; 1 \le i \le 12 \right\}$$

be a DSm vector space of refined labels over the field R.
Consider

$$W = \left\{ \begin{pmatrix} L_{a_1} & L_{a_2} \\ 0 & 0 \\ L_{a_3} & L_{a_4} \\ 0 & 0 \\ L_{a_5} & L_{a_6} \\ 0 & 0 \end{pmatrix} \middle| L_{a_i} \in L_R; 1 \le i \le 6 \right\} \subseteq V$$

be a DSm vector subspace of V over the field R.

Consider

$$\beta = \begin{pmatrix} 0 & 0 \\ L_{a_1} & L_{a_2} \\ 0 & 0 \\ L_{a_3} & 0 \\ 0 & 0 \\ 0 & 0 \end{pmatrix} \in V;$$

the subset S which spans W by no means can contain $\beta$ so the subset $S \cup \beta$ is a linearly independent subset of V.

*Example 2.36:* Let

$$V = \left\{ \begin{pmatrix} L_{a_1} & L_{a_2} & L_{a_3} \\ L_{a_4} & L_{a_5} & L_{a_6} \\ L_{a_7} & L_{a_8} & L_{a_9} \\ L_{a_{10}} & L_{a_{11}} & L_{a_{12}} \end{pmatrix} \middle| L_{a_i} \in L_R; 1 \le i \le 12 \right\}$$

be a DSm vector space of refined labels over the field R.

$$W = \left\{ \begin{pmatrix} L_{a_1} & L_{a_2} & L_{a_3} \\ 0 & 0 & 0 \\ L_{a_4} & L_{a_5} & L_{a_6} \\ 0 & 0 & 0 \end{pmatrix} \middle| L_{a_i} \in L_R; 1 \le i \le 6 \right\} \subseteq V$$

be a proper DSm vector subspace of V over R.

Now DSm dimension of V is 12 and the DSm dimension of W as a subspace of V over R is 6. Clearly dim W < dim V.

*Example 2.37:* Let

$$V = \left\{ \begin{pmatrix} L_{a_1} & L_{a_2} & L_{a_3} & L_{a_4} & L_{a_5} & L_{a_6} \\ L_{a_7} & L_{a_8} & L_{a_9} & L_{a_{10}} & L_{a_{11}} & L_{a_{12}} \\ L_{a_{13}} & L_{a_{14}} & L_{a_{15}} & L_{a_{16}} & L_{a_{17}} & L_{a_{18}} \end{pmatrix} \middle| L_{a_i} \in L_R; 1 \le i \le 18 \right\}$$

be a DSm vector space of refined labels of dimension 18 over the field R.

Consider

$$T = \left\{ \begin{pmatrix} L_{a_1} & 0 & 0 & 0 & 0 & 0 \\ 0 & 0 & 0 & 0 & 0 & 0 \\ 0 & 0 & 0 & 0 & 0 & 0 \end{pmatrix}, \begin{pmatrix} 0 & 0 & 0 & 0 & L_{a_2} & 0 \\ 0 & 0 & 0 & 0 & 0 & 0 \\ 0 & 0 & 0 & 0 & 0 & 0 \end{pmatrix}, \right.$$

$$\begin{pmatrix} 0 & 0 & 0 & 0 & 0 & 0 \\ 0 & L_{a_3} & 0 & 0 & 0 & 0 \\ 0 & 0 & 0 & 0 & 0 & 0 \end{pmatrix}, \begin{pmatrix} 0 & 0 & 0 & 0 & 0 & 0 \\ 0 & 0 & 0 & 0 & 0 & 0 \\ L_{a_4} & 0 & 0 & 0 & 0 & 0 \end{pmatrix},$$

$$\left. \begin{pmatrix} 0 & 0 & 0 & 0 & 0 & 0 \\ 0 & 0 & 0 & 0 & 0 & 0 \\ 0 & 0 & 0 & 0 & 0 & L_{a_5} \end{pmatrix}, \begin{pmatrix} 0 & 0 & 0 & 0 & 0 & L_{a_6} \\ 0 & 0 & 0 & 0 & 0 & 0 \\ 0 & 0 & 0 & 0 & 0 & 0 \end{pmatrix} \right\} \subseteq V$$

is a set of linearly independent subset of V of refined labels we see |T| = 6 but DSm dim V = 18 and clearly T forms a part of the basis of V over R.

Now having seen examples of these theorems and corollaries we now proceed onto define and give examples of more new concepts on DSm vector space of refined labels over the field R.

*Example 2.38:* Let

$$V = \left\{ \begin{pmatrix} L_{a_1} & L_{a_2} & L_{a_3} \\ L_{a_4} & L_{a_5} & L_{a_6} \\ L_{a_7} & L_{a_8} & L_{a_9} \end{pmatrix} \middle| L_{a_i} \in L_R; 1 \leq i \leq 9 \right\}$$

be a DSm vector space of refined labels over R. Consider

$$W_1 = \left\{ \begin{pmatrix} 0 & L_{a_1} & L_{a_2} \\ 0 & 0 & L_{a_3} \\ 0 & L_{a_4} & 0 \end{pmatrix} \middle| L_{a_i} \in L_R; 1 \leq i \leq 3 \right\}$$

and

$$W_2 = \left\{ \begin{pmatrix} L_{a_1} & L_{a_2} & 0 \\ 0 & L_{a_3} & 0 \\ L_{a_4} & 0 & L_{a_5} \end{pmatrix} \middle| L_{a_i} \in L_R; 1 \le i \le 5 \right\}$$

be two DSm vector subspaces of refined labels over the field R. Now DSm dimension of V is 9. DSm dimension of $W_1$ is four where as DSm dimension of $W_2$ is 5. Clearly

$$W_1 \cap W_2 = \left\{ \begin{pmatrix} 0 & L_{a_1} & 0 \\ 0 & 0 & 0 \\ 0 & 0 & 0 \end{pmatrix} \middle| L_{a_i} \in L_R \right\} \subseteq V$$

is a DSm subspace of V over R. Further dim $W_1$ + dim $W_2$ = dim ($W_1 \cap W_2$) + dim ($W_1 + W_2$) for DSm dim ($W_1 \cap W_2$) = one, DSm dim $W_1$ = 4 and DSm dim $W_2$ = 5. Thus 5 + 4 = 1 + dim ($W_1 + W_2$) we see DSm dim ($W_1 + W_2$) = 9 – 1 = 8 as the DSm subspace generated by $W_1$ and $W_2$ is given by

$$T = \left\{ \begin{pmatrix} L_{a_1} & L_{a_2} & L_{a_3} \\ 0 & L_{a_4} & L_{a_5} \\ L_{a_6} & L_{a_7} & L_{a_8} \end{pmatrix} \middle| L_{a_i} \in L_R; 1 \le i \le 8 \right\} \subseteq V$$

and dimension of T is 8. Hence the claim. This result or conclusion of this example can be stated as a Theorem and the proof is direct.

**THEOREM 2.16:** *Let $W_1$ and $W_2$ be any two finite dimensional DSm vector subspaces of a DSm vector space V of refined labels over the reals R. Then dim $W_1$ + dim $W_2$ = dim ($W_1 \cap W_2$) + dim ($W_1 + W_2$).*

Let V and W be DSm vector spaces of refined labels over the field R. A DSm linear transformation from V into W is a function T from V into W such that
$$T(c\alpha + \beta) = cT(\alpha) + T(\beta) = c(T_\alpha) + (T_\beta)$$
for all $\alpha, \beta \in V$ and $c \in R$. We know this but we have not defined DSm rank space of T and DSm null space of T in the formal way through we have used these concepts informally.

Let V and W be any two DSm vector spaces of refined labels over the field R. T a DSm linear transformation from V into W. The null space or more formally the DSm null space of T is the set of all refined label vector $\alpha$ in V such that $T\alpha = 0$.

If V is finite dimensional DSm vector space, the DSm rank of T is the dimension of the DSm range of T and the DSm nullity of T is the dimension of the DSm null space of T.

We will first illustrate this situation by an example before we proceed onto give a theorem.

***Example 2.39:*** Let

$$V = \left\{ \begin{pmatrix} L_{a_1} & L_{a_2} \\ L_{a_3} & L_{a_4} \\ L_{a_5} & L_{a_6} \\ L_{a_7} & L_{a_8} \\ L_{a_9} & L_{a_{10}} \end{pmatrix} \middle| L_{a_i} \in L_R; 1 \le i \le 10 \right\}$$

be a DSm vector space of refined labels over the field R.
Let

$$W = \left\{ \begin{pmatrix} L_{a_1} & L_{a_2} & L_{a_3} & L_{a_4} & L_{a_5} \\ L_{a_6} & L_{a_7} & L_{a_8} & L_{a_9} & L_{a_{10}} \end{pmatrix} \middle| L_{a_i} \in L_R; 1 \le i \le 10 \right\}$$

be a DSm vector space of refined labels over the field R. Suppose $T : V \to W$ be a DSm linear transformation of V to W given by

$$T\left( \begin{pmatrix} L_{a_1} & L_{a_2} \\ L_{a_3} & L_{a_4} \\ L_{a_5} & L_{a_6} \\ L_{a_7} & L_{a_8} \\ L_{a_9} & L_{a_{10}} \end{pmatrix} \right) = \begin{pmatrix} L_{a_1} & 0 & L_{a_5} & 0 & L_{a_9} \\ L_{a_6} & 0 & L_{a_6} & 0 & L_{a_{10}} \end{pmatrix}.$$

Clearly DSm null space of T is

$$\left\{ \begin{pmatrix} 0 & 0 \\ L_{a_3} & L_{a_4} \\ 0 & 0 \\ L_{a_7} & L_{a_8} \\ 0 & 0 \end{pmatrix} \middle| L_{a_3}, L_{a_4}, L_{a_7}, L_{a_8} \in L_R \right\}.$$

Thus the DSm null space of T is a subspace of V. Now DSm nullity T = 4 and DSm rank T = 6. Thus DSm (T) + DSm nullity (T) = DSm dim V = 10 = 6 + 4.

Now we give the following theorem and expect the reader to prove it.

**THEOREM 2.17:** *Let V and W be any two DSm vector spaces of refined labels over the field R and let T be a DSm linear transformation from V into W. Suppose V is finite dimensional then, DSm rank (T) + DSm nullity (T) = DSm dim V.*

Further we see if T and P are any two DSm linear transformations of DSm vector spaces of refined labels V and W defined over R. Then (T + P) is a function defined by (T + P) $\alpha$ = T$\alpha$ + P$\alpha$ for all $\alpha \in$ V is a DSm linear transformation of V into W. If c is any element of the field R the function cT defined by (cT) ($\alpha$) = c(T$\alpha$) is a DSm linear transformation from V into W. Also the set of all DSm linear transformations from V into W together with addition and multiplication defined above is again a vector space over the field F.

DSm L (V, W) = {T : V $\to$ W, where V and W are DSm vector spaces over the field R}, denotes the collection of all DSm linear transformations from V to W.

Now as in case of usual vector spaces we in case of DSm vector spaces also define composition of DSm linear transformations. Let V, W and Z be vector spaces of refined labels over the field of reals. Let T : V $\to$ W be a DSm linear transformation, P : W $\to$ Z be a DSm linear transformation of W to Z. We can define the composed function PT defined by PT $(\alpha) = P(T(\alpha))$ is a DSm linear transformation from V into Z. This can be proved by any interested reader for any general DSm vector spaces V, W and Z. Now we know if W in a DSm linear transformation from V to W is taken as V that $T : V \to V$ we call T a DSm linear operator on V. We have the following properties to be true also in case of DSm linear operators.

**THEOREM 2.18:** *Let V be DSm vector space of refined labels over the field R. Let $T_1$, $T_2$ and S be DSm linear operators on V and let c be an element of R.*

(a)  $IS = SI = S$ for all $S : V \to V$ where I is the DSm identity operator on V.

(b)  $S(T_1 + T_2) = ST_1 + ST_2$ and $(T_1 + T_2)S = T_1S + T_2S$

(c)  $c(ST_1) = (cS)T = S(cT)$.

The proof is a matter of routine and hence is left as an exercise to the reader. As in case of usual vector spaces one can easily prove the following theorems.

**THEOREM 2.19:** *Let T be a DSm linear transformation of the DSm vector space V into W then T is non singular if and only if T carries each linearly independent subset of V onto a linearly independent subset of W.*

**THEOREM 2.20**: *Let V and W be any two finite dimensional DSm vector spaces over the field R such that DSm dim V = DSm dim W. If T is a DSm linear transformation from V into W then the following are equivalent.*

(i)  *T is invertible.*
(ii)  *T is non singular.*
(iii)  *T is onto that is range of T is W.*

Next we proceed onto define DSm linear functionals on the DSm vector space V over the real field R. This concept is also defined in an analogous way. Let V be a DSm vector space over the real field R. Let f be a function from V to R such that
$$f(c\alpha + \beta) = c f(\alpha) + f(\beta)$$
where $\alpha, \beta \in V$ are labels and $c \in R$ is a real.

Let V* = {collection of a DSm linear functionals of V}
    = DSm L (V, R).

Then as in case of usual vector spaces in case of DSm vector spaces also we have V* and V are of same dimension if V is finite dimensional.

We will first illustrate this by one or two examples.

***Example 2.40:*** Let $V = (L_R \times L_R \times L_R) = \{(L_a, L_b, L_c) / L_a, L_b, L_c \in L_R\}$ be a DSm vector space of refined labels over R.

Define $f : V \to R$ by for $v \in V$.
$$f(v) = f((L_a, L_b, L_c)) = (r_1 + r_2 + r_3)$$
where
$$r_1 = \frac{a}{m+1}, \; r_2 = \frac{b}{m+1} \text{ and } r_3 = \frac{c}{m+1}$$
such that $L_a = r_1$, $L_b = r_2$ and $L_c = r_3$; f is a linear functional on V.

***Example 2.41:*** Let
$$V = \left\{ \begin{pmatrix} L_a & L_b & L_c \\ L_d & L_e & L_f \end{pmatrix} \text{ where } L_a, L_b, L_c, L_d, L_e, L_f \in L_R \right\}$$
be a DSm vector space of refined labels over the reals R. Define for any $v = \begin{pmatrix} L_{a_1} & L_{a_2} & L_{a_3} \\ L_{a_4} & L_{a_5} & L_{a_6} \end{pmatrix} \in V$, $f : V \to R$ by $f(v) = r_1 + r_2 + r_3 + r_4 + r_5 + r_6$ where $L_{a_i} = r_i$; $1 \le i \le 6$, f is a linear functional from V to R.

Now having seen examples linear functionals one can define the DSm basis of V and V* and give examples of them. We can define as in case of vector spaces the hyperspace for DSm vector spaces of refined labels over the field R.

Let V be a DSm vector space of refined labels over the field R. Let V be a DSm dimension n-dimensional space over R. Let $W \subseteq V$, if W is a DSm (n–1) - dimensional subspace of V over R then we define W to be a DSm hypersubspace of refined labels of V over R.

We will illustrate this situation by some examples.

**Example 2.42:** Let $V = \{(L_R, L_R, L_R, L_R, L_R) = (L_a, L_b, L_c, L_d, L_e) \mid L_a, L_b, L_c, L_d, L_e, \in L_R\}$ be a DSm vector space of refined labels over the field R. Consider

$$W = \left\{ \begin{pmatrix} L_{a_1} & L_{a_2} & 0 & L_{a_3} & L_{a_4} \end{pmatrix} \middle| L_{a_i} \in L_R; 1 \leq i \leq 4 \right\} \subseteq V,$$

is a DSm refined label subspace of dimension four of V over R.

**Example 2.43:** Let

$$V = \left\{ \begin{pmatrix} L_{a_1} & L_{a_2} \\ L_{a_3} & L_{a_4} \end{pmatrix} \middle| L_{a_i} \in L_R; 1 \leq i \leq 4 \right\}$$

be a DSm vector space of refined labels over R. Let

$$W = \left\{ \begin{pmatrix} 0 & L_{a_1} \\ L_{a_2} & L_{a_3} \end{pmatrix} \middle| L_{a_i} \in L_R; 1 \leq i \leq 3 \right\} \subseteq V.$$

W is a DSm hyper subspace of V of refined labels over R of dimension 3. We see V is of DSm dimension 4. Consider

$$T = \left\{ \begin{pmatrix} L_{a_1} & L_{a_2} \\ L_{a_3} & 0 \end{pmatrix} \middle| L_{a_i} \in L_R; 1 \leq i \leq 3 \right\} \subseteq V$$

is again a DSm vector hyper subspace of V of DSm dimension three over R.

However T and W are isomorphic as DSm vector spaces of refined labels.

Now we define some DSm space of polynomial rings of finite degree. Let

$$V = \left\{ \sum_{i=0}^{m} L_{a_i} x^i \, \middle| \, m < \infty, \, L_{a_i} \in L_R \right\}$$

be a DSm vector space of refined labels over R. Consider

$$W = \left\{ \sum_{i=0}^{m-1} L_{a_i} x^i \, \middle| \, L_{a_i} \in L_R \right\} \subseteq V;$$

W is a DSm hyper subspace of V of DSm dimension m and DSm dimension of V is m + 1.

Now having see examples of DSm hyper subspaces of refined labels we now proceed onto define the concept of DSm annihilator of refined labels of annihilator of a subset of a DSm vector space of refined labels over R.

Let V be a DSm vector space of refined labels over the field R and S be a proper subset of V, the annihilator of S is the set $S^o$ of DSm linear functionals f on V such that f(a) = 0 for every a in S.

It is obvious to the reader $S^o$ is a DSm subspace of $V^*$, whether S is a DSm subspace of V or not. If S = {zero vector alone} then $S^o = V^*$. If S = V then $S^o$ is the zero subspace of $V^*$.

The following results are obvious and hence is left as an exercise for the reader to prove.

**THEOREM 2.21:** *Let V be a DSm finite dimensional vector space of refined labels over the field R. Let W be a DSm vector subspace of V over R.*
*Then DSm dim W + DSm dim $W^o$ = DSm dim V.*

**COROLLARY 2.3:** *If W is a k-dimensional subspace of refined labels of a DSm n-dimensional vector space of refined labels of the finite dimensional vector space V of refined labels then W is the intersection of (n–k) DSm hyper subspaces in V.*

**COROLLARY 2.4:** *If $W_1$ and $W_2$ are DSm vector subspaces of refined labels of V over the field R then $W_1 = W_2$ if and only if $W_1^o = W_2^o$.*

Interested reader is expected to give examples of these results.

We can as in case of vector spaces define the notion of double dual in case of DSm vector space of refined labels over R.

Let V be a DSm vector space of refined labels over the field R. V* be the DSm dual space of V over R. Consider $V^{**}$ the dual of $V^*$.

If α is a refined label vector in V then α induces a DSm linear functional $L_\alpha$ on $V^*$ defined by

$$L_\alpha(f) = f(\alpha), \text{ f in } V^*$$

The fact that $L_\alpha$ is a linear is just a reformulation of the definition of DSm linear operator in $V^*$.

$$\begin{aligned} L_\alpha(cf + g) &= (cf + g)(\alpha) \\ &= cf(\alpha) + g(\alpha) \\ &= cf(\alpha) + g(\alpha) \\ &= cL_\alpha(f) + L_\alpha(g). \end{aligned}$$

If V is DSm finite dimensional and $\alpha \neq 0$ then $L_\alpha \neq 0$, in otherwords there exists a linear functional f such that $f(\alpha) \neq 0$. The proof is left as an exercise to the reader.

The following theorem is an interesting consequence.

**THEOREM 2.22:** *Let V be a DSm finite dimensional vector space of refined labels over the field of refined labels over the field R. For each vector $\alpha$ in V define $L_\alpha(f) = f(\alpha)$, f in $V^*$. The mapping $\alpha \to L_\alpha$ is then an isomorphism of V onto $V^{**}$.*

This proof also is simple and hence is left as an exercise to the reader.

**COROLLARY 2.5:** *Let V be a DSm finite dimensional vector space of refined labels over the field V. Each basis for $V^*$ is the dual of some basis for V.*

This proof is also direct and is analogous to the one done for vector spaces.

**THEOREM 2.23**: *If S is any subset of the DSm finite dimensional vector space of refined labels of V over the field R then $(S^o)^o$ is the subspace spanned by S.*

This proof is also direct by using the results.

DSm dim W + DSm dim $W^o$ = DSm dim V and DSm dim $W^o$ + DSm dim $W^{oo}$ = DSm dim $V^*$ and using the fact DSm dim V = DSm dim $V^*$ we get the result. It is just interesting to note that if V is a DSm vector space of refined labels a DSm hyperspace in V is a maximal proper subspace of V.

*Example 2.44:* Let

$$V = \left\{ \begin{bmatrix} L_{a_1} & L_{a_2} \\ L_{a_3} & L_{a_4} \\ L_{a_5} & L_{a_6} \\ L_{a_7} & L_{a_8} \\ L_{a_9} & L_{a_{10}} \end{bmatrix} \middle| L_{a_i} \in L_R; 1 \leq i \leq 10 \right\}$$

be a DSm vector space of refined labels of dimension 10.
Consider

$$H = \left\{ \begin{bmatrix} L_{a_1} & L_{a_2} \\ L_{a_3} & L_{a_4} \\ L_{a_5} & L_{a_6} \\ L_{a_7} & L_{a_8} \\ L_{a_9} & 0 \end{bmatrix} \middle| L_{a_i} \in L_R; 1 \leq i \leq 9 \right\} \subseteq V,$$

H is a DSm vector subspace of V of refined labels of dimension 9. Clearly the maximum DSm dimension this vector space V can have for its vector subspaces W is 9 hence W is a maximum proper subspace of V.

**THEOREM 2.24:** *If f is a non zero DSm linear functional on the vector space of refined labels V, then the null space of f is a DSm hypersubspace of refined labels in V. Conversely every DSm hyperspace of refined labels in V is the null space of a non zero linear functional on V.*

The proof is also direct and hence is left as an exercise for the reader to prove.

**THEOREM 2.25:** *If g and f are two DSm linear functionals on a DSm vector space of refined labels over the field R, then g is a*

scalar multiple of f if and only if the null space of g contains the null space of f that is if and only if f ($\alpha$) = 0 implies g ($\alpha$) = 0.

The proof is also direct and hence is left for the reader to prove.

Recall if f is a DSm linear functional on a DSm vector space of refined labels over R then $N_f$ is the null space of f and $N_f$ is the hyperspace in V.

**THEOREM 2.26:** *Let g, $f_1$, $f_2$, …, $f_r$ be DSm linear functionals on the DSm vector space V of refined labels over R with respective null spaces $N_1$, $N_2$, …, $N_r$. Then g is a linear combination of $f_1$, $f_2$, …, $f_r$ if and only if N contains the intersection $N_1 \cap N_2 \cap … \cap N_r$.*

This proof is also simple and direct and hence is left as an exercise to the reader.

Now as in case of vector spaces we can in case of DSm vector spaces also define the notion of DSm transpose of a linear transformation. Let V and W be two DSm vector spaces of refined labels over the field R and T a DSm linear transformation from V into W. Then T induces a linear transformation from $W^*$ into $V^*$ as follows. Suppose g is a linear functional on W, and let f($\alpha$) = g(T$\alpha$) for each $\alpha$ in V, then f($\alpha$) = g(T$\alpha$) defines a function f from V into R viz., the composition of T, a function from V into W, with g a function from W into F. Since both T and g are linear f is also linear functional on V. Thus T provides us with a rule $T^t$ which associates with each linear functional g on W a linear functional f = $T^t$ g on V defined by f = $T^t$ g on V. Note $T^t$ is actually a linear transformation from $W^*$ into $V^*$.

Now in view of this we have the following theorem.

**THEOREM 2.27:** *Let V and W be two DSm vector spaces of refined labels over the field R. For each linear transformation T from V into W there is a unique linear transformation $T^t$ from $W^*$ into $V^*$ such that $(T^t_g)$ ($\alpha$) = g ($T_\alpha$) for every g in $W^*$ and $\alpha$ in V.*

We also call the DSm transpose of a linear transformation $T^t$ of the linear transformation as adjoint of T.

The following theorem is also direct and hence is left for the reader as a simple exercise.

**THEOREM 2.28:** *Let V and W be any two DSm vector spaces of refined labels over the real field R. Let T be a DSm linear transformation of V into W. The DSm null space of $T^t$ is the annihilator of the range of T. If V and W are finite dimensional then*
  (i)   *rank $(T^t)$ = rank T,*
  (ii)  *the range of $T^t$ is the annihilator of the null space of T.*

We have defined the notion of DSm polynomial with refined label coefficients.
That is
$$V = \left\{ \sum_{i=0} L_{a_i} x^i \,\bigg|\, L_{a_i} \in L_R \right\}$$
where $L_R$ is the field of refined labels.
Let
$$p(x) = \sum_{i=0}^{n} L_{a_i} x^i = L_{a_0} + L_{a_1} x + \ldots + L_{a_n} x^n$$
where $L_{a_0}, L_{a_1}, \ldots, L_{a_n} \in L_R$ where x is a variable. If $L_{a_n} \neq 0$ then we define p (x) is of degree n.

As in case of usual polynomial with real coefficients we in case of polynomials with refined coefficients also have the following result.
If
$$f(x) = \sum_{i=0}^{m} L_{a_i} x^i \text{ and } d(x) = \sum_{i=0}^{n} L_{b_i} x^i$$
are DSm polynomials in
$$L_R[x] = \left\{ \sum_{i=0}^{\infty} L_{a_i} x^i \,\bigg|\, L_{a_i} \in R \right\},$$

m < ∞ and n < ∞ such that deg d < deg f. Then there exists a DSm polynomial g(x) in $L_R[x]$ such that either f – dg = 0 or deg (f – dg) < deg f.

Also the following theorem is direct can be proved as in case of usual polynomial rings.

**THEOREM 2.29**: *If f and d are polynomials over the DSm field $L_R$ that is f, d ∈ $L_R$ [x] and d is different from zero, then there exists refined label coefficient polynomials q, r ∈ $L_R$ [x] such that*

    (i)    *f = dq + r*
    (ii)    *either r = 0 or deg r < deg d.*

*The refined label coefficient polynomials*

$$q(x) = \sum_i L_{a_i} x^i \text{ and } r(x) = \sum_i L_{r_i} x^i$$

*given in (1) and (ii) are unique.*

Now as in case of usual polynomials with real coefficients define in case of polynomials with refined labels coefficients define in case of polynomials with refined label coefficients the notion of quotients, divides and multiple.

Let

$$L_R[x] = \left\{ \sum_i L_{a_i} x^i \,\middle|\, L_{a_i} \in L_R \right\}$$

be a refined label polynomial over the DSm field $L_R$ of refined labels. Let

$$d(x) = \sum_i L_{d_i} x^i$$

be a non zero polynomial over the field $L_R$. If f is in $L_R[x]$ then there is atmost one polynomial

$$q(x) = \sum_i L_{q_i} x^i$$

in $L_R$ [x] such that f(x) = dq. If such a q(x) exists we say that d(x) divides f(x), that f(x) is a multiple of d(x) we call q(x) the quotient of f(x) and d(x). We also write q(x) = f(x) / d(x).

Let $f(x) = \sum_{i} L_{a_i} x^i \in L_R[x]$ be a polynomial over the field $L_R$ of refined labels and let $L_c \in L_R$ be an element of $L_R$. Then f is divisible by $x - L_c$ if and only if $f(L_c) = 0$,. Let $R = L_R$ be the refined field. An element $L_c \in L_R$ is said to be a refined label root or a zero refined label of a given refined coefficients polynomial $f(x)$ over $L_R$ if $f(L_c) = 0$. Now we just assume that the polynomial ring $L_R[x]$ acts live usual rings under differentiation.

For instance if

$$p(x) = \sum_{i}^{7} L_{a_i} x^i = L_{a_0} + L_{a_1} x + L_{a_2} x^2 + \ldots + L_{a_7} x^7$$

is in $L_R[x]$ then we differentiate $p(x)$ as

$$\frac{dp(x)}{dx} = \frac{d}{dx} L_{a_0} + L_{a_1} x + L_{a_2} x^2 + \ldots + L_{a_7} x^7$$
$$= 0 + L_{a_1} x + 2 L_{a_1} x + 3 L_{a_3} x^2 + \ldots + 7 L_{a_7} x^6.$$

$L_{a_1} + L_{2a_2} x + L_{3a_3} x^2 + \ldots + L_{7a_7} x^6$ and $\frac{dp(x)}{dx}$ is in $L_R[x]$, so we see derivative of a refined coefficient polynomial is again a refined coefficient polynomial as basically $L_R \cong R$ (R reals).

So we can differentiate a refined coefficient any desired number of times or even until it is zero as $\frac{d(L_{a_i})}{dx} = 0$ for all $L_{a_i} \in L_R$. Now once we use the concept of derivatives in the refined coefficient polynomials we can have the concept of Taylor's formula to be true and the proof of which is analogous to the proof of the classical Taylor's formula.

Further as $L_R$ is isomorphic with R the reals we see we have the notion of binomial theorem is also true.

That is $L_a, L_b \in L_R$ then

$$(L_a + L_b)^n = \sum_{k=0}^{n} \binom{n}{k} (L_a)^{n-k} (L_b)^k$$

where from [34-5]

so that
$$L_a^n = L_{a^n/(m+1)^{n-1}}$$

$$(L_a + L_b)^n = \sum_{k=0}^{n} \binom{n}{k} (L_a)^{n-k} (L_b)^k$$

$$= (L_a)^n + n(L_a)^{n-1}(L_b) + \binom{n}{2}(L_a)^{n-2}(L_b)^2 + \ldots +$$

$$\binom{n}{t}(L_a)^{n-t}(L_b)^t + \ldots + (L_b)^n$$

$$= L_{a^n/(m+1)^{n-1}} + nL_{a^{n-1}/(m+1)^{n-2}} \cdot L_b + \frac{n(n-1)}{1.2} \cdot L_{a^{n-2}/(m+1)^{n-3}}(L_b)^2$$

$$+ \ldots + \binom{n}{t} L_{a^{n-t}/(m+1)^{n-t-1}} \cdot L_{b^t/(m+1)^{t-1}} + \ldots + L_{b^n/(m+1)^{n-1}}$$

$$= L_{a^n/(m+1)^{n-1}} + \frac{nL_{a^{n-1} \cdot b}}{(m+1)^{n-2}} + L_{\frac{n(n-1)}{1.2}} \cdot L_{\frac{a^{n-2}b^2}{(m+1)^{n-3}}} + \ldots +$$

$$L_{\binom{n}{t}\frac{a^{n-t}b^t}{(m+1)^{n-2}}} + \ldots + L_{b^n/(m+1)^{n-1}}.$$

Thus
$$(L_a + L_b)^2 = (L_a)^2 + (L_b)^2 + 2L_a L_b$$
$$= L_{a^2/m+1} + 2L_{ab/m+1} + L_{b^2/m+1}$$
$$= L_{a^2/m+1} + L_{2ab/m+1} + L_{b^2/m+1}.$$

Now
$$(L_a + L_b)^3 = (L_a)^3 + 3(L_a)^2(L_b) + 3(L_a)(L_b)^2 + (L_b)^3$$
$$= L_{a^3/(m+1)^2} + 3L_{a^2b/(m+1)^2} + 3L_{ab^2/(m+1)^2} + L_{a^3/(m+1)^2}$$

and so on.

Now
$$(L_a + L_b)^4 =$$
$$(L_a)^4 + \frac{4.3.2}{1.2.3}(L_a)^3 L_b + \frac{4.3}{1.2}(L_a)^2(L_b)^2 + \frac{4.3.2}{1.2.3}(L_a)(L_b)^3 + (L_b)^4$$

$$= L_{a^4/(m+1)^3} + 4L_{a^3/(m+1)^2}L_b + 6L_{a^2/m+1}L_{b^2/m+1} + 4L_a L_{b^3/(m+1)^2} + L_{b^4/(m+1)^3}$$

$$= L_{a^4/(m+1)^3} + L_{4ab^3/(m+1)^3} + L_{6a^2b^2/(m+1)^3} + L_{4ab^3/(m+1)^3} + L_{b^4/(m+1)^3}$$

Thus by the very definition of the DSm refined $L_R$ we see the notion of differentiation and binomial theorem can be got with appropriate modifications. Thus we are in a state to give the Taylors formula for polynomial with coefficients from the refined label field $L_R$.

**THEOREM 2.30**: *(Taylors formula) Let $L_R$ be the DSm refined field of characteristic zero, $L_c$ an element in $L_R$ and n a positive integer. If $f(x) = \sum L_{a_i} x^i$ is a polynomial over $L_R$ with $\deg f \leq n$ then $f(x) = \sum_{k=0}^{n} \frac{(D^k f)}{L^k}(L_c)(x - L_c)^k$*

*Hint:* $D, D^2, \ldots, D^n$ are the differential operators of the DSm polynomial

$$f(x) = \sum_{i=0}^{n} L_{a_i} x^i.$$

Now we have just shown

$$(L_a + L_b)^n = \sum_{k=0}^{n} \binom{n}{k} L_a^{n-k} L_b^k$$

$$= L_{a^n/(m+1)^{n-1}} + L_{\frac{na^{n-1}.b}{(m+1)^{n-1}}} + L_{\frac{n(n-1)}{2} \times \frac{a^{n-2}.b^2}{(m+1)^{n-1}}} + L_{b^n/(m+1)^{n-1}}$$

Now
$$\begin{aligned}
x^n &= [L_c + (x - L_c)]^n \\
&= \sum_{k=0}^{n} \binom{n}{k}(L_c)^{n-k}(x - L_c)^k \\
&= (L_c)^n + n.(L_c)^{n-1}(x - L_c) + \ldots + (x - L_c)^n \\
&= L_{c^n/(m+1)^{n-1}} + L_{nc^{n-1}/(m+1)^{n-1}}(x - L_c) + \ldots + (x - L_c)^n.
\end{aligned}$$

If $f = \sum_{k=0}^{n} L_{a_k} x^k$ then

$$D^* f (L_c) = \sum_k L_{a_n} (D^k x^n)(L_c)$$

and

$$\sum_{k=0}^{n} D^k f \frac{L_c}{\lfloor k} (x - L_c)^k = \sum_k \sum_p L_{a_p} \frac{D^k x^p}{\lfloor k} (L_c)(x - L_c)^k$$
$$= \sum_p L_{a_p} x^p = f.$$

We say $L_c$ is a multiple root or of multiplicity r if $(x - L_c)^r$ divides $f(x) = \sum_{i=0}^{n} L_{a_i} x^i$. Clearly $r \leq n$.

Using this simple concept we have the following theorem.

**THEOREM 2.31:** *Let $L_R$ be the DSm field of refined labels of characteristic zero. $f = f(x) = \sum_{i=0}^{n} L_{a_i} x^i$ be a polynomial with refined coefficients with deg $f \leq n$. Then the scalar $L_c$ is a root of f of multiplicity r if and only if $(D^k f)(L_c) = 0$, $0 \leq k \leq r - 1$ and $(D_f^r)(L_c) \neq 0$.*

The proof is analogous to the proof of polynomial with real coefficients with some appropriate changes.

Now we proceed onto define the concept of ideals in $L_R[x]$.

Let $L_R$ be the field of refined labels, $L_R[x]$ be the polynomial in the variable x with coefficients from $L_R$, the DSm refined field of labels. An ideal in $L_R[x]$ is a subspace M of $L_R[x]$ ($L_R[x]$ is a vector space over R and $L_R \cong R$) such that fg is in M when every $f \in L_R[x]$ and $g \in M$.

If
$$d(x) = \sum L_{d_i} x^i$$
is a polynomial in $L_R[x]$ the set $M = d(x) L_R[x]$ of all multiples df of d by any arbitrary f in $L_R[x]$ is an ideal. For M is non empty as
$$d(x) = \sum L_{d_i} x^i \in M.$$

If
$$f(x) = \sum L_{a_i} x^i$$
and
$$g(x) = \sum L_{b_i} x^i$$
are in $L_R[x]$ and $L_c \in L_R[x]$ then $L_c(df) - dg = d(L_c f - g)$ is in M so that M is a DSm subspace. Finally M contains $(df) g = d(fg)$ as well. Thus M is a DSm ideal and is called / defined as the DSm principal ideal generated by d.

We will give some examples of DSm ideals in $L_R[x]$.

***Example 2.45:*** $L_R[x]$ be the polynomial ring with refined label coefficients from $L_R$. Consider $M = \langle (L_{m+1} x^2 + L_a) \rangle$, be the ideal generated by the polynomial $L_{m+1} x^2 + L_a$ as $L_{m+1}$ is the identity element in $L_R$.

***Example 2.46:*** Let $L_R[x]$ be the polynomial ring with refined label coefficients. Let M be the ideal generated by $p(x) = L_{m+1} x^4 + L_a x^2 + L_a$ where $L_a = L_{m+1}$.

Both the ideals given in examples 2.46 and 2.47 are principal ideals.

We can say if $L_R[x]$ is a polynomial ring with refined coefficients and if M is an ideal then M is generated by a polynomial $L_{m+1} x^n + L_{a_{n-1}} x^{n-1} + \ldots + L_1 x + L_0$ where $L_{m+1}$ is the unit in $L_R$.

Now as in case of usual polynomials we can define greatest common divisors to the polynomials in $L_R [x]$.

Let $p_1(x), p_2(x), \ldots, p_n(x)$ be polynomials with refined label coefficients from $L_R$ not all of which is zero; that is $p_i(x) \in L_R[x]; 1 \leq i \leq n$.
The monic generator
$$d(x) = \sum_i L_{d_i} x^i$$
of the ideal $p_1(x) L_R [x] + p_2(x) L_R [x] + \ldots + p_n(x) L_R [x]$ is called the greatest common divisor (g.c.d) of $p_1(x), p_2(x), \ldots, p_n(x)$. This terminology is justified by the following corollary.

**COROLLARY 2.6:** *If $p_1, p_2, \ldots, p_n$ are polynomials over the field $L_R$, not all of which are zero, there is a unique monic polynomial d in $L_R [x]$ such that*
   *(a) d is in the ideal generated by $p_1, \ldots, p_n$;*
   *(b) d divides each of the polynomials $p_i$.*

Any polynomial satisfying (a) and (b) necessarily satisfies (c). d is divisible by every polynomial which divides each of the polynomials $p_1, p_2, \ldots, p_n$. We say the polynomials $p_1(x), p_2(x), \ldots, p_n(x)$ are relatively prime if their greatest common divisor is one or equivalently if the ideal they generate is all of $L_R [x]$.

Now we just describe the concept of reducibility in case of $L_R [x]$. Let $L_R$ be the field of refined labels. A polynomial $f(x)$ in $L_R[x]$ is said to be reducible over $L_R$ if there exists polynomials $g(x), h(x) \in L_R[x]$ in $L_R [x]$ of degree $\geq 1$ such that $f = gh$ and if not f is said to be irreducible over $L_R$. A non scalar irreducible polynomial over $L_R$ is called a prime polynomial over $L_R$ and we say it is a prime in $L_R [x]$. We just mention a theorem the proof of which is left as an exercise to the reader.

**THEOREM 2.32**: *Let f (x) be a polynomial over the DSm refined labels $L_R$ with derivative f ´(x). Then f is a product of distinct irreducible polynomials over $L_R$ if and only if f and f ´ are relatively prime.*

Let A be a n × n matrix with entries from the DSm refined label field $L_R$. Let $A = (L_{a_{ij}})_{n \times n}$ where $L_{a_{ij}} \in L_R$, a characteristic value of A in $L_R$ is a refined label $L_c$ in $L_R$ such that the matrix $(A - L_c I)$ is singular non invertible.

Since $L_c$ is a characteristic value of $A = (L_{a_{ij}})$ if and only if det $(A - L_c I) = 0$ or equivalently if and only if det $(L_c I - A) = 0$, we form the matrix $(xI - A)$ with polynomial entries, and consider the polynomial $f(x) = \det(xI - A)$. Clearly the characteristic value of A in $L_R$ are just the scalars $L_c$ in $L_R$ such that $f(L_c) = 0$.

For this reason f is called the characteristic polynomial of $A = (L_{a_{ij}})$. It is important to note that f is a monic polynomial which has degree exactly n. This is easily seen from the formula for the determinant of a matrix in terms of its entries.

We will illustrate this by some examples.

*Example 2.47:* Let
$$A = \begin{bmatrix} L_a & L_b \\ L_c & L_d \end{bmatrix}$$

where $L_a, L_b, L_c, L_d \in L_R$. The characteristic polynomial for A is det $(xI - A)$

$$\begin{aligned}
&= \begin{vmatrix} xL_{m+1} - L_a & -L_b \\ -L_c & xL_{m+1} - L_d \end{vmatrix} \\
&= (xL_{m+1} - L_a)(xL_{m+1} - L_d) - L_b L_c \\
&= x^2 L_{m+1} - xL_a L_{m+1} - x L_d L_{m+1} + L_a L_d - L_b L_c \\
&= x^2 L_{m+1} - x L_a - x L_d + L_{ad/m+1} - L_{bc/m+1} \\
&= L_{m+1} x^2 - x (L_{a+d}) + L_{(ad-b)/m+1} \\
&= 0.
\end{aligned}$$

We can solve for x

$$x = \frac{L_{a+d} \pm \sqrt{L_{(a+d)^2} - 4L_{a+d} L_{ad-bc/m+1}}}{2L_{m+1}}$$

$$= \frac{L_{a+d} \pm \sqrt{L_{a^2+d^2+2ad/m+1} - 4L_{(a+d)(ad-bc)/(m+1)^2}}}{2L_{m+1}}$$

knowing the values of a, b, c and d one can solve for x in terms of elements from $L_R$.

As in case of usual linear operators we can in case linear operators in DSm vector space V, can say a DSm linear operator on V is diagonalizable if there is a basis for V each vector of which is a characteristic vector of T.

Now V is a DSm vector space over R and $L_R \cong R$ so we can with appropriate modifications arrive at the DSm characteristic values related with a DSm linear operator on V.

Suppose that

$$T\alpha = L_c^\alpha = \frac{c}{m+1} \alpha.$$

If f is any polynomial with coefficients from $L_R$ then $f(T)\alpha = f(L_c)\alpha$.

Let T be a linear operator on the DSm finite dimensional space V. If $L_{c_1}, L_{c_2}, ..., L_{c_k}$ be the distinct characteristic values of T and let $W_i$ be the DSm space of DSm characteristic vectors associated with the DSm characteristic value $c_i$. If $W = W_1 + \ldots + W_k$, then dim $W$ = dim $W_1$ + dim $W_2$ + … + dim $W_k$.

Let T be a DSm linear operator on a finite dimensional DSm vector space over the field R. The minimal polynomial for T is the unique monic generator of the ideal of polynomials over R which annihilate T.

**THEOREM 2.33:** *(Cayley-Hamilton) Let T be a DSm linear operator on a DSm finite dimensional vector space V. If f is the characteristic polynomial for T then f (T) = 0; in other words the minimal polynomial divides the characteristic polynomial for T.*

We define DSm invariant subspaces or invariant DSm subspaces. Let V be a DSm vector space over the reals R and T be a DSm linear operator on V. If W is a DSm subspace of V, we say W is DSm invariant under T if each vector $\alpha$ in W, the vector $T\alpha$ is in W that is if T (W) is contained in W.

Let W be a DSm invariant subspace for the DSm linear operator for T and let $\alpha$ be an element in V. The T-conductor of $\alpha$ into W is a set $S_T(\alpha; W)$ which consists of all polynomials g (over the scalar field) such that $g(T)\alpha$ is in W. The unique monic generator of the DSm ideal $S(\alpha; W)$ is also called the DSm T-conductor of $\alpha$ into W (the T-annihilator in case W = {0}).

In an analogous way we have for any DSm finite dimensional vector space over the real field R. Let T be a DSm linear operator on V such that the minimal polynomial for T is a product of linear factors.

$$p = (x - c_1)^{r_1}(x - c_2)^{r_2}...(x - c_k)^{r_k} \ ; c_i \in R; 1 \le i \le k.$$

Let W be a proper DSm subspace of V which is invariant under T. There exists a vector $\alpha$ in V such that
    (a) $\alpha$ is not in W
    (b) $(T - cI)\alpha$ is in W for some DSm characteristic value c of the DSm operator T.

Several results in this direction can be studied for DSm vector spaces over R with simple appropriate operations.

**Chapter Three**

# SPECIAL DSM VECTOR SPACES

In this chapter we define new types of DSm linear algebras and DSm vector spaces over sets, semigroups etc.

We illustrate them with examples.

Let $L_R = \{L_a \mid a \in R\}$ be the collection of all refined labels. $L_R$ is the DSm real field of refined labels. Also $L_R$ is isomorphic with the field of reals as fields.

$$L_r = \frac{r}{m+1}$$ that is for every $L_r$ in $L_R$ there exists a unique r in R such that $r = \frac{a}{m+1}$. For every r in R there exists $L_a$ in $L_R$ such that $L_a = L_{r(m+1)}$ such that $r = L_a$.

**DEFINITION 3.1**: *Let S be a subset of reals R. V be a subset labels (the subset can be from ordinary labels or from the set of refined labels). We say V is a refined label set vector space over the set S or DSm set vector space of refined labels over the set S (V is a ordinary labels set vector space over the set S) if for all v $\in$ V and for all s $\in$ S; vs and sv $\in$ V.*

We give examples of this.

***Example 3.1:*** Let $V = \{L_0, L_{m+1}, L_a, L_b\}$ be a set of refined labels. Let $S = \{0, 1\}$ be a set. V is a set refined label space over the set S (or set refined label space over the set S).

***Example 3.2:*** Let $V = \{0, L_1, L_2, L_5, L_6\}$ be a set of ordinary labels. Suppose $S = \{0, 1\}$ be a subset of R. V is a set ordinary label space over the set S.

***Example 3.3:*** Let

$$V = \left\{ \begin{pmatrix} L_a & L_b \\ L_c & L_d \end{pmatrix}, (L_a \ L_b), (L_c \ L_d \ L_e) \middle| L_a, L_b, L_c, L_d, L_e \in L_R \right\}$$

be a set of refined labels. V is a set vector space of refined labels over the set $S = \{0,1\}$.

***Example 3.4:*** Let $S = \{0, 1, 2, \ldots, 25\}$ be a set.

$$V = \left\{ \begin{pmatrix} L_a \\ L_b \\ L_c \\ L_d \end{pmatrix}, (L_a \ L_b \ L_c), \sum_{i=0}^{5} L_{a_i} x^i \middle| L_a, L_b, L_c, L_d, L_{a_i} \in L_R ; 0 \le i \le 5 \right\}$$

V is a set vector space of refined labels over the set S.

***Example 3.5:*** Let

$$V = \left\{ \begin{pmatrix} L_{a_1} & L_{a_2} & L_{a_3} \\ L_{a_4} & L_{a_5} & L_{a_6} \\ L_{a_7} & L_{a_8} & L_{a_9} \\ L_{a_{10}} & L_{a_{11}} & L_{a_{12}} \\ L_{a_{13}} & L_{a_{14}} & L_{a_{16}} \end{pmatrix}, (L_{a_1} \ L_{a_2} \ L_{a_3} \ L_{a_4}) \middle| L_{a_i} \in L_R ; 1 \le i \le 16 \right\}$$

be a set vector space of refined labels over the set $S = \{0, 12, \ldots, \infty\}$.

*Example 3.6:* Let

$$V = \left\{ \sum_{i=1}^{8} L_{a_i} x^i, \begin{pmatrix} L_{a_1} & L_{a_2} \\ L_{a_3} & L_{a_4} \\ L_{a_5} & L_{a_6} \end{pmatrix}, (L_{a_1}\ L_{a_2}\ L_{a_3}\ L_{a_4}\ L_{a_5}) \,\middle|\, L_{a_i} \in L_R; 0 \leq i \leq 8 \right\}$$

be a set vector space of refined labels over the set $S = Z^+ \cup \{0\}$.

*Example 3.7:* Let

$$V = \left\{ \begin{bmatrix} L_{a_1} & L_{a_2} & L_{a_3} \\ L_{a_4} & L_{a_5} & L_{a_6} \end{bmatrix}, (L_{a_1}\ L_{a_2}\ ,...,\ L_{a_8}) \,\middle|\, L_{a_i} \in L_R; 1 \leq i \leq 8 \right\}$$

be a set vector space of refined labels over the set $S = Q^+ \cup \{0\}$.

*Example 3.8:* Let

$$V = \left\{ \begin{bmatrix} L_{a_1} & L_{a_2} \\ L_{a_3} & L_{a_4} \\ L_{a_5} & L_{a_6} \\ L_{a_7} & L_{a_8} \end{bmatrix}, \begin{pmatrix} L_{a_1} & L_{a_2} \\ L_{a_3} & L_{a_4} \end{pmatrix}, (L_{a_1}, L_{a_2}, L_{a_3}, L_{a_4}, L_{a_5}, L_{a_6}) \,\middle|\, \begin{array}{l} L_{a_i} \in L_R; \\ 1 \leq i \leq 8 \end{array} \right\}$$

be a set vector space of refined labels over the set $S = 3Z^+ \cup \{0\}$.

*Example 3.9:* Let

$$V = \left\{ \begin{bmatrix} L_{a_1} \\ L_{a_2} \\ L_{a_3} \\ L_{a_4} \\ L_{a_5} \\ L_{a_6} \\ L_{a_7} \end{bmatrix}, (L_{a_1}, L_{a_2}, L_{a_3}, L_{a_4}) \,\middle|\, L_{a_i} \in L_R; 1 \leq i \leq 7 \right\}$$

be a set vector space of refined labels over the set $S = 15Z^+ \cup \{0\}$.

*Example 3.10:* Let

$$V = \left\{ (L_{a_1}\ L_{a_2}\ L_{a_3}), \begin{bmatrix} L_{a_1} & L_{a_2} & L_{a_3} & L_{a_4} \\ L_{a_5} & L_{a_6} & L_{a_7} & L_{a_8} \\ L_{a_9} & 0 & L_{a_{10}} & L_{a_{11}} \\ 0 & L_{a_{12}} & 0 & L_{a_{13}} \end{bmatrix} \middle| L_{a_i} \in L_R; 1 \leq i \leq 13 \right\}$$

be a set of refined labels over the set $S = 9Z^+ \cup \{0\}$.

*Example 3.11:* Let

$$V = \left\{ \sum_{i=0}^{19} L_{a_i} x^i, \begin{pmatrix} L_{a_1} & L_{a_2} & L_{a_3} & L_{a_4} & L_{a_5} \\ L_{a_6} & L_{a_7} & L_{a_8} & L_{a_9} & L_{a_{10}} \end{pmatrix} \middle| L_{a_i} \in L_R; 0 \leq i \leq 19 \right\}$$

be a set vector space of refined labels over the set $S = 8Z^+ \cup \{0\}$.

*Example 3.12:* Let

$$V = \left\{ \begin{bmatrix} L_{a_1} & L_{a_2} & L_{a_3} \\ L_{a_4} & L_{a_5} & L_{a_6} \\ L_{a_7} & L_{a_8} & L_{a_9} \end{bmatrix}, \sum_{i=0}^{8} L_{a_i} x^i, \begin{pmatrix} L_{a_1} \\ L_{a_2} \\ L_{a_3} \\ \vdots \\ L_{a_{14}} \end{pmatrix} \middle| L_{a_i} \in L_R; 0 \leq i \leq 14 \right\}$$

be a set vector space of refined labels over the set $S = 17Z^+ \cup \{0\}$.

*Example 3.13:* Let

$$P = \left\{ \begin{pmatrix} L_{a_1} & L_{a_2} & L_{a_3} \\ L_{a_4} & L_{a_5} & L_{a_6} \end{pmatrix}, \begin{pmatrix} L_{a_1} \\ L_{a_2} \\ L_{a_3} \\ L_{a_4} \\ L_{a_5} \end{pmatrix}, \begin{bmatrix} L_{a_1} & L_{a_2} & L_{a_3} \\ L_{a_4} & L_{a_5} & L_{a_6} \\ L_{a_7} & L_{a_8} & L_{a_9} \end{bmatrix} \middle| L_{a_i} \in L_R; 0 \leq i \leq 9 \right\}$$

be a set vector space of refined labels over the set $S = R^+ \cup \{0\}$.

*Example 3.14:* Let

$$V = \left\{ \begin{bmatrix} L_{a_1} & L_{a_2} \\ L_{a_3} & L_{a_4} \\ L_{a_5} & L_{a_6} \\ L_{a_7} & L_{a_8} \\ L_{a_9} & L_{a_{10}} \\ L_{a_{11}} & L_{a_{12}} \end{bmatrix}, \begin{bmatrix} L_{a_1} & L_{a_2} & L_{a_3} \\ L_{a_4} & L_{a_5} & L_{a_6} \\ L_{a_7} & L_{a_8} & L_{a_9} \end{bmatrix} \middle| L_{a_i} \in L_R; 1 \leq i \leq 12 \right\}$$

be a set vector space of refined labels over the set $S = Q^+ \cup \{0\}$.

*Example 3.15:* Let

$$V = \left\{ \begin{bmatrix} L_{a_1} & L_{a_2} & L_{a_3} \\ L_{a_4} & L_{a_5} & L_{a_6} \\ L_{a_7} & L_{a_8} & L_{a_9} \\ L_{a_{10}} & L_{a_{11}} & L_{a_{12}} \\ L_{a_{13}} & L_{a_{14}} & L_{a_{15}} \end{bmatrix}, \begin{bmatrix} L_{a_1} & L_{a_2} & L_{a_3} & L_{a_4} \\ L_{a_5} & L_{a_6} & L_{a_7} & L_{a_8} \end{bmatrix}, \begin{bmatrix} L_{a_1} \\ L_{a_2} \\ L_{a_3} \\ L_{a_4} \end{bmatrix} \middle| \begin{array}{l} L_{a_i} \in L_R; \\ 1 \leq i \leq 15 \end{array} \right\}$$

be a set vector space of refined labels over set $S = 5Z^+ \cup \{0\}$.

*Example 3.16:* Let

$$V = \left\{ \begin{bmatrix} L_{a_1} & L_{a_2} & L_{a_3} \\ L_{a_4} & L_{a_5} & L_{a_6} \\ \vdots & \vdots & \vdots \\ L_{a_{26}} & L_{a_{27}} & L_{a_{28}} \end{bmatrix}, \sum_{i=0}^{20} L_{a_i} x^i, \begin{bmatrix} L_{a_1} & L_{a_2} & L_{a_3} \\ L_{a_4} & L_{a_5} & L_{a_6} \\ L_{a_7} & L_{a_8} & L_{a_9} \end{bmatrix} \middle| L_{a_i} \in L_R; 0 \leq i \leq 27 \right\}$$

be a set vector space of refined labels over the set $S = R^+ \cup \{0\}$.

Now having seen examples of set vector space of refined labels over the set S, we now proceed onto define substructures in them.

*Example 3.17:* Let

$$V = \left\{ \sum_{i=0}^{8} L_{a_i} x^i, \begin{bmatrix} L_{a_1} \\ L_{a_2} \\ L_{a_3} \\ L_{a_4} \end{bmatrix}, \begin{bmatrix} L_{a_1} & L_{a_2} & L_{a_3} & L_{a_4} \\ L_{a_5} & L_{a_6} & L_{a_7} & L_{a_8} \end{bmatrix} \middle| L_{a_i} \in L_R; 0 \le i \le 8 \right\}$$

be a set vector space refined labels over the set $S = 5Z^+ \cup \{0\}$.

Consider

$$W = \left\{ \begin{bmatrix} L_{a_1} \\ L_{a_2} \\ L_{a_3} \\ L_{a_4} \end{bmatrix} \middle| L_{a_i} \in L_R; 1 \le i \le 4 \right\} \subseteq V,$$

W is a set vector space of refined labels over set $S = 5Z^+ \cup \{0\}$.

**DEFINITION 3.2**: *Let V be a DSm set vector space of refined labels over the set S. Suppose $W \subseteq V$; if W is a set vector space of refined labels over the set S then W is a DSm set vector subspace of refined labels over the set S of V.*

We will illustrate this situation by some examples.

*Example 3.18*: Let

$$V = \left\{ \sum_{i=0}^{9} L_{a_i} x^i, \begin{bmatrix} L_{a_1} \\ L_{a_2} \\ \vdots \\ L_{a_8} \end{bmatrix}, \begin{bmatrix} L_{a_1} & L_{a_2} & L_{a_3} & L_{a_4} \\ L_{a_5} & L_{a_6} & L_{a_7} & L_{a_8} \\ L_{a_9} & L_{a_{10}} & L_{a_{11}} & L_{a_{12}} \end{bmatrix} \middle| L_{a_i} \in L_R; 0 \le i \le 12 \right\}$$

be a set vector space refined labels over the set $S = 3Z^+ \cup \{0\}$.

Consider

$$W = \left\{ \sum_{i=0}^{9} L_{a_i} x^i, \begin{bmatrix} L_{a_1} \\ L_{a_2} \\ \vdots \\ L_{a_8} \end{bmatrix} \middle| L_{a_i} \in L_R; 0 \le i \le 9 \right\} \subseteq V$$

a set vector subspace of refined labels over the set $S = 3Z^+ \cup \{0\}$.

*Example 3.19*: Let

$$V = \left\{ \sum_{i=0}^{12} L_{a_i} x^i, \begin{bmatrix} L_{a_1} & L_{a_2} & L_{a_3} \\ L_{a_4} & L_{a_5} & L_{a_6} \\ L_{a_7} & L_{a_8} & L_{a_9} \end{bmatrix} \middle| L_{a_i} \in L_R; 0 \le i \le 12 \right\}$$

be a DSm set vector space of refined labels over the set $S = 3Z^+ \cup \{0\}$.

Take

$$W = \left\{ \sum_{i=0}^{7} L_{a_i} x^i, \begin{bmatrix} L_{a_1} & L_{a_2} & L_{a_3} \\ 0 & L_{a_4} & 0 \\ 0 & 0 & L_{a_5} \end{bmatrix} \middle| L_{a_i} \in L_R; 0 \le i \le 7 \right\} \subseteq V,$$

W is a set vector subspace of refined labels over S of V.

*Example 3.20:* Let

$$V = \left\{ \begin{bmatrix} L_{a_1} \\ L_{a_2} \\ L_{a_3} \\ L_{a_4} \end{bmatrix}, (L_{a_1}, L_{a_2}, L_{a_3}, ..., L_{a_9}) \middle| L_{a_i} \in L_R; 1 \le i \le 9 \right\}$$

be a set vector space of refined labels over the set $S = Z^+ \cup \{0\}$.

Consider

$$W = \left\{ \begin{bmatrix} L_{a_1} \\ 0 \\ L_{a_2} \\ 0 \end{bmatrix}, (L_{a_1}, 0, L_{a_2}, 0, L_{a_3}, 0, L_{a_4}, 0, L_{a_5}) \middle| L_{a_i} \in L_R; 1 \leq i \leq 5 \right\}$$

$\subseteq V$ be a set vector subspace of refined labels of V over the set $S = Z^+ \cup \{0\}$.

*Example 3.21:* Let

$$V = \left\{ \begin{pmatrix} L_{a_1} & L_{a_2} \\ L_{a_3} & L_{a_4} \end{pmatrix}, \begin{bmatrix} L_{a_1} & L_{a_2} & L_{a_3} \\ L_{a_4} & 0 & 0 \\ L_{a_5} & L_{a_6} & L_{a_7} \\ 0 & 0 & L_{a_8} \\ 0 & L_{a_9} & 0 \end{bmatrix} \middle| L_{a_i} \in L_R; 1 \leq i \leq 9 \right\}$$

be a set vector space of refined labels over the set $S = 5Z^+ \cup \{0\}$. Consider

$$W = \left\{ \begin{pmatrix} 0 & L_{a_1} \\ L_{a_2} & 0 \end{pmatrix}, \begin{bmatrix} L_{a_1} & L_{a_2} & L_{a_3} \\ 0 & 0 & 0 \\ L_{a_4} & L_{a_5} & L_{a_6} \\ 0 & 0 & 0 \\ 0 & L_{a_7} & 0 \end{bmatrix} \middle| L_{a_i} \in L_R; 1 \leq i \leq 7 \right\} \subseteq V$$

be a set vector subspace of refined labels of V over the set $S = 5Z^+ \cup \{0\}$.

*Example 3.22:* Let

$$V = \left\{ \sum_{i=0}^{3} L_{a_i} x^i, \begin{bmatrix} L_{a_1} & L_{a_2} & L_{a_3} & L_{a_4} & L_{a_5} \\ L_{a_6} & L_{a_7} & L_{a_8} & L_{a_9} & L_{a_{10}} \\ L_{a_{12}} & L_{a_{13}} & L_{a_{14}} & L_{a_{15}} & L_{a_{16}} \end{bmatrix} \middle| L_{a_i} \in L_R; 0 \leq i \leq 16 \right\}$$

be a set vector space of refined labels over set $S = 20Z^+ \cup \{0\}$.

Consider

$$W = \left\{ \sum_{i=0}^{3} L_{a_i} x^i, \begin{bmatrix} L_{a_1} & L_{a_2} & L_{a_3} & L_{a_4} & L_{a_5} \\ 0 & 0 & 0 & 0 & 0 \\ L_{a_6} & L_{a_7} & L_{a_8} & L_{a_9} & L_{a_{10}} \end{bmatrix} \middle| L_{a_i} \in L_R; 0 \le i \le 16 \right\} \subseteq V$$

be a set vector subspace of refined labels of V over the set S = $20Z^+ \cup \{0\}$.

We can have several DSm set vector subspaces of refined labels over the set S.

We say for any DSm set vector space of refined labels V over the set S we say V is a direct sum of DSm set vector subspace of refined labels $W_1, W_2, \ldots, W_k$ if $V = W_1 + W_2 + \ldots + W_k$ and $W_i \cap W_j = 0$ or $\phi$ if $i \ne j$. If on the other hand $W_1, W_2, \ldots, W_k$ are DSm vector subspace of refined labels of V if $V = W_1 + \ldots + W_k$ but $W_i \cap W_j \ne 0$ if $i \ne j$ then we say V is just a sum of the DSm set vector subspaces $W_1, W_2, \ldots, W_k$.

We will illustrate this situation by some examples.

***Example 3.23:*** Let

$$V = \left\{ \begin{pmatrix} L_{a_1} & L_{a_2} \\ L_{a_3} & L_{a_4} \end{pmatrix}, \sum_{i=0}^{8} L_{a_i} x^i, (L_{a_1}, L_{a_2}, \ldots, L_{a_8}) \middle| L_{a_i} \in L_R; 0 \le i \le 8 \right\}$$

be a DSm set vector space over the set $3Z^+ \cup \{0\}$.

Consider

$$W_1 = \left\{ \begin{pmatrix} L_{a_1} & L_{a_2} \\ L_{a_3} & L_{a_4} \end{pmatrix} \middle| L_{a_i} \in L_R; 1 \le i \le 4 \right\} \subseteq V,$$

$$W_2 = \left\{ \sum_{i=0}^{\infty} L_{a_i} x^i \middle| L_{a_i} \in L_R; 0 \le i \le 8 \right\} \subseteq V$$

and $W_3 = \left\{ (L_{a_1}, L_{a_2}, L_{a_3}, \ldots, L_{a_8}) \middle| L_{a_i} \in L_R; 1 \le i \le 8 \right\} \subseteq V$.

We see $W_1$, $W_2$ and $W_3$ are set vector subspaces of refined labels of V over the set $S = 3Z^+ \cup \{0\}$. Further $V = W_1 + W_2 + W_3$ with $W_i \cap W_j = \phi$, if $i \ne j$; $1 \le i, j \le 3$.

*Example 3.24:* Let

$$V = \left\{ \sum_{i=0}^{4} L_{a_i} x^i, \begin{bmatrix} L_{a_1} & L_{a_2} & L_{a_3} & L_{a_4} \\ L_{a_5} & L_{a_6} & L_{a_7} & L_{a_8} \end{bmatrix}, \begin{bmatrix} L_{a_1} & L_{a_2} & L_{a_3} \\ L_{a_4} & L_{a_5} & L_{a_6} \\ L_{a_7} & L_{a_8} & L_{a_9} \\ L_{a_{10}} & L_{a_{11}} & L_{a_{12}} \\ L_{a_{13}} & L_{a_{14}} & L_{a_{15}} \end{bmatrix} \middle| \begin{array}{l} L_{a_i} \in L_R; \\ 0 \leq i \leq 15 \end{array} \right\}$$

be a DSm set vector space of refined labels over the set $S = Q^+ \cup \{0\}$. Consider

$$W_1 = \left\{ \sum_{i=0}^{3} L_{a_i} x^i, \begin{bmatrix} L_{a_1} & L_{a_2} & L_{a_3} & L_{a_4} \\ L_{a_5} & L_{a_6} & L_{a_7} & L_{a_8} \end{bmatrix} \middle| L_{a_i} \in L_R; 0 \leq i \leq 8 \right\} \subseteq V,$$

$$W_2 = \left\{ \sum_{i=0}^{2} L_{a_i} x^i, \begin{bmatrix} L_{a_1} & L_{a_2} & L_{a_3} \\ \vdots & \vdots & \vdots \\ L_{a_{13}} & L_{a_{14}} & L_{a_{15}} \end{bmatrix} \middle| L_{a_i} \in L_R; 0 \leq i \leq 15 \right\} \subseteq V$$

and

$$W_3 = \left\{ \sum_{i=0}^{4} L_{a_i} x^i, \begin{bmatrix} L_{a_1} & L_{a_2} & L_{a_3} & L_{a_4} \\ L_{a_5} & L_{a_6} & L_{a_7} & L_{a_8} \end{bmatrix} \middle| L_{a_i} \in L_R; 0 \leq i \leq 8 \right\} \subseteq V,$$

be DSm subspaces of V of refined labels. Clearly $V = W_1 + W_2 + W_3$ but $W_i \cap W_j \neq \phi$ or $(0)$. Hence V is only a sum of $W_1$, $W_2$ and $W_3$ and not a direct sum.

*Example 3.25:* Let

$$V = \left\{ \begin{bmatrix} L_{a_1} & L_{a_2} & L_{a_3} \\ L_{a_4} & L_{a_5} & L_{a_6} \\ L_{a_7} & L_{a_8} & L_{a_9} \end{bmatrix}, \begin{bmatrix} L_{a_1} & L_{a_2} \\ L_{a_3} & L_{a_4} \\ L_{a_5} & L_{a_6} \\ L_{a_7} & L_{a_8} \\ L_{a_9} & L_{a_{10}} \\ L_{a_{11}} & L_{a_{12}} \end{bmatrix}, \begin{pmatrix} L_{a_1} & L_{a_2} & L_{a_3} & L_{a_4} & L_{a_5} & L_{a_6} \\ L_{a_7} & L_{a_8} & L_{a_9} & L_{a_{10}} & L_{a_{11}} & L_{a_{12}} \end{pmatrix} \middle| \begin{array}{l} L_{a_i} \in L_R; \\ 1 \leq i \leq 12 \end{array} \right\}$$

be a DSm vector space of refined labels over the set $S = 13Z^+ \cup \{0\}$.

Consider

$$W_1 = \left\{ \begin{bmatrix} L_{a_1} & L_{a_2} & L_{a_3} \\ L_{a_4} & L_{a_5} & L_{a_6} \\ L_{a_7} & L_{a_8} & L_{a_9} \end{bmatrix} \middle| L_{a_i} \in L_R; 1 \leq i \leq 9 \right\} \subseteq V,$$

$$W_2 = \left\{ \begin{pmatrix} L_{a_1} & L_{a_2} & L_{a_3} & L_{a_4} & L_{a_5} & L_{a_6} \\ L_{a_7} & L_{a_8} & L_{a_9} & L_{a_{10}} & L_{a_{11}} & L_{a_{12}} \end{pmatrix}, \begin{bmatrix} L_{a_1} & 0 & 0 \\ 0 & L_{a_2} & 0 \\ 0 & 0 & L_{a_3} \end{bmatrix} \middle| \begin{matrix} L_{a_i} \in L_R; \\ 1 \leq i \leq 12 \end{matrix} \right\} \subseteq V$$

and

$$W_3 = \left\{ \begin{bmatrix} L_{a_1} & L_{a_2} & L_{a_3} \\ 0 & 0 & 0 \\ L_{a_4} & 0 & L_{a_5} \end{bmatrix}, \begin{bmatrix} L_{a_1} & L_{a_2} \\ L_{a_3} & L_{a_4} \\ L_{a_5} & L_{a_6} \\ L_{a_7} & L_{a_8} \\ L_{a_9} & L_{a_{10}} \\ L_{a_{11}} & L_{a_{12}} \end{bmatrix} \middle| L_{a_i} \in L_R; 1 \leq i \leq 12 \right\} \subseteq V,$$

be DSm vector subspaces refined labels over the set $S = 13Z^+ \cup \{0\}$. Clearly $V = W_1 + W_2 + W_3$ and $W_i \cap W_j \neq (0)$ or $\phi$ if $i \neq j$; $1 \leq i, j \leq 3$.

Now having seen examples of DSm subspaces of V and direct sum and sum of subspaces of V.

**DEFINITION 3.3**: *Let V be a DSm set vector space of refined labels over the set S. We say a proper subset $B \subseteq V$ is said to be a linearly independent set if $x, y \in B$ then $x \neq sy$ or $y \neq s'x$ for any s and s' in S. If the set B is not linearly independent then we say B is a linearly dependent set.*

We will illustrate this situation by an example.

*Example 3.26:* Let

$$V = \left\{ \begin{pmatrix} L_{a_1} & L_{a_2} \\ L_{a_3} & L_{a_4} \end{pmatrix}, (L_{a_1}, L_{a_2}, L_{a_3}), \begin{bmatrix} L_{a_1} \\ L_{a_2} \\ L_{a_3} \\ L_{a_4} \end{bmatrix}, \sum_{i=0}^{6} L_{a_i} x^i \,\middle|\, \begin{array}{l} L_{a_i} \in L_R; \\ 0 \le i \le 4 \end{array} \right\}$$

be a DSm set vector space of refined labels over the set $S = 3Z^+ \cup \{0\}$. We see V has linearly independent sets.

For take

$$B = \left\{ \begin{bmatrix} L_{a_1} \\ L_{a_2} \\ L_{a_3} \\ L_{a_4} \end{bmatrix}, (L_{a_1}, L_{a_2}, L_{a_3}), \begin{pmatrix} L_{a_1} & L_{a_2} \\ L_{a_3} & L_{a_4} \end{pmatrix} \,\middle|\, L_{a_i} \in L_R; 1 \le i \le 4 \right\} \subseteq V.$$

Clearly $(L_{a_1}, L_{a_2}, L_{a_3}) \ne s \begin{pmatrix} L_{a_1} & L_{a_2} \\ L_{a_3} & L_{a_4} \end{pmatrix}$ for any $s \in S$. We can give several such subsets in V which are linearly independent.

**DEFINITION 3.4**: *Let V and W be any two DSm set vector spaces of refined labels over the same set S. A map T from V to W is said to be a set linear transformation if*

$$T(v) = \omega$$
$$T(sv) = s\omega = sT(v).$$

*for all $v, \omega \in V$ and $s \in S$.*

We will illustrate this situation by some examples.

*Example 3.27:* Let

$$V = \left\{ \sum_{i=0}^{5} L_{a_i} x^i, \begin{pmatrix} L_{a_1} & L_{a_2} \\ L_{a_3} & L_{a_4} \end{pmatrix} \,\middle|\, L_{a_i} \in L_R; 0 \le i \le 5 \right\}$$

and

$$W = \left\{ \begin{pmatrix} L_{a_1} & L_{a_2} & L_{a_3} \\ L_{a_4} & L_{a_5} & L_{a_6} \end{pmatrix}, (L_{a_1} \ L_{a_2} \ L_{a_3} \ L_{a_4}) \middle| L_{a_i} \in L_R ; 1 \leq i \leq 6 \right\}$$

be two DSm set vector spaces over the set $S = R^+ \cup \{0\}$ of refined labels.

Choose $T : V \to W$ a map such that

$$T\left(\sum_{i=0}^{5} L_{a_i} x^i\right) = \begin{pmatrix} L_{a_0} & L_{a_1} & L_{a_2} \\ L_{a_3} & L_{a_4} & L_{a_5} \end{pmatrix}$$

and

$$T\begin{pmatrix} L_{a_1} & L_{a_2} \\ L_{a_3} & L_{a_4} \end{pmatrix} = (L_{a_1} \ L_{a_2} \ L_{a_3} \ L_{a_4});$$

T is easily verified to be a DSm set vector linear transformation of refined labels.

*Example 3.28:* Let

$$V = \left\{ \begin{pmatrix} L_{a_1} & L_{a_2} \\ L_{a_3} & L_{a_4} \\ L_{a_5} & L_{a_6} \\ L_{a_7} & L_{a_8} \\ L_{a_9} & L_{a_{10}} \end{pmatrix}, \begin{pmatrix} L_{a_1} & L_{a_3} & L_{a_4} \\ L_{a_2} & 0 & L_{a_5} \end{pmatrix} \middle| L_{a_i} \in L_R ; 1 \leq i \leq 10 \right\}$$

and

$$W = \left\{ \begin{pmatrix} L_{a_1} & L_{a_2} \\ L_{a_3} & L_{a_4} \end{pmatrix}, (L_{a_1}, L_{a_2}, ..., L_{a_{10}}) \middle| L_{a_i} \in L_R ; 1 \leq i \leq 10 \right\}$$

be two DSm set vector spaces of refined labels over the set $S = Z^+ \cup \{0\}$.

Define a map $T : V \to W$ by

$$T\left(\begin{pmatrix} L_{a_1} & L_{a_3} & L_{a_4} \\ L_{a_2} & 0 & L_{a_5} \end{pmatrix}\right) = \begin{pmatrix} L_{a_1} & L_{a_2} \\ L_{a_4} & L_{a_6} \end{pmatrix}.$$

and

$$T\left(\begin{pmatrix} L_{a_1} & L_{a_2} \\ L_{a_3} & L_{a_4} \\ L_{a_5} & L_{a_6} \\ L_{a_7} & L_{a_8} \\ L_{a_9} & L_{a_{10}} \end{pmatrix}\right) = (L_{a_1}, L_{a_2}, ..., L_{a_{10}})$$

T is a DSm set linear transformation of V into W.

Now we can define T to be set DSm linear operator if V = W.

We will just illustrate this situation by some examples.

*Example 3.29:* Let

$$V = \left\{ \sum_{i=0}^{15} L_{a_i} x^i, \begin{pmatrix} L_{a_1} & L_{a_2} \\ L_{a_3} & L_{a_4} \end{pmatrix}, 4 \times 4 \text{ refined label matrix}, \left(L_{a_1}, L_{a_2}, L_{a_3}, L_{a_4}\right) \right.$$

$\left. L_{a_i} \in L_R; 0 \le i \le 15 \right\}$ be a DSm set vector space of refined labels on the set $S = 3Z^+ \cup \{0\}$.

Define $T : V \to V$ by

$$T\left(\sum_{i=0}^{15} L_{a_i} x^i\right) = \begin{bmatrix} L_{a_1} & L_{a_2} & L_{a_3} & L_{a_4} \\ L_{a_5} & L_{a_6} & L_{a_7} & L_{a_8} \\ L_{a_9} & L_{a_{10}} & L_{a_{11}} & L_{a_{12}} \\ L_{a_{13}} & L_{a_{14}} & L_{a_{15}} & L_{a_{16}} \end{bmatrix},$$

$$T\left(\begin{pmatrix} L_{a_1} & L_{a_2} \\ L_{a_3} & L_{a_4} \end{pmatrix}\right) = \left(L_{a_1}, L_{a_2}, L_{a_3}, L_{a_4}\right),$$

$$T\left(\begin{bmatrix} L_{a_0} & L_{a_2} & L_{a_3} & L_{a_4} \\ L_{a_5} & L_{a_6} & L_{a_7} & L_{a_8} \\ L_{a_9} & L_{a_{10}} & L_{a_{11}} & L_{a_{12}} \\ L_{a_{13}} & L_{a_{14}} & L_{a_{15}} & L_{a_1} \end{bmatrix}\right) = \left(\sum_{i=0}^{15} L_{a_i} x^i\right)$$

and

$$T\left((L_{a_1}, L_{a_2}, L_{a_3}, L_{a_4})\right) = \begin{pmatrix} L_{a_1} & L_{a_2} \\ L_{a_3} & L_{a_4} \end{pmatrix};$$

then T is a DSm set linear operator on the set vector space of refined labels V over the set S. We see as in case of usual vector spaces study the algebraic structure enjoyed by the set of all set linear transformations from V to W. Likewise the reader is expected to study the algebraic structure enjoyed by the set of all set linear operators of V to V of the DSm set vector space of refined labels over the set S.

Now having seen examples of them we can proceed onto define the notion of DSm set linear algebraic of refined labels over the set S or set linear algebra of refined labels over the set S.

**DEFINITION 3.5**: *Let V be a DSm set vector space of refined labels over the set S. If V is closed with respect to the operation of addition and if for s ∈ S and v, w ∈ V we have s (v + w) = sv + sw then we define V to be a DSm set linear algebra of refined labels over the set S.*

We will first illustrate this situation by some examples.

*Example 3.30:* Let

$$V = \left\{ \begin{bmatrix} L_{a_1} & L_{a_2} \\ L_{a_3} & L_{a_4} \\ L_{a_5} & L_{a_6} \end{bmatrix} \middle| L_{a_i} \in L_R; 1 \leq i \leq 6 \right\}$$

be a DSm set linear algebra of refined labels over the set $S = Z^+ \cup \{0\}$.

*Example 3.31:* Let
$$V = \left\{ \sum_{i=0}^{8} L_{a_i} x^i \,\middle|\, L_{a_i} \in L_R ; 0 \le i \le 8 \right\}$$
be a set linear algebra of refined labels over the set $S = 5Z^+ \cup 3Z^+ \cup \{0\}$.

*Example 3.32:* Let
$$W = \left\{ \sum_{i=0}^{\infty} L_{a_i} x^i \,\middle|\, L_{a_i} \in L_R ; 0 \le i \le \infty \right\}$$
be a set linear algebra of refined labels over the set $S = 7Z^+ \cup 3Z^+ \cup 8Z^+ \cup \{0\}$.

Now we have the following theorem.

**THEOREM 3.1**: *Let V be a set linear algebra of refined labels over the set S. V is a set vector space of refined labels over the set S. However if V is a set vector space of refined labels V need not in general be a set linear algebra of refined labels over S.*

The proof is direct and hence is left as an exercise to the reader.
   Now we will give examples of set linear subalgebras of refined labels.

*Example 3.33:* Let
$$V = \left\{ \sum_{i=0}^{20} L_{a_i} x^i \,\middle|\, L_{a_i} \in L_R \right\}$$
be a set linear algebra of refined labels over the set $S = 3Z^+ \cup \{0\}$. Let
$$W = \left\{ \sum_{i=0}^{10} L_{a_i} x^i \,\middle|\, L_{a_i} \in L_R \right\} \subseteq V;$$
W is a set linear subalgebra of V of refined labels over R.

*Example 3.34:* Let

$$V = \left\{ \begin{pmatrix} L_{a_1} & L_{a_2} \\ L_{a_3} & L_{a_4} \\ L_{a_5} & L_{a_6} \\ L_{a_7} & L_{a_8} \\ L_{a_9} & L_{a_{10}} \end{pmatrix} \middle| L_{a_i} \in L_R; 1 \le i \le 10 \right\}$$

be a set linear algebra of refined labels over the set $S = Q^+ \cup \{0\}$. Choose

$$W = \left\{ \begin{pmatrix} L_{a_1} & 0 \\ 0 & L_{a_2} \\ 0 & L_{a_3} \\ L_{a_4} & 0 \\ L_{a_5} & 0 \end{pmatrix} \middle| L_{a_i} \in L_R; 1 \le i \le 5 \right\} \subseteq V.$$

W is a set linear subalgebra of refined labels over the set $S = Q^+ \cup \{0\}$ of V.

*Example 3.35:* Let

$$M = \left\{ \begin{pmatrix} L_{a_1} & L_{a_2} & L_{a_3} \\ L_{a_4} & L_{a_5} & L_{a_6} \\ L_{a_7} & L_{a_8} & L_{a_9} \\ L_{a_{10}} & L_{a_{11}} & L_{a_{12}} \end{pmatrix} \middle| L_{a_i} \in L_R; 1 \le i \le 12 \right\}$$

be a set linear algebra of refined labels over $S = Q^+ \cup \{0\}$.

Choose

$$W = \left\{ \begin{pmatrix} L_{a_1} & 0 & L_{a_2} \\ 0 & L_{a_3} & 0 \\ L_{a_4} & 0 & L_{a_5} \\ 0 & L_{a_6} & 0 \end{pmatrix} \middle| L_{a_i} \in L_R; 1 \le i \le 6 \right\} \subseteq M.$$

W is a set linear subalgebra of refined labels over the set $S = Q^+ \cup \{0\}$.

Now having seen examples of set linear subalgebras of refined labels and set vector subspaces of refined labels we now proceed onto define subset linear subalgebras of a set linear algebra of refined labels and subset vector subspaces of set vector spaces of refined labels over a subset of the set over which these structures are defined.

**DEFINITION 3.6**: *Let V be a DSm set linear algebra of refined labels over the set S. Let $W \subseteq V$ be a proper subset of V and $T \subseteq S$ be a proper subset of S. If W is a DSm set linear algebra of refined labels over the set T; then we define W to be a DSm subset linear subalgebra of refined labels over the subset T of S of the set linear algebra V over S.*

We can define analogously define subset vector subspaces of a refined labels of a set vector space of refined labels defined over the set S.

We will illustrate both the situations by some examples.

*Example 3.36:* Let

$$V = \left\{ \sum_{i=0}^{8} L_{a_i} x^i, \begin{bmatrix} L_{a_1} & L_{a_2} & L_{a_3} & L_{a_4} \\ L_{a_5} & L_{a_6} & L_{a_7} & L_{a_8} \\ L_{a_9} & L_{a_{10}} & L_{a_{11}} & L_{a_{12}} \end{bmatrix}, \begin{bmatrix} L_{a_1} \\ L_{a_2} \\ L_{a_3} \\ L_{a_4} \end{bmatrix} \middle| L_{a_i} \in L_R; 0 \le i \le 12 \right\}$$

be a set vector space of refined labels over the set $S = Z^+ \cup \{0\}$.
Choose

$$W = \left\{ \begin{bmatrix} L_{a_1} \\ 0 \\ L_{a_2} \\ 0 \end{bmatrix}, \sum_{i=0}^{5} L_{a_i} x^i \middle| L_{a_i} \in L_R; 0 \le i \le 5 \right\} \subseteq V$$

and $T = 8Z^+ \cup \{0\} \subseteq Z^+ \cup \{0\} = S$; clearly W is a subset vector subspace of V of refined labels over the subset T of S.

*Example 3.37:* Let

$$M = \left\{ \begin{pmatrix} L_{a_1} & L_{a_2} \\ L_{a_3} & L_{a_4} \end{pmatrix}, \begin{pmatrix} L_{a_1} & L_{a_3} & 0 & L_{a_5} & L_{a_7} & 0 \\ L_{a_2} & L_{a_4} & 0 & L_{a_6} & 0 & L_{a_8} \end{pmatrix}, \begin{bmatrix} L_{a_1} & L_{a_2} \\ 0 & L_{a_3} \\ L_{a_4} & 0 \\ 0 & L_{a_5} \\ L_{a_6} & 0 \end{bmatrix} \middle| L_{a_i} \in L_R; 1 \le i \le 8 \right\}$$

be a set vector space of refined labels over the set $S = Q^+ \cup \{0\}$. Let

$$W = \left\{ \begin{pmatrix} L_{a_1} & 0 \\ L_{a_2} & 0 \end{pmatrix}, \begin{bmatrix} L_{a_1} & 0 \\ 0 & 0 \\ L_{a_2} & 0 \\ 0 & 0 \\ L_{a_3} & 0 \end{bmatrix} \middle| L_{a_i} \in L_R; 1 \le i \le 3 \right\} \subseteq V,$$

be a subset vector subspace of refined labels over the subset $T = Z^+ \cup \{0\} \subseteq Q^+ \cup \{0\} = S$.

*Example 3.38:* Let

$$V = \left\{ \sum_{i=0}^{25} L_{a_i} x^i \middle| L_{a_i} \in L_R; 0 \le i \le 25 \right\}$$

be a set linear algebra of refined labels over the set $S = 3Z^+ \cup 5Z^+ \cup 7Z^+ \cup \{0\}$. Let

$$W = \left\{ \sum_{i=0}^{10} L_{a_i} x^i \middle| L_{a_i} \in L_R; 0 \le i \le 10 \right\} \subseteq V;$$

W is a subset linear subalgebra of V over the subset $T = 3Z^+ \cup 7Z^+ \cup \{0\} \subseteq S$ of refined labels.

*Example 3.39:* Let
$$V = \left\{ \begin{pmatrix} L_{a_1} & L_{a_2} & L_{a_3} \\ 0 & L_{a_4} & 0 \\ L_{a_5} & 0 & L_{a_6} \\ 0 & 0 & 0 \\ L_{a_7} & L_{a_8} & L_{a_9} \\ 0 & L_{a_{10}} & 0 \end{pmatrix} \middle| L_{a_i} \in L_R ; 1 \le i \le 10 \right\}$$
be a set linear algebra of refined labels over the set $S = Q^+ \cup \{0\}$.

Consider
$$W = \left\{ \begin{pmatrix} 0 & 0 & 0 \\ 0 & L_{a_1} & 0 \\ L_{a_2} & 0 & L_{a_3} \\ 0 & 0 & 0 \\ 0 & 0 & 0 \\ 0 & 0 & 0 \end{pmatrix} \middle| L_{a_i} \in L_R ; 1 \le i \le 3 \right\} \subseteq V,$$

W is a subset linear subalgebra of refined labels of V over the subset $T = 3Z^+ \cup 5Z^+ \cup 7Z^+ \cup \{0\}$ of $S = Q^+ \cup \{0\}$.

Now for set linear algebras of refined labels we can define the notion of set linear transformations and set linear operator.

We will just give one or two examples interested reader is expect to study them as it is only direct.

*Example 3.40:* Let
$$V = \left\{ \sum_{i=0}^{7} L_{a_i} x^i \middle| L_{a_i} \in L_R ; 0 \le i \le 7 \right\}$$
be a set linear algebra of refined labels over the set $S = Q^+ \cup \{0\}$.

$$W = \left\{ \begin{pmatrix} L_{a_1} & L_{a_2} \\ L_{a_3} & L_{a_4} \\ L_{a_5} & L_{a_6} \\ L_{a_7} & L_{a_8} \end{pmatrix} \middle| L_{a_i} \in L_R; 1 \leq i \leq 8 \right\}$$

be a set linear algebra of refined labels over the set $S = Q^+ \cup \{0\}$.

Define $T : V \to W$ where

$$T(v) = T\left( \sum_{i=0}^{7} L_{a_i} x^i \right) = \begin{pmatrix} L_{a_1} & L_{a_2} \\ L_{a_3} & L_{a_4} \\ L_{a_5} & L_{a_6} \\ L_{a_7} & L_{a_8} \end{pmatrix}$$

where $v \in V$.

It is easily verified T is a set linear transformation of V into W of refined label vector spaces.

*Example 3.41:* Let

$$V = \left\{ \begin{pmatrix} L_a & L_b \\ L_c & L_d \end{pmatrix}, \begin{pmatrix} L_{a_1} & L_{a_2} & L_{a_3} & L_{a_4} & L_{a_5} \\ L_{a_6} & L_{a_7} & L_{a_8} & L_{a_9} & L_{a_{10}} \end{pmatrix} \middle| \begin{matrix} L_{a_i} \in L_R; 1 \leq i \leq 10, \\ L_a, L_b, L_c, L_d \in L_R \end{matrix} \right\}$$

be a set vector space of refined labels over the set $S = Q^+ \cup \{0\}$.

$$W = \left\{ \sum_{i=0}^{9} L_{a_i} x^i, (L_{a_1} \quad L_{a_2} \quad L_{a_3} \quad L_{a_4} \quad L_{a_5} \quad L_{a_6}) \middle| L_{a_i} \in L_R; 0 \leq i \leq 9 \right\}$$

be a set vector space of refined labels over the set $S = Q^+ \cup \{0\}$.

Define $T : V \to W$ by

$$T\left( \begin{pmatrix} L_a & L_b \\ L_c & L_d \end{pmatrix} \right) = (L_a, L_b, 0, 0, L_c, L_d)$$

and

$$T\left( \begin{pmatrix} L_{a_1} & L_{a_2} & L_{a_3} & L_{a_4} & L_{a_5} \\ L_{a_6} & L_{a_7} & L_{a_8} & L_{a_9} & L_{a_0} \end{pmatrix} \right) = \sum_{i=0}^{9} L_{a_i} x^i .$$

T is a set linear transformation of V into W of refined labels over the set $S = Q^+ \cup \{0\}$.

It is pertinent to mention here that all set linear transformation of V into W need not be invertible. Likewise we can define this notion for set linear algebra of refined labels.

Now we proceed onto define the notion of DSm semigroup vector space of refined labels defined over the semigroup S and discuss a few of its properties.

**DEFINITION 3.7:** *Let V be a DSm set vector space over the set S if S is an additive semigroup and if the following conditions hold good.*
  *(a) $sv \in V$ for all $s \in S$ and $v \in V$.*
  *(b) $0.v = 0 \in V$ for all $v \in V$ and $0 \in S$, 0 is zero vector.*
  *(c) $(s_1 + s_2) v = s_1 v + s_2 v$ for all $s_1, s_2 \in S$ and $v \in V$.*
*Then we define V to be a DSm semigroup vector space over the semigroup S of refined labels.*

We will illustrate this situation by some examples.

*Example 3.42:* Let

$$V = \left\{ \sum_{i=0}^{20} L_{a_i} x^i, \begin{bmatrix} L_a & L_b & L_c \\ L_d & L_e & L_f \\ L_g & L_h & L_k \end{bmatrix} \middle| \begin{array}{l} L_{a_i}, L_a, L_b, L_c, L_d, L_e, L_f, \\ L_g, L_h, L_k \in L_R; 0 \leq i \leq 20 \end{array} \right\}$$

be a DSm semigroup vector space of refined labels over the semigroup $S = 7Z^+ \cup \{0\}$.

*Example 3.43:* Let

$$V = \left\{ \begin{pmatrix} L_{a_1} & L_{a_2} \\ L_{a_3} & L_{a_4} \\ \vdots & \vdots \\ L_{a_{13}} & L_{a_{14}} \end{pmatrix}, \begin{pmatrix} L_{a_1} & L_{a_2} & L_{a_3} & L_{a_4} \\ L_{a_5} & L_{a_6} & L_{a_7} & L_{a_8} \end{pmatrix} \middle| L_{a_i} \in L_R; 1 \leq i \leq 14 \right\}$$

be a semigroup vector space of refined labels over the semigroup $S = Q^+ \cup \{0\}$.

*Example 3.44:* Let

$$V = \left\{ \begin{bmatrix} L_{a_1} & L_{a_2} & L_{a_3} \\ L_{a_4} & L_{a_5} & L_{a_6} \\ L_{a_7} & L_{a_8} & L_{a_9} \end{bmatrix}, \begin{bmatrix} L_{a_1} & L_{a_2} & L_{a_3} & L_{a_4} & L_{a_5} \\ L_{a_6} & L_{a_7} & L_{a_8} & L_{a_9} & L_{a_{10}} \\ L_{a_{11}} & L_{a_{12}} & L_{a_{13}} & L_{a_{14}} & L_{a_{15}} \\ L_{a_{16}} & L_{a_{17}} & L_{a_{18}} & L_{a_{19}} & L_{a_{20}} \\ L_{a_{21}} & L_{a_{22}} & L_{a_{23}} & L_{a_{24}} & L_{a_{25}} \end{bmatrix}, \left( L_{a_1} L_{a_2} L_{a_3} \right) \middle| \begin{array}{l} L_{a_i} \in L_R; \\ 1 \leq i \leq 25 \end{array} \right\}$$

be a semigroup vector space of refined labels over the semigroup $S = Z^+ \cup \{0\}$.

*Example 3.45:* Let

$$V = \left\{ \begin{bmatrix} L_{a_1} & L_{a_2} \\ 0 & L_{a_3} \\ L_{a_4} & 0 \\ 0 & L_{a_5} \\ L_{a_6} & 0 \\ 0 & L_{a_7} \\ L_{a_8} & 0 \\ 0 & L_{a_9} \end{bmatrix}, \begin{bmatrix} L_{a_1} & L_{a_2} & L_{a_3} & L_{a_4} \\ L_{a_5} & L_{a_6} & L_{a_7} & L_{a_8} \end{bmatrix} \middle| L_{a_i} \in L_R; 1 \leq i \leq 9 \right\}$$

be a semigroup vector space of refined labels over the semigroup $S = 5Z^+ \cup \{0\}$.

Now we can define substructures of two types.

**DEFINITION 3.8:** *Let V be a semigroup vector space of refined labels over the semigroup S under addition with zero. If $W \subseteq V$ (W is proper subset of V) is a semigroup vector space of refined labels over the semigroup S then we define W to be a semigroup vector subspace of V over the semigroup S. If $W \subseteq V$ is such* that for some subsemigroup T ⊆ S, W is a semigroup vector subspace of refined labels of V over the subsemigroup T of S then we define W to be subsemigroup vector subspace of refined labels over the subsemigroup T of the semigroup S.

We will illustrate this by some simple examples.

*Example 3.46:* Let

$$V = \left\{ \sum_{i=0}^{20} L_{a_i} x^i, \begin{bmatrix} L_{a_1} & L_{a_2} & L_{a_3} & L_{a_4} \\ L_{a_5} & L_{a_6} & L_{a_7} & L_{a_8} \\ L_{a_9} & L_{a_{10}} & L_{a_{11}} & L_{a_{12}} \end{bmatrix}, \begin{bmatrix} L_{a_1} \\ L_{a_2} \\ L_{a_3} \\ L_{a_4} \\ L_{a_5} \end{bmatrix} \middle| L_{a_i} \in L_R; 0 \le i \le 12 \right\}$$

be a DSm semigroup vector space of refined labels over the semigroup $S = Z^+ \cup \{0\}$.

Consider

$$W = \left\{ \sum_{i=0}^{20} L_{a_i} x^i, \begin{bmatrix} 0 \\ L_{a_1} \\ 0 \\ L_{a_2} \\ 0 \end{bmatrix} \middle| L_{a_i} \in L_R; 0 \le i \le 20 \right\} \subseteq V$$

be a DSm semigroup vector subspace of refined labels of V over the semigroup $S = Z^+ \cup \{0\}$.

$$P = \left\{ \begin{bmatrix} L_{a_1} & 0 & L_{a_2} & 0 \\ 0 & L_{a_3} & 0 & L_{a_4} \\ L_{a_5} & 0 & L_{a_6} & 0 \end{bmatrix}, \begin{bmatrix} L_{a_1} \\ L_{a_2} \\ L_{a_3} \\ 0 \\ 0 \end{bmatrix} \middle| L_{a_i} \in L_R; 1 \le i \le 6 \right\}$$

be a DSm semigroup vector subspace of refined labels over the semigroup S of V.

*Example 3.47:* Let

$$V = \left\{ \begin{bmatrix} L_{a_1} & L_{a_2} \\ 0 & L_{a_3} \\ L_{a_4} & 0 \\ 0 & L_{a_5} \\ L_{a_6} & 0 \\ 0 & L_{a_7} \\ L_{a_8} & 0 \end{bmatrix}, \begin{bmatrix} L_{a_1} & 0 & L_{a_2} & 0 & L_{a_3} & 0 & L_{a_4} \\ 0 & L_{a_5} & 0 & L_{a_6} & 0 & L_{a_7} & 0 \\ L_{a_8} & 0 & L_{a_9} & 0 & L_{a_{10}} & 0 & L_{a_{11}} \end{bmatrix} \middle| \begin{array}{l} L_{a_i} \in L_R; \\ 1 \le i \le 11 \end{array} \right\}$$

be a DSm semigroup vector space of refined labels over the semigroup $S = Q^+ \cup \{0\}$.

Consider

$$W = \left\{ \begin{bmatrix} 0 & L_{a_2} \\ 0 & L_{a_3} \\ 0 & 0 \\ 0 & L_{a_1} \\ L_{a_4} & 0 \\ 0 & 0 \\ L_{a_5} & 0 \end{bmatrix}, \begin{bmatrix} 0 & 0 & L_{a_1} & 0 & 0 & 0 & L_{a_2} \\ 0 & L_{a_3} & 0 & 0 & 0 & L_{a_4} & 0 \\ L_{a_5} & 0 & L_{a_6} & 0 & 0 & 0 & L_{a_7} \end{bmatrix} \middle| \begin{array}{l} L_{a_i} \in L_R; \\ 1 \le i \le 7 \end{array} \right\}$$

$\subseteq V$, W is a DSm semigroup vector subspace of refined labels over the semigroup $S = Q^+ \cup \{0\}$.

*Example 3.48:* Let

$$V = \left\{ \begin{bmatrix} L_{a_1} & L_{a_2} \\ L_{a_3} & L_{a_4} \\ L_{a_5} & L_{a_6} \end{bmatrix}, (L_{a_1}\ L_{a_2}\ L_{a_3}\ L_{a_4}\ L_{a_5}\ L_{a_6}), \begin{bmatrix} L_{a_1} & L_{a_2} & L_{a_3} \\ L_{a_4} & L_{a_5} & L_{a_6} \\ L_{a_7} & L_{a_8} & L_{a_9} \end{bmatrix} \middle| \begin{array}{l} L_{a_i} \in L_R; \\ 1 \le i \le 9 \end{array} \right\}$$

be a DSm semigroup vector space of refined labels over the semigroup $S = Z^+ \cup \{0\}$.

Consider

$$W = \left\{ \begin{bmatrix} L_{a_1} & 0 \\ 0 & L_{a_2} \\ 0 & L_{a_3} \end{bmatrix}, (L_{a_1}, L_{a_2}, 0, 0, L_{a_3}, L_{a_4}) \, \middle| \, L_{a_i} \in L_R; 1 \le i \le 4 \right\}$$

$\subseteq V$ and $P = 8Z^+ \cup \{0\} \subseteq S$; P is a subsemigroup of the semigroup S and W is a semigroup vector space of refined labels over the semigroup P. Hence W is a DSm subsemigroup vector subspace of refined labels V over the subsemigroup P of S.

*Example 3.49:* Let

$$V = \left\{ \begin{bmatrix} L_{a_1} & L_{a_2} & L_{a_3} \\ L_{a_4} & L_{a_5} & L_{a_6} \\ L_{a_7} & L_{a_8} & L_{a_9} \\ L_{a_{10}} & L_{a_{11}} & L_{a_{12}} \\ L_{a_{13}} & L_{a_{14}} & L_{a_{15}} \\ L_{a_{16}} & L_{a_{17}} & L_{a_{18}} \\ L_{a_{19}} & L_{a_{20}} & L_{a_{21}} \end{bmatrix}, \begin{bmatrix} L_{a_1} & 0 & L_{a_2} & 0 & L_{a_3} & 0 & L_{a_4} & 0 & L_{a_5} \\ L_{a_6} & L_{a_7} & L_{a_8} & L_{a_9} & L_{a_{10}} & L_{a_{11}} & L_{a_{12}} & L_{a_{13}} & L_{a_{14}} \\ 0 & L_{a_{15}} & 0 & L_{a_{16}} & 0 & L_{a_{17}} & 0 & L_{a_{18}} & 0 \end{bmatrix} \right.$$

$L_{a_i} \in L_R$; $1 \le i \le 21 \}$ be a DSm semigroup vector space of refined labels over the semigroup $S = Q^+ \cup \{0\}$.

Consider

$$W = \left\{ \begin{bmatrix} L_{a_1} & 0 & L_{a_2} \\ 0 & L_{a_3} & 0 \\ L_{a_4} & 0 & L_{a_5} \\ 0 & L_{a_6} & 0 \\ L_{a_7} & 0 & L_{a_8} \\ 0 & L_{a_9} & 0 \\ L_{a_{10}} & 0 & L_{a_{11}} \end{bmatrix}, \begin{bmatrix} 0 & 0 & L_{a_2} & 0 & L_{a_4} & 0 & 0 & 0 & L_{a_6} \\ 0 & L_{a_1} & 0 & 0 & 0 & L_{a_5} & 0 & L_{a_7} & 0 \\ 0 & 0 & 0 & L_{a_3} & 0 & 0 & 0 & 0 & 0 \end{bmatrix} \, \middle| \, \begin{array}{l} L_{a_i} \in L_R; \\ 1 \le i \le 11 \end{array} \right\} \subseteq V$$

is a DSm subsemigroup vector subspace of refined labels of V over the subsemigroup $T = Z^+ \cup \{0\} \subseteq S$.

*Example 3.50:* Let

$$V = \left\{ \begin{bmatrix} L_{a_1} & L_{a_2} \\ 0 & L_{a_3} \\ 0 & 0 \\ L_{a_4} & 0 \\ 0 & L_{a_5} \\ L_{a_6} & 0 \\ 0 & L_{a_7} \\ L_{a_8} & 0 \end{bmatrix}, (L_{a_1}, L_{a_2}, ..., L_{a_9}), \sum_{i=0}^{8} L_{a_i} x^i \,\middle|\, L_{a_i} \in L_R; 0 \leq i \leq 8 \right\}$$

be a DSm semigroup vector space of refined labels over the semigroup $S = Z^+ \cup \{0\}$.

Consider

$$W = \left\{ \sum_{i=0}^{4} L_{a_i} x^i, (L_{a_1}\, 0\, L_{a_2}\, 0\, L_{a_3}\, 0\, L_{a_4}\, 0\, L_{a_5}), \begin{bmatrix} L_{a_1} & 0 \\ 0 & 0 \\ 0 & 0 \\ L_{a_2} & 0 \\ 0 & 0 \\ L_{a_3} & 0 \\ 0 & 0 \\ L_{a_4} & 0 \end{bmatrix} \,\middle|\, \begin{array}{l} L_{a_i} \in L_R; \\ 0 \leq i \leq 5 \end{array} \right\}$$

$\subseteq V$ be a DSm subsemigroup vector subspace of refined labels over the subsemigroup $P = 8Z^+ \cup \{0\} \subseteq S$ of V over S.

Now we just define the notion of DSm semigroup linear algebra of refined labels over the semigroup S.

If V be a DSm semigroup vector space of refined labels over the semigroup S. If V is itself a semigroup under addition with zero then we define V to be a DSm semigroup linear algebra of refined labels over the semigroup S.

We will illustrate this situation by some examples.

*Example 3.51:* Let

$$V = \left\{ \begin{bmatrix} L_{a_1} & L_{a_2} \\ L_{a_3} & L_{a_4} \\ L_{a_5} & L_{a_6} \\ L_{a_7} & L_{a_8} \end{bmatrix} \middle| L_{a_i} \in L_R; 1 \le i \le 8 \right\}$$

be a semigroup linear algebra of refined labels over the semigroup $S = Z^+ \cup \{0\}$.

*Example 3.52:* Let

$$V = \left\{ \begin{bmatrix} L_{a_1} & L_{a_2} \\ L_{a_3} & L_{a_4} \\ L_{a_5} & L_{a_6} \\ L_{a_7} & L_{a_8} \\ L_{a_9} & L_{a_{10}} \\ L_{a_{11}} & L_{a_{12}} \\ L_{a_{13}} & L_{a_{14}} \\ L_{a_{15}} & L_{a_{16}} \end{bmatrix} \middle| L_{a_i} \in L_R; 1 \le i \le 16 \right\}$$

be a semigroup linear algebra of refined labels over the semigroup $S = R^+ \cup \{0\}$.

*Example 3.53:* Let

$$V = \left\{ \sum_{i=0}^{25} L_{a_i} x^i \middle| L_{a_i} \in L_R; 0 \le i \le 25 \right\}$$

be a semigroup linear algebra of refined labels over the semigroup $S = Z^+ \cup \{0\}$.

**THEOREM 3.2:** *Let V be a DSm semigroup linear algebra of refined labels over the semigroup S. V is a DSm semigroup vector space of refined labels over the semigroup S. If V is a DSm semigroup vector space of refined labels over the* *semigroup S then in general V is not DSm semigroup linear algebra of refined labels over the semigroup S.*

The proof is direct and hence is left as an exercise for the reader to prove.

**Example 3.54:** Let

$$V = \left\{ \sum_{i=0}^{8} L_{a_i} x^i, \begin{bmatrix} L_{a_1} & L_{a_2} \\ L_{a_3} & L_{a_4} \\ L_{a_5} & L_{a_6} \end{bmatrix}, (L_{a_1}, L_{a_2}, ..., L_{a_8}) \,\middle|\, L_{a_i} \in L_R; 0 \le i \le 8 \right\}$$

be a DSm semigroup vector space of refined labels over the semigroup $S = Z^+ \cup \{0\}$.
Let

$$W = \left\{ (L_{a_1} \, 0 \, L_{a_2} \, 0 \, L_{a_3} \, 0 \, L_{a_4} \, 0), \sum_{i=0}^{3} L_{a_i} x^i, \begin{bmatrix} L_{a_1} & 0 \\ 0 & L_{a_2} \\ L_{a_3} & 0 \end{bmatrix} \,\middle|\, \begin{matrix} L_{a_i} \in L_R; \\ 0 \le i \le 4 \end{matrix} \right\}$$

$\subseteq V$. W is only a pseudo DSm set vector subspace of refined labels over the set $P = 3Z^+ \cup 8Z^+ \cup 13Z^+ \cup \{0\} \subseteq Z^+ \cup \{0\} = S$ of V.

**Example 3.55:** Let

$$V = \left\{ \begin{bmatrix} L_{a_1} \\ L_{a_2} \\ \vdots \\ L_{a_{15}} \end{bmatrix}, \sum_{i=0}^{20} L_{a_i} x^i, \begin{bmatrix} L_{a_1} & L_{a_2} \\ L_{a_3} & L_{a_4} \end{bmatrix} \,\middle|\, L_{a_i} \in L_R; 0 \le i \le 20 \right\}$$

be a DSm semigroup vector space of refined labels over the semigroup $S = Z^+ \cup \{0\}$.
Consider

$$W = \left\{ \sum_{i=0}^{10} L_{a_i} x^i, \begin{bmatrix} L_{a_1} & L_{a_2} \\ 0 & 0 \end{bmatrix}, \begin{bmatrix} L_{a_1} \\ 0 \\ L_{a_2} \\ 0 \\ L_{a_3} \\ 0 \\ L_{a_4} \\ 0 \\ L_{a_5} \\ 0 \\ L_{a_6} \\ 0 \\ L_{a_7} \\ 0 \\ L_{a_8} \end{bmatrix} \middle| L_{a_i} \in L_R; 0 \leq i \leq 10 \right\}$$

$\subseteq V$, $W$ is a pseudo DSm set vector subspace of refined labels of V over the set $P = \{3Z^+ \cup 2Z^+ \cup 11Z^+ \cup \{0\}\} \subseteq S$.

*Example 3.56:* Let

$$V = \left\{ \begin{bmatrix} L_{a_1} & L_{a_2} \\ L_{a_3} & L_{a_4} \\ L_{a_5} & L_{a_6} \\ L_{a_7} & L_{a_8} \\ L_{a_9} & L_{a_{10}} \end{bmatrix} \middle| L_{a_i} \in L_R; 1 \leq i \leq 10 \right\}$$

be a DSm semigroup linear algebra of refined labels over the semigroup $S = Q^+ \cup \{0\}$.

Let

$$W = \left\{ \begin{bmatrix} L_{a_1} & 0 \\ 0 & L_{a_2} \\ L_{a_3} & 0 \\ 0 & L_{a_4} \\ L_{a_5} & 0 \end{bmatrix} \middle| L_{a_i} \in L_R; 1 \leq i \leq 5 \right\} \subseteq V,$$

be a DSm semigroup linear subalgebra of refined labels over the semigroup $S = Q^+ \cup \{0\}$.

*Example 3.57:* Let

$$V = \left\{ \sum_{i=0}^{28} L_{a_i} x^i \middle| L_{a_i} \in L_R; 0 \leq i \leq 28 \right\}$$

be a DSm semigroup linear algebra over the semigroup $S = Q^+ \cup \{0\}$. Consider

$$W = \left\{ \sum_{i=0}^{12} L_{a_i} x^i \middle| L_{a_i} \in L_R; 0 \leq i \leq 12 \right\} \subseteq V$$

is a DSm semigroup linear subalgebra of refined labels over the semigroup S of V.

*Example 3.58:* Let

$$V = \left\{ \begin{pmatrix} L_{a_1} & L_{a_2} & L_{a_3} \\ L_{a_4} & L_{a_5} & L_{a_6} \\ L_{a_7} & L_{a_8} & L_{a_9} \\ L_{a_{10}} & L_{a_{11}} & L_{a_{12}} \end{pmatrix} \middle| L_{a_i} \in L_R; 1 \leq i \leq 12 \right\}$$

be a DSm semigroup linear algebra of refined labels over the semigroup $S = Z^+ \cup \{0\}$.

$$W = \left\{ \begin{pmatrix} L_{a_1} & 0 & L_{a_2} \\ 0 & L_{a_3} & 0 \\ L_{a_4} & 0 & L_{a_5} \\ 0 & L_{a_6} & 0 \end{pmatrix} \middle| L_{a_i} \in L_R; 1 \leq i \leq 6 \right\} \subseteq V,$$

be a DSm semigroup linear subalgebra of refined labels of V over the semigroup S.

*Example 3.59:* Let V = {All $10 \times 10$ refined labels from $L_R$} be a DSm semigroup linear algebra of refined labels over the semigroup $S = R^+ \cup \{0\}$. Let W = {$10 \times 10$ diagonal matrices with entries from $L_R$} $\subseteq$ V be the DSm semigroup linear subalgebra of refined labels over the semigroup $S = R^+ \cup \{0\}$ of V.

We can see for these semigroup linear algebras of refined labels also we can define the notion of pseudo DSm subset linear subalgebra of refined labels. This task is left as an exercise to the reader. However we give examples of them.

*Example 3.60:* Let

$$V = \left\{ \begin{bmatrix} L_{a_1} & L_{a_2} \\ L_{a_3} & L_{a_4} \\ L_{a_5} & L_{a_6} \\ L_{a_7} & L_{a_8} \\ L_{a_9} & L_{a_{10}} \\ L_{a_{11}} & L_{a_{12}} \end{bmatrix} \middle| L_{a_i} \in L_R ; 1 \leq i \leq 12 \right\}$$

be a DSm semigroup linear algebra of refined labels over the semigroup $S = Z^+ \cup \{0\}$. Consider $T = 3Z^+ \cup 2Z^+ \cup 7Z^+ \cup \{0\} \subseteq S$ a proper subset of S. Take

$$W = \left\{ \begin{bmatrix} L_{a_1} & 0 \\ 0 & L_{a_2} \\ L_{a_3} & 0 \\ 0 & L_{a_4} \\ L_{a_5} & 0 \\ 0 & L_{a_6} \end{bmatrix} \middle| L_{a_i} \in L_R ; 1 \leq i \leq 6 \right\} \subseteq V,$$

W is a pseudo DSm set linear subalgebra of refined labels over the set T of S of V.

*Example 3.61:* Let
$$V = \left\{ \sum_{i=0}^{140} L_{a_i} x^i \,\middle|\, L_{a_i} \in L_R; 0 \le i \le 140 \right\}$$
be a DSm semigroup linear algebra of refined labels over the semigroup $S = Q^+ \cup \{0\}$. Consider
$$W = \left\{ \sum_{i=0}^{20} L_{a_i} x^i \,\middle|\, L_{a_i} \in L_R; 0 \le i \le 20 \right\} \subseteq V,$$
W is a pseudo DSm subset linear subalgebra of V over the subset $T = 17Z^+ \cup 19Z^+ \cup \{0\} \subseteq S$.

*Example 3.62:* Let
$$V = \left\{ \begin{bmatrix} L_{a_1} & L_{a_2} & L_{a_3} \\ L_{a_4} & L_{a_5} & L_{a_6} \\ L_{a_7} & L_{a_8} & L_{a_9} \end{bmatrix} \,\middle|\, L_{a_i} \in L_R; 1 \le i \le 9 \right\}$$
be a DSm semigroup linear algebra of refined labels over the semigroup $S = Z$. Consider $T = 3Z^+ \cup 5Z \subseteq S$ a proper subset of S. Take
$$W = \left\{ \begin{bmatrix} L_{a_1} & 0 & L_{a_2} \\ 0 & L_{a_3} & 0 \\ L_{a_4} & 0 & L_{a_5} \end{bmatrix} \,\middle|\, L_{a_i} \in L_R; 1 \le i \le 5 \right\} \subseteq V,$$
W is a pseudo DSm linear subalgebra of refined labels of V over the subset T of S.

Now having seen examples we can define linear transformation linear operator on these structures in an analogous way with appropriate modifications.

However we recall just the notion of generating subset. Let V be a DSm semigroup of vector space of refined labels over the semigroup S. Let $T = \{v_1, v_2, \ldots, v_n\} \subseteq V$ be a subset of V we say T generates the semigroup vector space of refined labels V over S if every element v ∈ V can be got as v = sv$_i$; v$_i$ ∈ T, s ∈ V. We can as in case of set vector spaces define the notion of direct union or sum and just sum or union or pseudo direct union of subspaces in case of semigroup vector space of linear algebras or semigroup linear algebra of refined labels over the semigroup S.

Let V be a semigroup vector space (linear algebra) of refined labels over the semigroup S. Suppose $W_1, W_2, \ldots, W_n$ be n semigroup vector subspaces (linear subalgebra) of refined labels over the semigroup S, such that V = ∪ W$_i$, and W$_i$ ∩ W$_j$ = φ or {0} if i ≠ j then we say V is the direct union of the semigroup vector subspaces (linear subalgebras) of the semigroup vector space (or linear algebra) of refined labels over the semigroup S.

We will illustrate this situation by some examples.

*Example 3.63:* Let

$$V = \left\{ \sum_{i=0}^{5} L_{a_i} x^i, \begin{bmatrix} L_{a_1} \\ L_{a_2} \\ L_{a_3} \\ L_{a_4} \end{bmatrix}, (L_{a_1}, L_{a_2}, \ldots, L_{a_9}), \begin{pmatrix} L_{a_1} & L_{a_2} \\ L_{a_3} & L_{a_4} \end{pmatrix} \middle| L_{a_i} \in L_R; 0 \le i \le 9 \right\}$$

be a semigroup vector space of refined labels over the semigroup S = Z. Consider

$$W_1 = \left\{ \sum_{i=0}^{5} L_{a_i} x^i \middle| L_{a_i} \in L_R; 0 \le i \le 5 \right\} \subseteq V,$$

$$W_2 = \left\{ \begin{bmatrix} L_{a_1} \\ L_{a_2} \\ L_{a_3} \\ L_{a_4} \end{bmatrix} \middle| L_{a_i} \in L_R; 1 \le i \le 4 \right\} \subseteq V,$$

$$W_3 = \{(L_{a_1}, L_{a_2}, ..., L_{a_9}) \mid L_{a_i} \in L_R; 1 \le i \le 9\}$$

and

$$W_4 = \left\{ \begin{pmatrix} L_{a_1} & L_{a_2} \\ L_{a_3} & L_{a_4} \end{pmatrix} \middle| L_{a_i} \in L_R; 1 \le i \le 4 \right\} \subseteq V$$

be the four semigroup vector subspaces of V of refined labels over the semigroup S = Z. Cleary $V = \bigcup_{i=1}^{4} W_i$ and $W_i \cap W_j = \phi$ if $i \ne j$, $1 \le i, j \le 4$.

Thus V is a direct union of $W_1$, $W_2$, $W_3$ and $W_4$ over S.

*Example 3.64:* Let

$$V = \left\{ \begin{bmatrix} L_{a_1} & L_{a_2} \\ L_{a_3} & L_{a_4} \\ \vdots & \vdots \\ L_{a_{21}} & L_{a_{22}} \end{bmatrix}, \begin{pmatrix} L_{a_1} & L_{a_2} & L_{a_3} & L_{a_4} & L_{a_5} \\ L_{a_6} & L_{a_7} & L_{a_8} & L_{a_9} & L_{a_{10}} \end{pmatrix}, \begin{bmatrix} L_{a_1} & L_{a_2} & L_{a_3} \\ L_{a_4} & L_{a_5} & L_{a_6} \\ L_{a_7} & L_{a_8} & L_{a_9} \end{bmatrix} \right.$$

$L_{a_i} \in L_R; 1 \le i \le 22$} be a DSm semigroup vector space of refined labels over the semigroup $S = 3Z^+ \cup \{0\}$.

Consider

$$W_1 = \left\{ \begin{bmatrix} L_{a_1} & L_{a_2} \\ L_{a_3} & L_{a_4} \\ \vdots & \vdots \\ L_{a_{21}} & L_{a_{22}} \end{bmatrix} \middle| L_{a_i} \in L_R; 1 \le i \le 22 \right\} \subseteq V$$

a DSm semigroup vector subspace of refined labels over the semigroup S.

$$W_2 = \left\{ \begin{pmatrix} L_{a_1} & L_{a_2} & L_{a_3} & L_{a_4} & L_{a_5} \\ L_{a_6} & L_{a_7} & L_{a_8} & L_{a_9} & L_{a_{10}} \end{pmatrix} \middle| \begin{array}{l} L_{a_i} \in L_R; \\ 1 \le i \le 10 \end{array} \right.$$

⊆ V be a DSm semigroup vector subspace of V of refined labels over the semigroup S.

$$W_3 = \left\{ \begin{bmatrix} L_{a_1} & L_{a_2} & L_{a_3} \\ L_{a_4} & L_{a_5} & L_{a_6} \\ L_{a_7} & L_{a_8} & L_{a_9} \end{bmatrix} \middle| L_{a_i} \in L_R; 1 \leq i \leq 9 \right\} \subseteq V$$

be a DSm semigroup vector subspace of V of refined labels over the semigroup S. Clearly $V = \cup W_i = W_1 \cup W_2 \cup W_3$ and $W_i \cap W_j = \phi$, if $i \neq j$; $1 \leq i, j \leq 3$.

Suppose we have a DSm semigroup vector space of refined labels V over the semigroup S. Let $W_1, W_2, \ldots, W_n$ be semigroup vector subspaces of refined labels of V over the semigroup S such that $V = \bigcup_{i=1}^{n} W_i$ with $W_i \cap W_j \neq \phi$ or $\{0\}$ if $i \neq j$; $1 \leq i, j \leq n$, then we define V to be a pseudo direct sum of semigroup vector subspaces of refined labels of V over the semigroup S.

We will illustrate this situation by some examples.

*Example 3.65:* Let

$$V = \left\{ \sum_{i=0}^{20} L_{a_i} x^i, \begin{pmatrix} L_{a_1} & L_{a_2} & L_{a_{23}} \\ L_{a_3} & L_{a_4} & L_{a_{24}} \\ \vdots & \vdots & \vdots \\ L_{a_{21}} & L_{a_{22}} & L_{a_{33}} \end{pmatrix}, \begin{bmatrix} L_{a_1} & L_{a_2} & L_{a_3} \\ L_{a_4} & L_{a_5} & L_{a_6} \end{bmatrix}, \begin{bmatrix} L_{a_1} & L_{a_2} & L_{a_3} & L_{a_4} \\ L_{a_5} & L_{a_6} & L_{a_7} & L_{a_8} \\ L_{a_9} & L_{a_{10}} & L_{a_{11}} & L_{a_{12}} \end{bmatrix} \right.$$

$L_{a_i} \in L_R; 0 \leq i \leq 33\}$ be a semigroup vector space of refined labels over the semigroup $S = Q$. Consider the following semigroup vector subspaces of V over the semigroup S.

$$W_1 = \left\{ \sum_{i=0}^{20} L_{a_i} x^i \,\Big|\, L_{a_i} \in L_R; 0 \leq i \leq 20 \right\} \subseteq V,$$

$$W_2 = \left\{ \begin{pmatrix} L_{a_1} & L_{a_2} & L_{a_{23}} \\ L_{a_3} & L_{a_4} & L_{a_{24}} \\ \vdots & \vdots & \vdots \\ L_{a_{21}} & L_{a_{22}} & L_{a_{33}} \end{pmatrix} \,\Big|\, L_{a_i} \in L_R; 1 \leq i \leq 33 \right\} \subseteq V,$$

$$W_3 = \left\{ \begin{bmatrix} L_{a_1} & L_{a_2} & L_{a_3} \\ L_{a_4} & L_{a_5} & L_{a_6} \end{bmatrix} \,\Big|\, L_{a_i} \in L_R; 1 \leq i \leq 6 \right\} \subseteq V$$

and

$$W_4 = \left\{ \begin{bmatrix} L_{a_1} & L_{a_2} & L_{a_3} & L_{a_4} \\ L_{a_5} & L_{a_6} & L_{a_7} & L_{a_8} \\ L_{a_9} & L_{a_{10}} & L_{a_{11}} & L_{a_{12}} \end{bmatrix} \,\Big|\, L_{a_i} \in L_R; 1 \leq i \leq 12 \right\} \subseteq V$$

be the four semigroup vector subspaces of refined labels of V over the semigroup $S = Q$. $V = W_1 \cup W_2 \cup W_3 \cup W_4$ and $W_i \cap W_j = \phi$ if $i \neq j$ $1 \leq i, j \leq 4$.

*Example 3.66:* Let

$$V = \left\{ \sum_{i=0}^{8} L_{a_i} x^i, \begin{pmatrix} L_{a_1} & L_{a_2} \\ L_{a_3} & L_{a_4} \\ L_{a_5} & L_{a_6} \end{pmatrix}, \begin{bmatrix} L_{a_1} \\ L_{a_2} \\ L_{a_3} \\ \vdots \\ L_{a_{10}} \end{bmatrix}, (L_{a_1}, L_{a_2}, \ldots, L_{a_{12}}) \,\Big|\, \begin{array}{l} L_{a_i} \in L_R; \\ 0 \leq i \leq 12 \end{array} \right\}$$

$\subseteq V$ be a semigroup vector subspace of V of refined labels over the semigroup $S = Q$.

Consider

$$W_1 = \left\{ \sum_{i=0}^{8} L_{a_i} x^i, (L_{a_1}, L_{a_2}, 0, 0, ..., 0, L_{a_{10}}, L_{a_{11}}, L_{a_{12}}) \middle| L_{a_i} \in L_R; 0 \le i \le 12 \right\}$$

be a semigroup vector subspace of refined labels over the semigroup S of V.

$$W_2 = \left\{ \begin{bmatrix} L_{a_1} \\ L_{a_2} \\ L_{a_3} \\ \vdots \\ L_{a_{10}} \end{bmatrix}, (L_{a_1}, L_{a_2}, ..., L_{a_{12}}) \middle| L_{a_i} \in L_R; 1 \le i \le 12 \right\} \subseteq V$$

be a semigroup subspace of V of refined labels over the semigroup S.

$$W_3 = \left\{ \begin{pmatrix} L_{a_1} & L_{a_2} \\ L_{a_3} & L_{a_4} \\ L_{a_5} & L_{a_6} \end{pmatrix}, \sum_{i=0}^{4} L_{a_i} x^i \middle| L_{a_i} \in L_R; 0 \le i \le 5 \right\} \subseteq V$$

be a semigroup vector subspace of refined labels of V over the semigroup S.

Clearly $V = W_1 \cup W_2 \cup W_3$ but $W_i \cap W_j \ne \phi$ or (0) for $i \ne j$ and $1 \le i, j \le 3$.

Hence V is the pseudo direct sum of subspaces of refined labels of V over S. Next we proceed onto define the notion of group linear algebra of refined labels and group vector space of refined labels over a group G.

**DEFINITION 3.9**: *Let V be a set with zero of refined labels which is non empty and G be a group under addition. We say V is a group vector space of refined labels over the group G if the following conditions are true*

(i) *For every $v \in V$ and $g \in G$ vg and gv are in V*
(ii) *$0.v = 0$ for every $v \in V$ and $0 \in G$.*

We will give some examples before we proceed onto define substructures in them.

*Example 3.67:* Let

$$V = \left\{ \sum_{i=0}^{9} L_{a_i} x^i, (L_{a_1}, L_{a_2}, ..., L_{a_{20}}) \middle| L_{a_i} \in L_R; 0 \leq i \leq 20 \right\}$$

be a group vector space of refined labels over the group $G = Q$.

*Example 3.68*: Let

$$V = \left\{ \begin{bmatrix} L_{a_1} \\ L_{a_2} \\ L_{a_3} \\ \vdots \\ L_{a_{25}} \end{bmatrix}, (L_{a_1}, L_{a_2}, ..., L_{a_{21}}), \begin{pmatrix} L_{a_1} & L_{a_2} \\ L_{a_3} & L_{a_4} \\ L_{a_5} & L_{a_6} \end{pmatrix} \middle| L_{a_i} \in L_R; 1 \leq i \leq 25 \right\}$$

be a group vector space of refined labels over the group $G = R$.

*Example 3.69:* Let

$$V = \left\{ \sum_{i=0}^{8} L_{a_i} x^i, \begin{bmatrix} L_{a_1} & L_{a_{11}} & L_{a_{21}} & L_{a_{31}} \\ L_{a_2} & L_{a_{12}} & L_{a_{22}} & L_{a_{32}} \\ \vdots & \vdots & \vdots & \vdots \\ L_{a_{10}} & L_{a_{20}} & L_{a_{30}} & L_{a_{40}} \end{bmatrix}, \begin{bmatrix} L_{a_1} & L_{a_2} & L_{a_3} \\ L_{a_4} & L_{a_5} & L_{a_6} \\ L_{a_7} & L_{a_8} & L_{a_9} \end{bmatrix} \middle| \begin{array}{l} L_{a_i} \in L_R; \\ 0 \leq i \leq 40 \end{array} \right\}$$

be a group vector space of refined labels over the group $G = Z$.

Now we will give some examples of substructures likes subvector spaces, pseudo set subvector spaces pseudo semigroup subvector spaces and subgroup vector subspaces of a group vector space of refined labels over a group G.

However the definition is a matter of routine and hence is left as an exercise to the reader.

*Example 3.70:* Let

$$V = \left\{ \sum_{i=0}^{25} L_{a_i} x^i, \begin{bmatrix} L_{a_1} & L_{a_2} & L_{a_3} & L_{a_4} & L_{a_5} \\ L_{a_6} & L_{a_7} & L_{a_8} & L_{a_9} & L_{a_{10}} \\ L_{a_{11}} & L_{a_{12}} & L_{a_{13}} & L_{a_{14}} & L_{a_{15}} \end{bmatrix} \middle| \begin{array}{l} L_{a_i} \in L_R; \\ 0 \le i \le 25 \end{array} \right\}$$

be a group vector space of refined labels over the group $G = Z$.

$$W = \left\{ \sum_{i=0}^{10} L_{a_i} x^i, \begin{bmatrix} L_{a_1} & 0 & 0 & 0 & L_{a_4} \\ L_{a_2} & 0 & 0 & 0 & L_{a_5} \\ L_{a_3} & 0 & 0 & 0 & L_{a_6} \end{bmatrix} \middle| L_{a_i} \in L_R; 0 \le i \le 10 \right\} \subseteq V$$

a group vector subspace of refined labels of V over the group $G = Z$.

Let

$$P = \left\{ \sum_{i=0}^{5} L_{a_i} x^i, \begin{bmatrix} L_{a_1} & 0 & 0 & L_{a_2} & 0 \\ 0 & 0 & L_{a_4} & 0 & 0 \\ L_{a_3} & 0 & 0 & 0 & L_{a_5} \end{bmatrix} \middle| L_{a_i} \in L_R \right\} \subseteq V$$

is a pseudo set vector subspace of refined labels of V over the set $S = 3Z^+ \cup 5Z^+ \cup \{0\} \subseteq Z = G$.

Consider

$$M = \left\{ \sum_{i=0}^{20} L_{a_i} x^i \middle| L_{a_i} \in L_R; 0 \le i \le 20 \right\} \subseteq V,$$

M is a pseudo semigroup vector subspace of refined labels over the semigroup $S = 3Z^+ \cup \{0\} \subseteq Z = G$.

Let

$$T = \left\{ \sum_{i=0}^{15} L_{a_i} x^i, \begin{bmatrix} L_{a_1} & 0 & L_{a_4} & 0 & L_{a_7} \\ 0 & L_{a_3} & 0 & L_{a_6} & 0 \\ L_{a_2} & 0 & L_{a_5} & 0 & L_{a_8} \end{bmatrix} \middle| L_{a_i} \in L_R; 0 \le i \le 15 \right\}$$

⊆ V be the subgroup vector subspace of V of refined labels over the subgroup H = 2Z ⊆ Z = G.

*Example 3.71:* Let

$$V = \left\{ \begin{bmatrix} L_{a_1} & L_{a_{15}} \\ L_{a_2} & L_{a_{16}} \\ \cdot & \cdot \\ \cdot & \cdot \\ \cdot & \cdot \\ L_{a_{14}} & L_{a_{28}} \end{bmatrix}, \begin{bmatrix} L_{a_1} & L_{a_2} & L_{a_3} \\ L_{a_4} & L_{a_5} & L_{a_6} \\ L_{a_7} & L_{a_8} & L_{a_9} \end{bmatrix}, (L_{a_1}, L_{a_2},...,L_{a_{11}}) \; \middle| \; \begin{matrix} L_{a_i} \in L_R; \\ 1 \le i \le 28 \end{matrix} \right\}$$

be a group linear subalgebras refined labels over the group G = R.

Consider

$$W = \left\{ (L_{a_1}, L_{a_2},...,L_{a_{11}}), \begin{bmatrix} 0 & 0 & 0 \\ L_{a_1} & L_{a_2} & 0 \\ L_{a_3} & 0 & L_{a_4} \end{bmatrix} \; \middle| \; L_{a_i} \in L_R; 1 \le i \le 11 \right\}$$

⊆ V be a group vector subspace of refined labels over the group G = R.

Take

$$P = \left\{ \begin{bmatrix} L_{a_1} & 0 \\ L_{a_2} & 0 \\ \cdot & \cdot \\ \cdot & \cdot \\ \cdot & \cdot \\ L_{a_{14}} & 0 \end{bmatrix}, \begin{bmatrix} L_{a_1} & 0 & 0 \\ L_{a_2} & 0 & L_{a_3} \\ 0 & L_{a_4} & 0 \end{bmatrix}, (L_{a_1}, L_{a_2},...,L_{a_{11}}) \; \middle| \; \begin{matrix} L_{a_i} \in L_R; \\ 1 \le i \le 14 \end{matrix} \right\}$$

⊆ V, P is a pseudo semigroup set vector subspace of refined labels over the semigroup $R^+ \cup \{0\} \subseteq G = R$.

Consider

$$M = \left\{ \begin{bmatrix} 0 & L_{a_1} \\ 0 & L_{a_2} \\ . & . \\ . & . \\ . & . \\ 0 & L_{a_{14}} \end{bmatrix}, (L_{a_1}, 0, L_{a_2}, 0, L_{a_4}, 0, L_{a_5}, 0, L_{a_6}, 0, L_{a_7}) \middle| \begin{array}{l} L_{a_i} \in L_R; \\ 1 \leq i \leq 14 \end{array} \right\}$$

⊆ V, M is a pseudo set vector subspace of V of refined labels over the subgroup $H = Q \subseteq R = G$. However V is not a simple group vector space of refined labels over the group $G = R$.

Now we can define the notion of linear transformation of group vector spaces V and W of refined labels over the group G. If $W = V$ we call the linear transformation as a linear operator.

Now we define group linear algebra of refined labels over the group G and illustrate them with examples.

**DEFINITION 3.10:** *Let V be a group vector space of refined labels over the group G. If V itself is a group under addition and*

$$a(v_1 + v_2) = av_1 + av_2$$

*and*

$$(a_1 + a_2) v = a_1 v + a_2 v$$

*for all v, $v_1$, $v_2 \in V$ and a, $a_1$, $a_2 \in G$ then we define V to be a group linear algebra of refined labels over the group G.*

We will illustrate this situation by some examples.

*Example 3.72:* Let
$$V = \left\{ \sum_{i=0}^{25} L_{a_i} x^i \,\middle|\, L_{a_i} \in L_R; 0 \leq i \leq 25 \right\}$$
be a group linear algebra of refined labels over the group G = Z.

*Example 3.73:* Let
$$V = \left\{ \begin{bmatrix} L_{a_1} & L_{a_2} \\ L_{a_3} & L_{a_4} \\ L_{a_5} & L_{a_6} \\ L_{a_7} & L_{a_8} \\ L_{a_9} & L_{a_{10}} \end{bmatrix} \,\middle|\, L_{a_i} \in L_R; 1 \leq i \leq 10 \right\}$$
be a group linear algebra of refined labels over the group G = R.

*Example 3.74 :* Let M = {All 20 × 20 matrices with entries from $L_R$} be a group linear algebra of refined labels over the group G = Q.

*Example 3.75:* Let
$$W = \left\{ \begin{bmatrix} L_{a_1} & L_{a_{11}} \\ L_{a_2} & L_{a_{12}} \\ . & . \\ . & . \\ L_{a_{10}} & L_{a_{20}} \end{bmatrix} \,\middle|\, L_{a_i} \in L_R; 1 \leq i \leq 20 \right\}$$
be a group linear algebra of refined labels over the group R = G.

Now we have the following interesting observation.

**THEOREM 3.3:** *Let V be a group linear algebra of refined labels over the group G, then V is a group vector space of refined labels over the group G. On the other hand if V is a group vector space of refined labels over the group G then V need not in general be a group vector space of refined labels over the group G.*

The proof is direct and hence left as an exercise for the reader to prove.

We say as in case of semigroup vector spaces when a proper subset P of a group G is a linearly independent set.

Let V be a group vector space of refined labels over the group G. $P \subseteq V$ be a proper subset of V; we say if for any pair of elements $p_1, p_2 \in p$ ($p_1 \neq p_2$) $p_1 = a\, p_2$ or $p_2 = a'\, p_1$ for some a, a' in G then we say P is a linearly dependent set. If for no pair of elements $p_1, p_2$ in P we have $p_1 = a\, p_2$ or $p_2 = a'\, p_1$ then we say P is a linearly independent subset of V.

*Example 3.76:* Let V be a group linear algebra where

$$V = \left\{ \begin{pmatrix} L_{a_1} & L_{a_2} \\ L_{a_3} & L_{a_4} \end{pmatrix} \middle| L_{a_i} \in L_R; 1 \leq i \leq 4 \right\}$$

over the group R = G.
Consider

$$P = \left\{ \begin{pmatrix} L_{a_1} & L_{a_2} \\ 0 & 0 \end{pmatrix}, \begin{pmatrix} L_{b_1} & L_{b_2} \\ 0 & 0 \end{pmatrix}, \begin{pmatrix} L_{a_1} & 0 \\ 0 & 0 \end{pmatrix}, \begin{pmatrix} L_{a_1} & 0 \\ 0 & L_{a_2} \end{pmatrix} \right\} \subseteq V.$$

P is a linearly dependent subset of V over the group G.

*Example 3.77:* Let

$$V = \left\{ \sum_{i=0}^{20} L_{a_i} x^i, \begin{bmatrix} L_{a_1} \\ L_{a_2} \\ L_{a_3} \\ \vdots \\ L_{a_{12}} \end{bmatrix}, (L_{a_1}, L_{a_2}, ..., L_{a_8}) \middle| L_{a_i} \in L_R; 0 \leq i \leq 20 \right\}$$

be a group vector space of refined labels over the group $G = Z$.

Consider

$$P = \left\{ \sum_{i=0}^{10} L_{a_i} x^i (L_{a_1}, L_{a_2}, ..., L_{a_8}) \middle| L_{a_i} \in L_R \right\} \subseteq V.$$

P is a linearly dependent subset of V.
For

$$(L_{a_1}, L_{a_2}, ..., L_{a_8}) = a(L'_{a_1}, L'_{a_2}, ..., L'_{a_8})$$

for $a \in G$. $L_{a_i}, L'_{a_j} \in L_R$.

*Example 3.78:* Let

$$V = \left\{ (L_{a_1}, L_{a_2}, ..., L_{a_9}) \middle| L_{a_i} \in L_R; 1 \leq i \leq 9 \right\}$$

be a group linear algebra of refined labels over the group $G = Q$.
Consider

$$P = \{(L_{a_1}, 0, ..., 0), (0, L_{a_2}, 0, ..., 0),$$
$$(0, 0, L_{a_3}, 0, ..., 0), (0, 0, 0, 0, ..., L_{a_8}, 0),$$
$$(0, 0, ..., 0, L_{a_9})$$

where $L_{a_i} \in L_R; 1 \leq i \leq 9\} \subseteq V$. P is a linearly independent subset of V.

Now having seen the concept of linearly dependent and independent set one can with appropriate modifications build linear transformations and linear operators on group vector spaces of refined labels over the group G.

Also the notions of direct sum and pseudo direct sum in case of group vector spaces of refined labels can be easily obtained without difficulty.

*Example 3.79:* Let

$$V = \left\{ \sum_{i=0}^{9} L_{a_i} x^i, (L_{a_1}, L_{a_2}, ..., L_{a_7}), \begin{bmatrix} L_{a_1} & L_{a_2} & L_{a_3} \\ L_{a_4} & L_{a_5} & L_{a_6} \\ L_{a_7} & L_{a_8} & L_{a_9} \end{bmatrix} \middle| \begin{array}{l} L_{a_i} \in L_R; \\ 0 \leq i \leq 10 \end{array} \right\}$$

be group vector space of refined labels over the group G.

Consider

$$W_1 = \left\{ \sum_{i=0}^{9} L_{a_i} x^i \middle| L_{a_i} \in L_R; 0 \leq i \leq 9 \right\} \subseteq V,$$

$$W_2 = \left\{ (L_{a_1}, L_{a_2}, ..., L_{a_7}) \middle| L_{a_i} \in L_R; 1 \leq i \leq 7 \right\} \subseteq V$$

and

$$W_3 = \left\{ \begin{bmatrix} L_{a_1} & L_{a_2} & L_{a_3} \\ L_{a_4} & L_{a_5} & L_{a_6} \\ L_{a_7} & L_{a_8} & L_{a_9} \end{bmatrix} \middle| L_{a_i} \in L_R; 1 \leq i \leq 9 \right\} \subseteq V$$

be group vector subspaces of V of refined labels over the group G.

Clearly $V = W_1 \cup W_2 \cup W_3$; $W_i \cap W_j = \phi$; $i \neq j$ $1 \leq i, j \leq 3$ so in a direct sum or union.

***Example 3.80:*** Let

$$V = \left\{ \sum_{i=0}^{8} L_{a_i} x^i, \begin{pmatrix} L_{a_1} & L_{a_2} & L_{a_3} \\ L_{a_4} & L_{a_5} & L_{a_6} \\ L_{a_7} & L_{a_8} & L_{a_9} \\ L_{a_{10}} & L_{a_{11}} & L_{a_{12}} \end{pmatrix}, \begin{pmatrix} L_{a_1} & L_{a_3} & L_{a_5} & L_{a_7} \\ L_{a_2} & L_{a_4} & L_{a_6} & L_{a_8} \end{pmatrix} \begin{pmatrix} L_{a_1} & L_{a_2} \\ L_{a_3} & L_{a_4} \end{pmatrix} \middle| \begin{array}{l} L_{a_i} \in L_R; \\ 0 \le i \le 12 \end{array} \right\}$$

be a group vector space of refined labels over the group G = R. Consider

$$W_1 = \left\{ \sum_{i=0}^{8} L_{a_i} x^i, \begin{pmatrix} L_{a_1} & 0 \\ 0 & L_{a_2} \end{pmatrix} \middle| L_{a_i} \in L_R; 0 \le i \le 8 \right\} \subseteq V$$

is a group vector subspace of refined labels of V over R = G.

$$W_2 = \left\{ \begin{pmatrix} L_{a_1} & L_{a_2} \\ L_{a_3} & L_{a_4} \end{pmatrix}, \begin{pmatrix} L_{a_1} & L_{a_3} & L_{a_5} & L_{a_7} \\ L_{a_2} & L_{a_4} & L_{a_6} & L_{a_8} \end{pmatrix} \middle| L_{a_i} \in L_R; 1 \le i \le 8 \right\}$$

$\subseteq V$ is again a group vector subspace of V of refined labels over the group G = R.

$$W_3 = \left\{ \begin{pmatrix} L_{a_1} & L_{a_2} & L_{a_3} \\ L_{a_4} & L_{a_5} & L_{a_6} \\ \vdots & \vdots & \vdots \\ L_{a_{10}} & L_{a_{11}} & L_{a_{12}} \end{pmatrix}, \begin{pmatrix} L_{a_1} & 0 \\ 0 & 0 \end{pmatrix} \middle| L_{a_i} \in L_R; 1 \le i \le 12 \right\} \subseteq V$$

be a group vector subspace of V of refined labels over the group G. Clearly $V = W_1 \cup W_2 \cup W_3$; $W_i \cap W_j \ne \phi$; $i \ne j$, $1 \le i, j \le 3$. Thus V is only a pseudo direct sum.

*Example 3.81:* Let

$$V = \left\{ \sum_{i=0}^{\infty} L_{a_i} x^i \,\middle|\, L_{a_i} \in L_R \right\}$$

be a group linear algebra of refined labels over the group $G = Q$. Let

$$W_1 = \left\{ \sum_{i=0}^{120} L_{a_i} x^i \,\middle|\, L_{a_i} \in L_R \right\} \subseteq V,$$

$$W_2 = \left\{ \sum_{i=121}^{2001} L_{a_i} x^i \,\middle|\, L_{a_i} \in L_R \right\} \subseteq V$$

and

$$W_3 = \left\{ \sum_{i=2002}^{\infty} L_{a_i} x^i \,\middle|\, L_{a_i} \in L_R \right\} \subseteq R$$

be group linear subalgebras of refined labels of V over the group $G = Q$.

Clearly $V = \bigcup_{i=1}^{3} W_i$ with $W_i \cap W_j = \phi$; $i \neq j$ $1 \leq i, j \leq 3$.

*Example 3.82:* Let

$$V = \left\{ \begin{bmatrix} L_{a_1} & L_{a_2} & L_{a_3} \\ L_{a_4} & L_{a_5} & L_{a_6} \\ L_{a_7} & L_{a_8} & L_{a_9} \\ L_{a_{10}} & L_{a_{11}} & L_{a_{12}} \\ L_{a_{13}} & L_{a_{14}} & L_{a_{15}} \end{bmatrix} \,\middle|\, L_{a_i} \in L_R; 1 \leq i \leq 15 \right\}$$

be a group linear algebra of refined labels over the group $G = R$.

Consider

$$W_1 = \left\{ \begin{bmatrix} L_{a_1} & L_{a_2} & L_{a_3} \\ 0 & 0 & 0 \\ 0 & 0 & 0 \\ 0 & 0 & 0 \\ L_{a_4} & L_{a_5} & L_{a_6} \end{bmatrix} \middle| L_{a_i} \in L_R; 1 \leq i \leq 6 \right\} \subseteq V,$$

$$W_2 = \left\{ \begin{bmatrix} 0 & 0 & 0 \\ L_{a_1} & L_{a_2} & L_{a_3} \\ 0 & 0 & 0 \\ 0 & 0 & 0 \\ L_{a_4} & L_{a_5} & L_{a_6} \end{bmatrix} \middle| L_{a_i} \in L_R; 1 \leq i \leq 6 \right\} \subseteq V,$$

$$W_3 = \left\{ \begin{bmatrix} L_{a_1} & L_{a_2} & L_{a_3} \\ L_{a_4} & L_{a_5} & L_{a_6} \\ 0 & 0 & 0 \\ 0 & 0 & 0 \\ 0 & 0 & 0 \end{bmatrix} \middle| L_{a_i} \in L_R; 1 \leq i \leq 6 \right\} \subseteq V$$

and

$$W_4 = \left\{ \begin{bmatrix} 0 & 0 & 0 \\ L_{a_1} & L_{a_2} & L_{a_3} \\ L_{a_4} & L_{a_5} & L_{a_6} \\ L_{a_7} & L_{a_8} & L_{a_9} \\ L_{a_{10}} & L_{a_{11}} & L_{a_{12}} \end{bmatrix} \middle| L_{a_i} \in L_R; 1 \leq i \leq 12 \right\} \subseteq V$$

be group linear subalgebras of refined labels of V over the group G = R. Clearly $W_i \cap W_j \neq \phi$; $i \neq j$ $1 \leq i, j \leq 4$ and $V = W_1 + W_2 + W_3 + W_4$. Thus V is a pseudo direct sum of group linear subalgebras of V over the group G.

**Chapter Four**

# DSM SEMIVECTOR SPACE OF REFINED LABELS

In this chapter we for the first time introduce the notion of semivector space of refined labels, semifield of refined labels and so on.

It is pertinent to state that in the DSm field of refined labels if we take only positive half of reals with 0 say $R^+ \cup \{0\}$ and find

$$L_{R^+ \cup \{0\}} = \left\{ \frac{r}{m+1} \,\middle|\, r \in R^+ \cup \{0\} \right\}$$

then we see $\left(L_{R^+ \cup \{0\}}, +, \times\right)$ to be the DSm semifield of refined labels. Clearly $L_{R^+ \cup \{0\}}$ is isomorphic with $R^+ \cup \{0\}$ and a label is equivalent to a positive real number since for a fixed $m \geq 1$ we have for every $L_a \in L_{R^+ \cup \{0\}}$ there exists a unique $r \in R^+ \cup \{0\}$, $r = \frac{a}{m+1}$ such that $L_a = r$ and reciprocally for every $r$ in $R^+ \cup \{0\}$ there exists a unique $L_a$ in $L_{R^+ \cup \{0\}}$; $L_a = L_{r(m+1)}$ such that $r$ $= L_a$ [34-5]. Further ($L_{R^+\cup\{0\}}$, +, ×, .) where for every '.' is a scalar multiplication $\alpha \in R^+ \cup \{0\}$ and $L_r \in L_{R^+\cup\{0\}}$ we have $\alpha.L_r = L_{\alpha.r} = \alpha.r/m+1 = \alpha.\dfrac{r}{m+1}$ is a semilinear algebra of refined labels over the semifield $R^+ \cup \{0\}$ called DSm semilinear algebra of refined labels.

We will define some more concepts.

Consider $X = \left\{(L_{a_1}, L_{a_2}, ..., L_{a_n}) \mid L_{a_i} \in L_{R^+\cup\{0\}}; 1 \le i \le n\right\}$; X is a DSm semiring of refined labels. X is not a semifield as it has zero divisors. However X is a commutative semiring with unit.

If we take $Y = \left\{(L_{a_1}, L_{a_2}, ..., L_{a_n})^t \mid L_{a_i} \in L_{R^+\cup\{0\}}\right\}$ then Y is not a semiring as product is not defined, but Y is a commutative semigroup under addition.

Consider

$$P = \left\{ \begin{bmatrix} L_{a_1} & L_{a_2} \\ L_{a_3} & L_{a_4} \\ \vdots & \vdots \\ L_{a_{n-1}} & L_{a_n} \end{bmatrix} \middle| L_{a_i} \in L_{R^+\cup\{0\}}; 1 \le i \le n \right\}$$

be a semigroup under addition P is not a DSm semiring as product cannot be defined on P. But if take

$$M = \left\{ \left(L_{a_{ij}}\right)_{n\times n} \middle| L_{a_{ij}} \in L_{R^+\cup\{0\}} \right\};$$

M is a semiring which is not commutative and is not a semifield.

Now let

$$S = \left\{ \sum_{i=0}^{\infty} L_{a_i} x^i \middle| L_{a_i} \in L_{R^+\cup\{0\}} \right\};$$

S is a semiring under usual addition and multiplication. Infact S is a semifield.

Now we will define DSm semivector spaces and DSm semilinear algebras over the semifield. It is pertinent to mention here that $L_{R^+ \cup \{0\}}$ the semifield is isomorphic with the semifield $R^+ \cup \{0\}$.

**DEFINITION 4.1:** *Let V be a semigroup of refined labels with respect to addition with identity zero. S be a semifield. If V is a semivector space over S then we define V to be DSm semivector space of refined labels over S (For more about semivector spaces refer ).*

We illustrate this situation by some examples.

***Example 4.1:*** Let $V = (L_{R^+ \cup \{0\}})$ be the DSm semivector space of refined labels over the semifield $Q^+ \cup \{0\}$.

***Example 4.2:*** Let

$$M = \left\{ \begin{bmatrix} L_{a_1} \\ L_{a_2} \\ L_{a_3} \end{bmatrix} \middle| L_{a_i} \in L_{R^+ \cup \{0\}} \right\}$$

is a DSm semivector space of refined labels over the semifield $S = Z^+ \cup \{0\}$.

***Example 4.3:*** Let

$$K = \left\{ \begin{bmatrix} L_{a_1} & L_{a_2} \\ L_{a_3} & L_{a_4} \\ \vdots & \vdots \\ L_{a_{15}} & L_{a_{16}} \end{bmatrix} \middle| L_{a_i} \in L_{R^+ \cup \{0\}}; 1 \le i \le 16 \right\}$$

is the DSm semivector space of refined labels over the semifield $S = Z^+ \cup \{0\}$.
   It is interesting to note K is not a semivector space over Z.

*Example 4.4:* Let

$$M = \left\{ \begin{bmatrix} L_{a_1} & L_{a_2} & L_{a_3} \\ L_{a_4} & L_{a_5} & L_{a_6} \\ L_{a_7} & L_{a_8} & L_{a_9} \end{bmatrix} \middle| L_{a_i} \in L_{R^+ \cup \{0\}}; 1 \leq i \leq 9 \right\}$$

be a DSm semivector space of refined labels over the semifield $S = Z^+ \cup \{0\}$.

*Example 4.5:* Let

$$P = \left\{ \sum_{i=0}^{25} L_{a_i} x^i \middle| L_{a_i} \in L_{R^+ \cup \{0\}}; 0 \leq i \leq 25 \right\}$$

be a DSm semivector space of refined labels over the semifield $S = Q^+ \cup \{0\}$.

Now we say a DSm semivector space of refined labels V over the semifield $Q^+ \cup \{0\}$ is a semilinear algebra if on V we can define a product and $a(v_1, v_2) = (av_1) v_2$ for every $a \in Q^+ \cup \{0\}$ (any semifield) and $v_1, v_2 \in V$.

We will give some examples of DSm semilinear algebra of refined labels over a semifield.

*Example 4.6:* Let

$$V = \left\{ \sum_{i=0}^{\infty} L_{a_i} x^i \middle| L_{a_i} \in L_{R^+ \cup \{0\}} \right\}$$

be a DSm semilinear algebra of refined labels over the semifield $S = Z^+ \cup \{0\}$.

*Example 4.7:* Let $M = \left\{ (L_{a_1}, L_{a_2}, ..., L_{a_{20}}) \middle| L_{a_i} \in L_{R^+ \cup \{0\}}; 1 \leq i \leq 20 \right\}$ be a semilinear algebra of refined labels over the semifield $Q^+ \cup \{0\}$.

*Example 4.8:* Let

$$T = \left\{ \begin{pmatrix} L_{a_1} & L_{a_2} & L_{a_3} & L_{a_4} \\ L_{a_5} & L_{a_6} & L_{a_7} & L_{a_8} \\ L_{a_9} & L_{a_{10}} & L_{a_{11}} & L_{a_{12}} \\ L_{a_{13}} & L_{a_{14}} & L_{a_{15}} & L_{a_{16}} \end{pmatrix} \middle| L_{a_i} \in L_{R^+ \cup \{0\}}; 1 \le i \le 16 \right\}$$

be a DSm semilinear algebra of refined labels over the semifield $S = R^+ \cup \{0\}$.

It is interesting to see all DSm semivector spaces in general are not DSm semilinear algebras of refined labels. For we see examples 4.2, 4.3 and 4.5 are DSm semivector spaces of refined labels but are not DSm semilinear algebras over their respective semifields.

In view of this we have the following theorem the proof of which is direct and hence is left as an exercise to the reader.

**THEOREM 4.1:** *Let V be a DSm semilinear algebra of refined labels over the semifield S. V is a DSm semivector space of refined labels over the semifield S. Suppose V is a DSm semivector space of refined labels over the semifield S then V need not in general be a DSm semilinear algebra over S.*

Now as in case of DSm linear algebras of refined labels we can define substructures and transformations in DSm semivector spaces and DSm semivector spaces over a semifield S.

*Example 4.9:* Let V be a semivector space of refined labels given by

$$\left\{ \begin{pmatrix} L_{a_1} & L_{a_2} & L_{a_3} & L_{a_4} & L_{a_5} \\ L_{a_6} & L_{a_7} & L_{a_8} & L_{a_9} & L_{a_{10}} \\ L_{a_{11}} & L_{a_{12}} & L_{a_{13}} & L_{a_{14}} & L_{a_{15}} \end{pmatrix} \middle| L_{a_i} \in L_{R^+ \cup \{0\}}; 1 \le i \le 15 \right\}$$

over the semifield $S = R^+ \cup \{0\}$.

Consider

$$W = \left\{ \begin{pmatrix} L_{a_1} & 0 & L_{a_2} & 0 & L_{a_3} \\ 0 & L_{a_4} & 0 & L_{a_5} & 0 \\ L_{a_6} & 0 & L_{a_7} & 0 & L_{a_8} \end{pmatrix} \middle| L_{a_i} \in L_{R^+ \cup \{0\}}; 1 \le i \le 8 \right\}$$

$\subseteq V$ is a semivector subspace of refined labels over $S = R^+ \cup \{0\}$ of V.

*Example 4.10:* Let

$$V = \left\{ \begin{pmatrix} L_{a_1} \\ L_{a_2} \\ \vdots \\ L_{a_{20}} \end{pmatrix} \middle| L_{a_i} \in L_{R^+ \cup \{0\}}; 1 \le i \le 20 \right\}$$

be a DSm semivector space of refined labels over the semifield $S = Q^+ \cup \{0\}$.

$$W = \left\{ \begin{pmatrix} L_{a_1} \\ 0 \\ L_{a_2} \\ \vdots \\ L_{a_{18}} \\ 0 \\ L_{a_{20}} \end{pmatrix} \middle| L_{a_i} \in L_{R^+ \cup \{0\}} \right\} \subseteq V;$$

W is a DSm semivector subspace of refined labels over S.

*Example 4:11:* Let

$$V = \left\{ \begin{pmatrix} L_{a_1} & L_{a_2} \\ L_{a_3} & L_{a_4} \end{pmatrix} \middle| L_{a_i} \in L_{R^+ \cup \{0\}}; 1 \le i \le 4 \right\}$$

be a DSm semilinear algebra of refined labels over the semifield $S = Z^+ \cup \{0\}$.

Consider

$$P = \left\{ \begin{pmatrix} L_{a_1} & 0 \\ 0 & L_{a_2} \end{pmatrix} \middle| L_{a_i} \in L_{R^+ \cup \{0\}}; 1 \leq i \leq 2 \right\} \subseteq V,$$

P is a DSm semilinear subalgebra of refined labels over the semifield $S = Z^+ \cup \{0\}$.

*Example 4.12:* Let

$$M = \left\{ \begin{bmatrix} L_{a_1} & L_{a_2} & L_{a_3} \\ L_{a_4} & L_{a_5} & L_{a_6} \\ L_{a_7} & L_{a_8} & L_{a_9} \end{bmatrix} \middle| L_{a_i} \in L_{R^+ \cup \{0\}}; 1 \leq i \leq 9 \right\}$$

be a DSm semilinear algebra of refined labels over the semifield $S = Q^+ \cup \{0\}$. Consider

$$P = \left\{ \begin{bmatrix} L_{a_1} & L_{a_2} & L_{a_3} \\ 0 & 0 & 0 \\ 0 & 0 & 0 \end{bmatrix} \middle| L_{a_i} \in L_{R^+ \cup \{0\}}; 1 \leq i \leq 3 \right\} \subseteq M;$$

P is a DSm semilinear subalgebra of refined labels over the semifield $S = Q^+ \cup \{0\}$.

*Example 4.13:* Let

$$W = \left\{ \sum_{i=0}^{25} L_{a_i} x^i \middle| L_{a_i} \in L_{R^+ \cup \{0\}} \right\}$$

be a DSm semivector space of refined labels over the semifield $S = Q^+ \cup \{0\}$. Clearly W is only a DSm semivector space of refined labels over $S = Q^+ \cup \{0\}$. Further W is not a DSm semilinear algebra of refined labels over S.

Now

$$T = \left\{ \sum_{i=0}^{10} L_{a_i} x^i \middle| L_{a_i} \in L_R \right\} \subseteq W;$$

T is a DSm semivector subspace of refined labels over $S = Q^+ \cup \{0\}$.

***Example 4.14:*** Let $V = \left\{ \sum_{i=0}^{\infty} L_{a_i} x^i \,\middle|\, L_{a_i} \in L_{R^+ \cup \{0\}} \right\}$
be a DSm semivector space of refined labels over $F = Z^+ \cup \{0\}$. V is also a DSm semilinear algebra of refined labels over $F = Z^+ \cup \{0\}$. Consider

$$T = \left\{ \sum_{i=0}^{\infty} L_{a_i} x^i \,\middle|\, L_{a_i} \in L_{R^+ \cup \{0\}} \right\} \subseteq V;$$

T is a DSm semilinear algebra of refined labels over $F = Z^+ \cup \{0\}$. Now consider V to be a DSm semilinear algebra (vector space) of refined labels over the semifield S. Let $W \subseteq V$ be a DSm semilinear subalgebra (vector subspace) of refined labels over the subsemifield $P \subseteq S$. We call W to be a DSm semilinear subalgebra of refined labels over the subsemifield P of S.

We will illustrate this by some simple examples.

***Example 4.15:*** Let

$$V = \left\{ \begin{bmatrix} L_{a_1} & L_{a_2} & L_{a_3} \\ L_{a_4} & L_{a_5} & L_{a_6} \\ L_{a_7} & L_{a_8} & L_{a_9} \\ L_{a_{10}} & L_{a_{11}} & L_{a_{12}} \\ L_{a_{13}} & L_{a_{14}} & L_{a_{15}} \end{bmatrix} \,\middle|\, L_{a_i} \in L_{R^+ \cup \{0\}}; 1 \le i \le 15 \right\}$$

be a DSm semivector space of refined labels over the semifield $S = Q^+ \cup \{0\}$.

Consider

$$M = \left\{ \begin{bmatrix} L_{a_1} & L_{a_2} & L_{a_3} \\ 0 & 0 & 0 \\ L_{a_4} & L_{a_5} & L_{a_6} \\ 0 & 0 & 0 \\ L_{a_7} & L_{a_8} & L_{a_9} \end{bmatrix} \,\middle|\, L_{a_i} \in L_{R^+ \cup \{0\}}; 1 \le i \le 9 \right\} \subseteq V,$$

M is a DSm semivector subspace of refined labels over the subsemifield $T = Z^+ \cup \{0\} \subseteq S = Q^+ \cup \{0\}$.

*Example 4.16:* Let
$$V = \left\{ \sum_{i=0}^{\infty} L_{a_i} x^i \,\middle|\, L_{a_i} \in L_{R^+ \cup \{0\}} \right\}$$
be a DSm semilinear algebra of refined labels over the semifield $S = R^+ \cup \{0\}$. Consider
$$W = \left\{ \sum_{i=0}^{\infty} L_{a_i} x^{2i} \,\middle|\, L_{a_i} \in L_{R^+ \cup \{0\}} \right\} \subseteq V,$$
W is a subsemifield semilinear subalgebra of refined labels over the subsemifield $T = Z^+ \cup \{0\} \subseteq R^+ \cup \{0\} = S$.

*Example 4.17:* Let
$$V = \left\{ \begin{bmatrix} L_{a_1} & L_{a_2} & L_{a_3} \\ L_{a_4} & L_{a_5} & L_{a_6} \\ L_{a_7} & L_{a_8} & L_{a_9} \end{bmatrix} \,\middle|\, L_{a_i} \in L_{R^+ \cup \{0\}}; 1 \leq i \leq 9 \right\}$$
be a DSm semilinear algebra of refined labels over the semifield $S = Z^+ \cup \{0\}$.

V has no subsemifield linear subalgebras of refined labels but has DSm semilinear subalgebras. For take
$$P = \left\{ \begin{bmatrix} L_{a_1} & L_{a_2} & L_{a_3} \\ 0 & L_{a_4} & L_{a_5} \\ 0 & 0 & L_{a_6} \end{bmatrix} \,\middle|\, L_{a_i} \in L_{R^+ \cup \{0\}} \right\} \subseteq V;$$
P is a DSm semilinear subalgebra of V over the semifield $S = Z^+ \cup \{0\}$.

Now having seen substructures we can write the DSm semivector space as a direct union or direct sum of semivector subspaces as well as pseudo direct sum or union of DSm semivector subspaces. Consider V the semivecctor space of refined labels over the semifield S. Suppose $W_1, W_2, \ldots, W_t$ be t

DSm semivector subspaces of V such that $V = \bigcup_{i=1}^{t} W_i$ with $W_i \cap W_j = \phi$ or (0) if $i \neq j$, $1 \leq i, j \leq t$.

V is a direct sum of semivector subspaces of refined labels over the semifield S. If $V = \bigcup_{i=1}^{t} W_i$ where $W_1, W_2, \ldots, W_t$ are semivector subspaces of refined labels over the semifield S and if $W_i \cap W_j \neq \phi$ if $i \neq j$, $1 \leq i, j \leq t$ then we define V to be a pseudo direct sum or pseudo direct union of semivector subspaces of V over S.

We will first illustrate this situation by some examples.

*Example 4.18:* Let

$$V = \left\{ \begin{bmatrix} L_{a_1} & L_{a_2} \\ L_{a_3} & L_{a_4} \\ L_{a_5} & L_{a_6} \\ L_{a_7} & L_{a_8} \end{bmatrix} \middle| L_{a_i} \in L_{R^+ \cup \{0\}}; 1 \leq i \leq 8 \right\}$$

be a DSm semivector space of refined labels over the field $S = R^+ \cup \{0\}$.

Consider

$$W_1 = \left\{ \begin{bmatrix} L_{a_1} & L_{a_2} \\ 0 & 0 \\ 0 & 0 \\ 0 & 0 \end{bmatrix} \middle| L_{a_i} \in L_{R^+ \cup \{0\}}; 1 \leq i \leq 2 \right\} \subseteq V,$$

$$W_2 = \left\{ \begin{bmatrix} 0 & 0 \\ L_{a_2} & L_{a_1} \\ 0 & 0 \\ L_{a_3} & L_{a_4} \end{bmatrix} \middle| L_{a_i} \in L_{R^+ \cup \{0\}}; 1 \leq i \leq 4 \right\} \subseteq V$$

and

$$W_3 = \left\{ \begin{bmatrix} 0 & 0 \\ 0 & 0 \\ L_{a_1} & L_{a_2} \\ 0 & 0 \end{bmatrix} \middle| L_{a_i} \in L_{R^+ \cup \{0\}}; 1 \leq i \leq 2 \right\} \subseteq V;$$

be a DSm semivector subspaces of V over the semifield $S = R^+ \cup \{0\}$.

Clearly $V = W_1 \cup W_2 \cup W_3$ and

$$W_i \cap W_j = \begin{bmatrix} 0 & 0 \\ 0 & 0 \\ 0 & 0 \\ 0 & 0 \end{bmatrix}$$

if $i \neq j$, $1 \leq i, j \leq 3$.

Take

$$P_1 = \left\{ \begin{bmatrix} L_{a_1} & L_{a_2} \\ 0 & 0 \\ L_{a_3} & 0 \\ 0 & 0 \end{bmatrix} \middle| L_{a_i} \in L_{R^+ \cup \{0\}}; 1 \leq i \leq 3 \right\} \subseteq V,$$

$$P_2 = \left\{ \begin{bmatrix} 0 & 0 \\ L_{a_1} & L_{a_2} \\ L_{a_3} & 0 \\ 0 & 0 \end{bmatrix} \middle| L_{a_i} \in L_{R^+ \cup \{0\}}; 1 \leq i \leq 3 \right\} \subseteq V,$$

$$P_3 = \left\{ \begin{bmatrix} 0 & 0 \\ 0 & L_{a_2} \\ 0 & 0 \\ L_{a_1} & L_{a_3} \end{bmatrix} \middle| L_{a_i} \in L_{R^+ \cup \{0\}}; 1 \leq i \leq 3 \right\} \subseteq V$$

and

$$P_4 = \left\{ \begin{bmatrix} L_{a_1} & 0 \\ 0 & 0 \\ L_{a_2} & L_{a_3} \\ 0 & 0 \end{bmatrix} \middle| L_{a_i} \in L_{R^+ \cup \{0\}}; 1 \leq i \leq 3 \right\} \subseteq V$$

be DSm semivector subspaces of refined labels of V over $S = R^+ \cup \{0\}$.

Clearly $V = \bigcup_{i=1}^{4} W_i$; $W_i \cap W_j \neq \begin{bmatrix} 0 & 0 \\ 0 & 0 \\ 0 & 0 \\ 0 & 0 \end{bmatrix}$ if $i \neq j$, $1 \leq i, j \leq 4$.

Thus V is a pseudo direct sum of DSm semivector subspaces of refined labels over the semifield $S = R^+ \cup \{0\}$.

*Example 4.19:* Let

$$V = \left\{ \begin{bmatrix} L_{a_1} & L_{a_2} \\ L_{a_3} & L_{a_4} \end{bmatrix} \middle| L_{a_i} \in L_{R^+ \cup \{0\}}; 1 \leq i \leq 4 \right\}$$

be a DSm semivector space of refined labels over the semifield $S = Z^+ \cup \{0\}$. Consider

$$W_1 = \left\{ \begin{bmatrix} 0 & L_{a_1} \\ 0 & 0 \end{bmatrix} \middle| L_{a_i} \in L_{R^+ \cup \{0\}} \right\} \subseteq V,$$

$$W_2 = \left\{ \begin{pmatrix} L_{a_1} & 0 \\ 0 & 0 \end{pmatrix} \middle| L_{a_i} \in L_{R^+ \cup \{0\}} \right\} \subseteq V,$$

$$W_3 = \left\{ \begin{pmatrix} 0 & 0 \\ L_{a_1} & 0 \end{pmatrix} \middle| L_{a_i} \in L_{R^+ \cup \{0\}} \right\} \subseteq V$$

and

$$W_4 = \left\{ \begin{pmatrix} 0 & 0 \\ 0 & L_{a_1} \end{pmatrix} \middle| L_{a_i} \in L_{R^+ \cup \{0\}} \right\} \subseteq V$$

are DSm semilinear subalgebras of V of refined labels over the semifield $S = R^+ \cup \{0\}$. Now

$$V = \bigcup_{i=1}^{4} W_i \text{ with } W_i \cap W_j = \begin{bmatrix} 0 & 0 \\ 0 & 0 \end{bmatrix}$$

if $i \neq j$, $1 \leq i, j \leq 4$. Hence V is a direct sum of semivector subspaces.

Consider

$$P_1 = \left\{ \begin{pmatrix} L_{a_1} & L_{a_2} \\ 0 & 0 \end{pmatrix} \middle| L_{a_i} \in L_{R^+ \cup \{0\}}; 1 \leq i \leq 2 \right\} \subseteq V,$$

$$P_2 = \left\{ \begin{pmatrix} 0 & L_{a_1} \\ L_{a_2} & 0 \end{pmatrix} \middle| L_{a_1}, L_{a_2} \in L_{R^+ \cup \{0\}} \right\} \subseteq V,$$

$$P_3 = \left\{ \begin{pmatrix} 0 & 0 \\ L_{a_1} & L_{a_2} \end{pmatrix} \middle| L_{a_1}, L_{a_2} \in L_{R^+ \cup \{0\}} \right\} \subseteq V,$$

$$P_4 = \left\{ \begin{pmatrix} L_{a_1} & L_{a_2} \\ L_{a_3} & 0 \end{pmatrix} \middle| L_{a_1}, L_{a_2} \in L_{R^+ \cup \{0\}}; 1 \leq i \leq 3 \right\} \subseteq V$$

and

$$P_5 = \left\{ \begin{pmatrix} L_{a_1} & 0 \\ 0 & L_{a_2} \end{pmatrix} \middle| L_{a_1}, L_{a_2} \in L_{R^+ \cup \{0\}} \right\} \subseteq V$$

are not DSm semilinear subalgebras but only pseudo DSm semivector subspaces of V and

$$V = \bigcup_{i=1}^{5} P_i \text{ with } P_i \cap P_j \neq \begin{bmatrix} 0 & 0 \\ 0 & 0 \end{bmatrix}$$

if $i \neq j$, $1 \leq i, j \leq 5$. Thus V is a pseudo direct sum of pseudo semivector subspaces of V over the semifield $R^+ \cup \{0\} = S$.

*Example 4.20:* Let

$$V = \left\{ \begin{bmatrix} L_{a_1} & L_{a_2} \\ L_{a_3} & L_{a_4} \\ L_{a_5} & L_{a_6} \\ L_{a_7} & L_{a_8} \\ L_{a_9} & L_{a_{10}} \end{bmatrix} \middle| L_{a_i} \in L_{R^+ \cup \{0\}}; 1 \le i \le 10 \right\}$$

be a DSm semivector space of refined labels over the semifield $S = Z^+ \cup \{0\}$.

Consider

$$W_1 = \left\{ \begin{bmatrix} L_{a_1} & L_{a_2} \\ L_{a_3} & L_{a_4} \\ L_{a_5} & L_{a_6} \\ L_{a_7} & L_{a_8} \\ L_{a_9} & L_{a_{10}} \end{bmatrix} \middle| L_{a_i} \in L_{R^+ \cup \{0\}}; 1 \le i \le 10 \right\}$$

be a DSm semivector space of refined labels over the semifield $S = Z^+ \cup \{0\}$.

Consider

$$W_1 = \left\{ \begin{bmatrix} L_{a_1} & L_{a_2} \\ L_{a_3} & 0 \\ 0 & 0 \\ 0 & 0 \\ 0 & 0 \end{bmatrix} \middle| L_{a_i} \in L_{R^+ \cup \{0\}}; 1 \le i \le 3 \right\} \subseteq V,$$

$$W_2 = \left\{ \begin{bmatrix} 0 & 0 \\ L_{a_1} & L_{a_2} \\ L_{a_3} & 0 \\ 0 & 0 \\ 0 & 0 \end{bmatrix} \middle| L_{a_i} \in L_{R^+ \cup \{0\}}; 1 \le i \le 3 \right\} \subseteq V,$$

$$W_3 = \left\{ \begin{bmatrix} 0 & L_{a_5} \\ 0 & L_{a_6} \\ 0 & 0 \\ L_{a_1} & L_{a_2} \\ L_{a_3} & L_{a_4} \end{bmatrix} \middle| L_{a_i} \in L_{R^+ \cup \{0\}}; 1 \le i \le 6 \right\} \subseteq V,$$

$$W_4 = \left\{ \begin{bmatrix} L_{a_1} & 0 \\ L_{a_2} & 0 \\ 0 & L_{a_3} \\ L_{a_4} & 0 \\ 0 & 0 \end{bmatrix} \middle| L_{a_i} \in L_{R^+ \cup \{0\}}; 1 \le i \le 4 \right\} \subseteq V$$

and

$$W_5 = \left\{ \begin{bmatrix} 0 & L_{a_1} \\ L_{a_2} & L_{a_3} \\ 0 & L_{a_4} \\ 0 & 0 \\ L_{a_5} & 0 \end{bmatrix} \middle| L_{a_i} \in L_{R^+ \cup \{0\}}; 1 \le i \le 5 \right\} \subseteq V$$

are DSm semivector subspaces of refined labels of V over the semifield $S = Z^+ \cup \{0\}$.

*Example 4.21:* Let

$$V = \left\{ \begin{bmatrix} L_{a_1} & L_{a_2} & L_{a_3} \\ L_{a_4} & L_{a_5} & L_{a_6} \\ L_{a_7} & L_{a_8} & L_{a_9} \end{bmatrix} \middle| L_{a_i} \in L_{R^+ \cup \{0\}}; 1 \le i \le 9 \right\}$$

be a DSm semilinear algebra of refined labels over the semifield $S = Q^+ \cup \{0\}$.

Consider

$$W_1 = \left\{ \begin{bmatrix} L_{a_1} & L_{a_2} & L_{a_3} \\ 0 & 0 & 0 \\ 0 & 0 & 0 \end{bmatrix} \middle| L_{a_i} \in L_{R^+ \cup \{0\}}; 1 \le i \le 3 \right\} \subseteq V,$$

$$W_2 = \left\{ \begin{bmatrix} 0 & 0 & 0 \\ L_{a_3} & 0 & 0 \\ 0 & 0 & 0 \end{bmatrix} \middle| L_{a_3} \in L_{R^+ \cup \{0\}} \right\} \subseteq V,$$

$$W_3 = \left\{ \begin{bmatrix} 0 & 0 & 0 \\ 0 & L_{a_4} & L_{a_5} \\ 0 & 0 & 0 \end{bmatrix} \middle| L_{a_4}, L_{a_5} \in L_{R^+ \cup \{0\}} \right\} \subseteq V,$$

$$W_4 = \left\{ \begin{bmatrix} 0 & 0 & 0 \\ 0 & 0 & 0 \\ L_{a_1} & L_{a_2} & L_{a_3} \end{bmatrix} \middle| L_{a_i} \in L_{R^+ \cup \{0\}}; 1 \le i \le 3 \right\} \subseteq V,$$

be a DSm semivector subspaces of V of refined labels over $Q^+ \cup \{0\} = S$. We know

$$V = \bigcup_{i=1}^{4} W_i \; ; W_i \cap W_j = \begin{pmatrix} 0 & 0 & 0 \\ 0 & 0 & 0 \\ 0 & 0 & 0 \end{pmatrix}$$

or (0) if $i \ne j$, $1 \le i, j \le 4$.

V is a direct sum of DSm semivector subspaces of V of refined labels over S.

***Example 4.22:*** Let $V = \left\{ (L_{a_1}, L_{a_2}, ..., L_{a_{10}}) \middle| L_{a_i} \in L_{R^+ \cup \{0\}}; 1 \le i \le 2 \right\}$ be a semivector space over the semifield $S = Z^+ \cup \{0\}$. Let

$$M_1 = \left\{ (L_{a_1}, L_{a_2}, 0...0) \middle| L_{a_1}, L_{a_2} \in L_{R^+ \cup \{0\}} \right\} \subseteq V,$$

$$M_2 = \left\{ (0, 0, L_{a_3}, L_{a_4}, 0, 0, ..., 0) \middle| L_{a_3}, L_{a_4} \in L_{R^+ \cup \{0\}} \right\} \subseteq V,$$

$$M_3 = \left\{(0,0,0,0L_{a_5}, L_{a_6}, 0,0,0,0) \,\middle|\, L_{a_5}, L_{a_6} \in L_{R^+ \cup \{0\}}\right\} \subseteq V,$$

$$M_4 = \left\{(0,0,0,0,0,0,0,0, L_{a_1}, L_{a_2}) \,\middle|\, L_{a_1}, L_{a_2} \in L_{R^+ \cup \{0\}}\right\} \subseteq V,$$

and

$$M_5 = \left\{(0,0,...,0, L_{a_1}, L_{a_2}, 0,0) \,\middle|\, L_{a_1}, L_{a_2} \in L_{R^+ \cup \{0\}}\right\} \subseteq V$$

be DSm semivector subspaces of V of refined labels over the semifield $S = Z^+ \cup \{0\}$.

Clearly

$$V = \bigcup_{i=1}^{5} M_5 \,;\, M_i \cap M_j = (0, 0, 0, \ldots, 0)$$

if $i \neq j$, $1 \leq i$, $j \leq 5$. Thus V is a direct sum of DSm semivector subspaces of refined labels over the semifield $S = Z^+ \cup \{0\}$.

Now we can define the linear transformation and linear operator of DSm semilinear algebras (semivector spaces) over the semifield $S = Z^+ \cup \{0\}$.

*Example 4.23:* Let

$$M = \left\{ \begin{bmatrix} L_{a_1} & L_{a_2} \\ L_{a_3} & L_{a_4} \end{bmatrix} \,\middle|\, L_{a_i} \in L_{R^+ \cup \{0\}}; 1 \leq i \leq 4 \right\}$$

be a DSm semivector space of refined labels over the semifield $S = Z^+ \cup \{0\}$.

Let

$$W = \left\{ (L_{a_1}, L_{a_2}, L_{a_3}, L_{a_4}) \,\middle|\, L_{a_i} \in L_{R^+ \cup \{0\}}; 1 \leq i \leq 4 \right\}$$

be a DSm semilinear algebra of refined labels over the semifield $S = Z^+ \cup \{0\}$. Let $T : M \to W$ be a map such that

$$T\left( \begin{bmatrix} L_{a_1} & L_{a_2} \\ L_{a_3} & L_{a_4} \end{bmatrix} \right) = (L_{a_1} + L_{a_2}, L_{a_2} + L_{a_3}, L_{a_4} + L_{a_3}, L_{a_1} + L_{a_4}).$$

It is easily verified that T is a linear transformation of M to W.

*Example 4.24:* Let

$$V = \left\{ \begin{bmatrix} L_{a_1} & L_{a_2} \\ L_{a_3} & L_{a_4} \\ L_{a_5} & L_{a_6} \\ L_{a_7} & L_{a_8} \end{bmatrix} \middle| L_{a_i} \in L_{R^+ \cup \{0\}}; 1 \le i \le 8 \right\}$$

be a DSm semivector space of refined labels over the semifield $S = Q^+ \cup \{0\}$. Define a map $T : V \to V$ by

$$T\left( \begin{bmatrix} L_{a_1} & L_{a_2} \\ L_{a_3} & L_{a_4} \\ L_{a_5} & L_{a_6} \\ L_{a_7} & L_{a_8} \end{bmatrix} \right) = \begin{bmatrix} L_{a_1} & 0 \\ L_{a_3} & 0 \\ L_{a_5} & 0 \\ L_{a_7} & 0 \end{bmatrix}.$$

T is a linear operator on V.

Suppose

$$W = \left\{ \begin{bmatrix} L_{a_1} & L_{a_2} \\ L_{a_3} & L_{a_4} \\ 0 & 0 \\ 0 & 0 \end{bmatrix} \middle| L_{a_i} \in L_{R^+ \cup \{0\}}; 1 \le i \le 4 \right\} \subseteq V$$

be a DSm semivector subspace of V of refined labels over $S = Q^+ \cup \{0\}$.

Define $T : V \to V$ by

$$T\left( \begin{bmatrix} L_{a_1} & L_{a_2} \\ L_{a_3} & L_{a_4} \\ L_{a_5} & L_{a_6} \\ L_{a_7} & L_{a_8} \end{bmatrix} \right) = \begin{bmatrix} L_{a_1} & L_{a_2} \\ L_{a_3} & L_{a_4} \\ 0 & 0 \\ 0 & 0 \end{bmatrix};$$

T is a projection of V onto W and $T^2 = T$.

Consider P a DSm semivector subspace of V where

$$P = \left\{ \begin{bmatrix} L_{a_1} & L_{a_2} \\ 0 & 0 \\ L_{a_3} & L_{a_4} \\ 0 & 0 \end{bmatrix} \middle| L_{a_i} \in L_{R^+ \cup \{0\}}; 1 \le i \le 4 \right\} \subseteq V$$

is a DSm semivector subspace of V over S of refined labels.

Define $T : V \to V$ by

$$T \left( \begin{bmatrix} L_{a_1} & L_{a_2} \\ L_{a_3} & L_{a_4} \\ L_{a_5} & L_{a_6} \\ L_{a_7} & L_{a_8} \end{bmatrix} \right) = \begin{bmatrix} L_{a_1} & 0 \\ L_{a_2} & 0 \\ 0 & L_{a_3} \\ 0 & L_{a_4} \end{bmatrix} ;$$

Clearly T is not a projection of V onto P.

Clearly kernel of

$$T = \left\{ \begin{bmatrix} L_{a_1} & 0 \\ L_{a_2} & 0 \\ 0 & L_{a_3} \\ 0 & L_{a_4} \end{bmatrix} \middle| L_{a_i} \in L_{R^+ \cup \{0\}}; 1 \le i \le 4 \right\} \subseteq V.$$

It can be easily verified ker T is a semivector subspace of refined labels of V over S.

Several properties can be derived for linear transformation of DSm semivector spaces refined labels over S. Now having defined DSm semivector spaces of refined labels over semifields we see most of the properties related with semivector spaces can be derived in case of DSm semivector spaces can be derived using appropriate modifications.

Now we proceed onto study the set of labels $\tilde{L} \triangleq \{L_0, L_1, \ldots, L_m, L_{m+1}\}$ whose indexes are positive integers between 0 and m+1 that is $0 = L_0 < L_1 < \ldots < L_m < L_{m+1} = 1 \equiv 1$ and $L_i = \dfrac{i}{m+1}$ for $i \in \{0, 1, 2, \ldots, m, m+1\}$.

**THEOREM 4.2**: *Let $0 = L_0 < L_1 < ... < L_m < L_{m+1} = 1$ and $L_i = \dfrac{i}{m+1}$ for $i \in \{0, 1, 2, ..., m, m+1\}$ be the set of labels. $\tilde{L}$ is a chain lattice with $\cap$ and $\cup$ (where $L_i \cap L_j = \min\{L_i, L_j\}$ and $L_i \cup L_j = \max\{L_i, L_j\}$).*

Proof is direct and hence left as an exercise to the reader.

We will now build different algebraic structures using the set of labels $\tilde{L}$.

Let us consider the matrix of finite set of labels or ordinary labels. $H = (L_{p_1}, ..., L_{p_t})$, $0 \le p_i < m+1$ is called the row ordinary labels matrix. $(L_9, L_2, L_0, L_5, L_7)$ where $m + 1 > 9$ is a row ordinary label matrix.

Similarly

$$C = \begin{bmatrix} L_{t_1} \\ L_{t_2} \\ \vdots \\ L_{t_r} \end{bmatrix}$$

where $0 \le t_i \le m+1$ is a column ordinary matrix.

Likewise

$$A = \begin{bmatrix} L_{11} & L_{12} & ... & L_{1n} \\ L_{21} & L_{22} & ... & L_{2n} \\ \vdots & \vdots & & \vdots \\ L_{p_1} & L_{p_2} & ... & L_{p_n} \end{bmatrix}$$

with $L_{ij} \in \tilde{L} = \{L_0, L_1, ..., L_m, L_{m+1}\}$; $0 \le i, j \le m+1$ is a $p \times n$ ordinary label matrix. If $n = p$ we call $A$ to be a square ordinary label matrix.

Similarly we can define $\sum_i L_i x^i$ where $L_i \in \{L_0, L_1, ..., L_m, L_{m+1}\}$ with $x$ a variable or an indeterminate is a polynomial in the variable $x$ with ordinary labels as coefficients. We can now define operations on them and give some algebraic structures on them.

Let $V = \{(L_{a_1},...,L_{a_n}) \mid L_{a_i} \in \{L_0, L_1, L_2,..., L_m, L_{m+1}\}; 1 \le i \le n\}$ be the set of all $1 \times n$ row ordinary label matrix. V is a semigroup under the operation max or min or is a semilattice (and a semigroup). Infact V is not a chain lattice.

Let

$$M = \left\{ \begin{bmatrix} L_{a_1} \\ L_{a_2} \\ \vdots \\ L_{a_r} \end{bmatrix} \middle| L_{a_i} \in \{L_0, L_1, L_2,..., L_m, L_{m+1}\}; 1 \le i \le r \right\}$$

M is also a semilattice / semigroup under $\cap$ (or $\cup$).

Consider

$$P = \left\{ \begin{bmatrix} L_{11} & L_{12} & ... & L_{1n} \\ L_{21} & L_{22} & ... & L_{2n} \\ \vdots & \vdots & & \vdots \\ L_{p_1} & L_{p_2} & ... & L_{p_n} \end{bmatrix} \middle| \begin{array}{l} L_{ij} \in \{L_0, L_1,..., L_{m+1}\}; \\ 0 \le i, j \le m+1 \end{array} \right\}$$

be the collection of $p \times n$ ordinary label matrices. P under $\cup$ and $\cap$ are semilattices or semigroups.

Similarly we can say

$$P = \left\{ \sum_i L_{a_i} x^i \middle| L_{a_i} \in \{0, L_1, L_2,..., L_m, L_{m+1}\} \right\}$$

are semigroups under $\cap$ or $\cup$.

We now show how the operations $\cup$ and $\cap$ are defined on P for one operation alone cannot be defined on P. Consider $p(x) = L_0 x + L_5 + L_2 x^8 + L_3 x^4$ and $q(x) = L_2 + L_1 x^2 + L_0 x^3$ in P.

$p(x) \cap q(x) = (L_5 + L_0 x + L_3 x^4 + L_2 x^8) \cap (L_2 + L_1 x^2 + L_0 x^3)$
$= L_5 \cap L_2 + L_0 \cap L_2 x L_3 \cap L_2 x^4 + L_2 \cap L_2 x^8 + L_5 \cap L_1 x^2 + L_0 \cap L_1 (x \cap x^2) + L_3 \cap L_1 (x^4 \cap x^2) + L_2 \cap L_1(x^8 \cap x^2) + L_5 \cap L_0 x^3 + L_0 \cap L_0 (x \cap x^3) + (L_3 \cap L_0) x^4 \cap x^3) + (L_2 \cap L_0) x^8 \cap x^3$
$= L_2 + L_0 x + L_2 x^4 + L_2 x^8 + L_1 x^2 + L_0 x + L_1 x^2 + L_1 x^2 + L_0 x^3 + L_0 x + L_0 x^3 + L_0 x^3 = L_2 + (L_0 + L_0 + L_0) x + (L_1 + L_2 + L_1) x^2 + (L_0 + L_0 + L_0) x^3 + L_2 x^4 + L_2 x^8$.

The set P cannot be identity $L_0 + L_0 + L_0$ or $L_1 + L_2 + L_1$ so we use '$\cup$' as the operation $L_0 + L_0 = L_0 \cup L_0 = L_0$.

$L_0 + L_0 + L_0 = L_0$ and $L_1 + L_2 + L_1 = L_1 \cup L_2 + L_1 = L_2 \cup L_1 = L_2$. Thus $p(x) q(x) = L_2 + L_0 x + L_2 x^2 + L_0 x^3 + L_2 x^4 + L_2 x^8$.

Now this element is recognalizable by P. Thus P is a semiring or a semifield under $\cup$ and $\cap$.

We will illustrate all these by some examples.

***Example 4.25:*** Let $L = \{L_0, L_1, \ldots, L_m, L_{m+1}\}$ be a semigroup under '$\cup$' L is a DSm semivector space of ordinary labels over the semifield $\{0, 1\}$.

***Example 4.26:*** Let $L = \{L_0, L_1, \ldots, L_m, L_{m+1}\}$ be a semigroup under '$\cap$'. L is a DSm semivector space of ordinary labels over the semifield $S = \{0, 1\}$ (S is lattice).

***Example 4.27:*** Let
$$W = \left\{ \left(L_{a_1}, L_{a_2}, L_{a_3}, L_{a_4}\right) \middle| L_{a_i} \in \{L_0, L_1, L_2, \ldots, L_m, L_{m+1}\}; 1 \leq i \leq 4 \right\}$$
be a DSm semivector space over $S = \{0, 1\}$ where $(W, \cup)$ is a semigroup (That if
$$x = \left(L_{a_1}, L_{a_2}, L_{a_3}, L_{a_4}\right) \text{ and } y = \left(L_{b_1}, L_{b_2}, L_{b_3}, L_{b_4}\right)$$
in W then
$$x \cup y = \left(L_{a_1} \cup L_{b_1}, L_{a_2} \cup L_{b_2}, L_{a_3} \cup L_{b_3}, L_{a_4} \cup L_{b_4}\right)$$
is in W.

***Example 4.28:*** Let
$$M = \left\{ \begin{bmatrix} L_{b_1} \\ L_{b_2} \\ L_{b_3} \\ L_{b_4} \end{bmatrix} \middle| L_{b_i} \in \{0 = L_0, L_1, L_2, \ldots, L_m, L_{m+1}\}; 1 \leq i \leq 4 \right\}$$

be a semigroup under $\cap$. M is a DSm semivector space of ordinary labels over $S = \{0, 1\}$ the semifield for

$$1 \cap \begin{bmatrix} L_{a_1} \\ L_{a_2} \\ L_{a_3} \\ L_{a_4} \end{bmatrix} = \begin{bmatrix} L_{a_1} \\ L_{a_2} \\ L_{a_3} \\ L_{a_4} \end{bmatrix}$$

and

$$0 \cap \begin{bmatrix} L_{a_1} \\ L_{a_2} \\ L_{a_3} \\ L_{a_4} \end{bmatrix} = \begin{bmatrix} 0 \\ 0 \\ 0 \\ 0 \end{bmatrix}.$$

*Example 4.29:* Let

$$L = \left\{ \begin{bmatrix} L_{a_1} & L_{a_2} \\ L_{a_3} & L_{a_4} \end{bmatrix} \middle| L_{a_i} \in \{0 = L_0, L_1, L_2, ..., L_m, L_{m+1} = 1\} \right\}$$

be a semigroup under '$\cup$'. L is a DSm semivector space of refined labels over the semifield $S = \{0, 1\}$.

*Example 4.30:* Let

$$M = \left\{ \begin{bmatrix} L_{a_1} & L_{a_2} & L_{a_3} \\ L_{a_4} & L_{a_5} & L_{a_6} \\ L_{a_7} & L_{a_8} & L_{a_9} \\ L_{a_{10}} & L_{a_{11}} & L_{a_{12}} \\ L_{a_{13}} & L_{a_{14}} & L_{a_{15}} \\ L_{a_{16}} & L_{a_{17}} & L_{a_{18}} \end{bmatrix} \middle| L_{a_i} \in \{0 = L_0, L_1, ..., L_m, L_{m+1}\}; 1 \leq i \leq 18 \right\}$$

is a semigroup under '$\cup$'; so M is a DSm semivector space of ordinary labels over the semifield $S = \{0, 1\}$.

***Example 4.31:*** Let
$$V = \left\{ \sum_i L_{a_i} x^i \, \bigg| \, L_{a_i} \in \{0, L_1, L_2, ..., L_m, L_{m+1}\} \right\}$$
be a semigroup under a special operation denoted by +. V is a semigroup semivector space of ordinary labels over $S = \{0, 1\}$.

For $p(x) = L_9 + L_1 x + L_0 x^2 + L_7 x^3$ and $q(x) = L_0 + L_2 x + L_3 x^3 + L_4 x^5$ in V we define special operation addition on V as follows:
$$p(x) + q(x) = (L_9 + L_0) + (L_1 + L_2) x + L_0 x^2 (L_3 + L_7)x^3 + L_4 x^5$$
where
$$L_9 + L_0 = L_9$$
$$L_2 + L_1 = L_2 \text{ and } L_3 + L_7 = L_7$$
that is $L_i + L_j = \max(i, j)$.

Thus $p(x) + q(x) = L_9 + L_2 x + L_0 x^2 + L_7 x^3 + L_4 x^5$.

Now this gives us a new class of DSm semivector spaces of ordinary labels where the DSm semivector spaces contain only a finite number of elements.

The notion of basis, linearly independent labels (which may or may not exists) linear transformation of DSm semivector space of ordinary labels, linear operators can be derived as in case of usual semivector spaces with some appropriate modifications. This task is left as an exercise to the reader. We will give one or two examples in this direction.

***Example 4.32:*** Let
$$V = \left\{ \begin{bmatrix} L_{a_1} & L_{a_2} \\ L_{a_3} & L_{a_4} \\ L_{a_5} & L_{a_6} \end{bmatrix} \, \bigg| \, L_{a_i} \in \{0 = L_0, L_1, L_2, ..., L_m, L_{m+1}\}; 1 \le i \le 6 \right\}$$
and
$$W = \left\{ (L_{a_1}, L_{a_2}, ..., L_{a_6}) \, \big| \, L_{a_i} \in \{L_0, L_1, L_2, ..., L_m, L_{m+1}\}; 1 \le i \le 6 \right\}$$
be DSm semivector spaces over the semifield $S = \{0, 1\}$ of ordinary labels.

Define T:V → W by

$$T\left(\begin{bmatrix} L_{a_1} & L_{a_2} \\ L_{a_3} & L_{a_4} \\ L_{a_5} & L_{a_6} \end{bmatrix}\right) = \left(L_{a_1}, L_{a_2}, ..., L_{a_6}\right)$$

T is a linear transformation of V into W.

*Example 4.33:* Let

$$V = \left\{ \sum_{i=0}^{8} L_{a_i} x^i \,\middle|\, L_{a_i} \in \{L_0, L_1, L_2, ..., L_m, L_{m+1}\}; 0 \le i \le 8 \right\}$$

be a DSm semivector space of ordinary labels over the semifield S = {0, 1}.

Let

$$W = \left\{ \begin{bmatrix} L_{a_1} & L_{a_2} & L_{a_3} \\ L_{a_4} & L_{a_5} & L_{a_6} \\ L_{a_7} & L_{a_8} & L_{a_9} \end{bmatrix} \,\middle|\, L_{a_i} \in \{L_0, L_1, L_2, ..., L_m, L_{m+1}\}; 1 \le i \le 9 \right\}$$

be a DSm semivector space of ordinary labels over the semifield S = {0, 1}.

Define T : V → W by

$$T\left(\sum_{i=0}^{8} L_{a_i} x^i\right) = \begin{bmatrix} L_{a_1} & L_{a_2} & L_{a_3} \\ L_{a_4} & L_{a_5} & L_{a_6} \\ L_{a_7} & L_{a_8} & L_{a_9} \end{bmatrix};$$

T is a linear transformation from V into W.

*Example 4.34:* Let

$$V = \left\{ \begin{bmatrix} L_{a_1} & L_{a_2} \\ L_{a_3} & L_{a_4} \\ \vdots & \vdots \\ L_{a_{15}} & L_{a_{16}} \end{bmatrix} \,\middle|\, L_{a_i} \in \{L_0, L_1, L_2, ..., L_m, L_{m+1}\}; 1 \le i \le 16 \right\}$$

be a DSm semivector space of ordinary labels over the semifield $S = \{0, 1\}$.

Define $T : V \to V$ by

$$T\left(\begin{bmatrix} L_{a_1} & L_{a_2} \\ L_{a_3} & L_{a_4} \\ \vdots & \vdots \\ L_{a_{15}} & L_{a_{16}} \end{bmatrix}\right) = \begin{bmatrix} L_{a_1} & 0 \\ 0 & L_{a_2} \\ L_{a_3} & 0 \\ 0 & L_{a_4} \\ L_{a_5} & 0 \\ 0 & L_{a_6} \\ L_{a_7} & 0 \\ 0 & L_{a_8} \end{bmatrix}$$

T is a linear operator on V.

*Example 4.35:* Let

$$V = \left\{ \begin{bmatrix} L_{a_1} & L_{a_2} & L_{a_3} \\ L_{a_4} & L_{a_5} & L_{a_6} \\ L_{a_7} & L_{a_8} & L_{a_9} \\ L_{a_{10}} & L_{a_{11}} & L_{a_{12}} \end{bmatrix} \middle| L_{a_i} \in \{L_0, L_1, L_2, ..., L_m, L_{m+1}\}; 1 \le i \le 12 \right\}$$

be a DSm semivector space of ordinary labels over the semifield $S = \{0, 1\}$.

Define $T : V \to V$ by

$$T\left(\begin{bmatrix} L_{a_1} & L_{a_2} & L_{a_3} \\ L_{a_4} & L_{a_5} & L_{a_6} \\ L_{a_7} & L_{a_8} & L_{a_9} \\ L_{a_{10}} & L_{a_{11}} & L_{a_{12}} \end{bmatrix}\right) = \begin{bmatrix} L_{a_1} & 0 & L_{a_2} \\ 0 & L_{a_3} & 0 \\ L_{a_4} & 0 & L_{a_5} \\ 0 & L_{a_6} & 0 \end{bmatrix}.$$

T is a linear operator on V. Now having seen examples linear operator and transformation of DSm semivector space of ordinary labels are now proceed onto define basis for them.

***Example 4.36:*** Let

$$V = \left\{ \begin{bmatrix} L_{a_1} \\ L_{a_2} \end{bmatrix} \,\middle|\, L_{a_i} \in \{L_0, L_1, L_2, ..., L_m, L_{m+1}\}; i = 1, 2 \right\}$$

be a DSm semivector space of ordinary labels over the semifield $S = \{0, 1\}$.

The basis of V is

$$B = \left\{ \begin{bmatrix} L_{a_i} \\ 0 \end{bmatrix}, \begin{bmatrix} 0 \\ L_{a_j} \end{bmatrix} \,\middle|\, L_{a_i}, L_{a_j} \in \{0, L_1, L_2, ..., L_m, L_{m+1}\} \right\}.$$

Thus B is of order $2(m + 2)$. The dimension of V is $2(m + 2)$ under the operation '$\cup$'. For

$$\begin{bmatrix} L_{a_i} \\ 0 \end{bmatrix} \cup \begin{bmatrix} 0 \\ L_{a_j} \end{bmatrix} = \begin{bmatrix} L_{a_i} \\ L_{a_j} \end{bmatrix}$$

and so on.

Interested reader can study basis and linearly independent elements of a DSm semivector space over the semifield $S = \{0, 1\}$ of ordinary labels.

Now we define special type of DSm semivector spaces of ordinary labels over the semifields.

Let V be the set of ordinary labels which is a semigroup under $\cup$ or $\cap$. Consider $\{L_0, L_1, ..., L_m, L_{m+1}\}$ is a semifield under $\cup$ and $\cap$ as S is a lattice. V is a special DSm semivector space of ordinary labels over the semifield $S = \{L_0, L_1, ..., L_m, L_{m+1}\}$.

We will give some examples.

***Example 4.37:*** Let

$$V = \left\{ \begin{bmatrix} L_{a_1} & L_{a_2} \\ L_{a_3} & L_{a_4} \\ L_{a_5} & L_{a_6} \end{bmatrix} \,\middle|\, L_{a_i} \in \{L_0, L_1, L_2, ..., L_m, L_{m+1}\}; 1 \leq i \leq 6 \right\}$$

be a semigroup under '$\cup$'. V is a special DSm semivector space of ordinary labels over the semifield $S = \{L_0, L_1, ..., L_m, L_{m+1}\}$.

*Example 4.38:* Let

$$M = \left\{ \begin{bmatrix} L_{a_1} \\ L_{a_2} \\ \vdots \\ L_{a_{20}} \end{bmatrix} \middle| L_{a_i} \in \{L_0, L_1, L_2, ..., L_m, L_{m+1}\}; 1 \leq i \leq 20 \right\}$$

be a special DSm semivector space of ordinary labels over the semifield $S = \{L_0, L_1, ..., L_m, L_{m+1}\}$.

*Example 4.39:* Let

$$T = \left\{ (L_{a_1}, L_{a_2}, ..., L_{a_{12}}) \middle| L_{a_i} \in \{L_0, L_1, L_2, ..., L_m, L_{m+1}\}; 1 \leq i \leq 12 \right\}$$

be a DSm special semivector space of ordinary labels over the semifield $L = \{L_0, L_1, ..., L_m, L_{m+1}\}$.

*Example 4.40:* Let

$$T = \left\{ \begin{bmatrix} L_{a_1} & L_{a_2} & L_{a_3} \\ L_{a_4} & L_{a_5} & L_{a_6} \\ L_{a_7} & L_{a_8} & L_{a_9} \end{bmatrix} \middle| L_{a_i} \in \{L_0, L_1, L_2, ..., L_m, L_{m+1}\}; 1 \leq i \leq 9 \right\}$$

be a DSm special semivector space of ordinary labels over the semifield $L = \{L_0, L_1, ..., L_m, L_{m+1}\}$.

We will now proceed on define substructures in them.

*Example 4.41:* Let

$$V = \left\{ \begin{bmatrix} L_{a_1} & L_{a_2} & L_{a_3} \\ L_{a_4} & L_{a_5} & L_{a_6} \\ L_{a_7} & L_{a_8} & L_{a_9} \\ L_{a_{10}} & L_{a_{11}} & L_{a_{12}} \end{bmatrix} \middle| L_{a_i} \in \{L_0, L_1, L_2, ..., L_m, L_{m+1}\}; 1 \leq i \leq 12 \right\}$$

be a DSm semivector space of ordinary labels over the semifield $S = \{0, 1\}$.

Consider

$$W = \left\{ \begin{bmatrix} L_{a_1} & L_{a_2} & L_{a_3} \\ 0 & 0 & 0 \\ L_{a_4} & L_{a_5} & L_{a_6} \\ 0 & 0 & 0 \end{bmatrix} \middle| L_{a_i} \in \{L_0, L_1, L_2, ..., L_m, L_{m+1}\}; 1 \le i \le 6 \right\} \subseteq V$$

is a DSm semivector subspace of refined labels over the semifield $S = \{0, 1\}$.

Consider

$$P = \left\{ \begin{bmatrix} L_{a_1} & 0 & L_{a_6} \\ L_{a_2} & L_{a_4} & 0 \\ L_{a_3} & 0 & L_{a_7} \\ 0 & L_{a_5} & 0 \end{bmatrix} \middle| L_{a_i} \in \{L_0, L_1, L_2, ..., L_m, L_{m+1}\}; 1 \le i \le 7 \right\} \subseteq V$$

is a DSm semivector subspace of refined labels over the semifield $S = \{0, 1\}$.

***Example 4.42:*** Let

$$T = \left\{ \sum_{i=0}^{12} L_{a_i} x^i \middle| L_{a_i} \in \{L_0, L_1, L_2, ..., L_m, L_{m+1}\}; 0 \le i \le 12 \right\}$$

be a DSm semivector space of ordinary labels over the semifield $S = \{0, 1\}$.

Take

$$P = \left\{ \sum_{i=0}^{5} L_{a_i} x^i \middle| L_{a_i} \in \{L_0, L_1, L_2, ..., L_m, L_{m+1}\}; 0 \le i \le 5 \right\} \subseteq T;$$

P is a DSm semivector subspace of ordinary labels over the semifield $S = \{0, 1\}$.

***Example 4.43:*** Let

$$W = \left\{ \begin{bmatrix} L_{a_1} & L_{a_2} \\ L_{a_3} & L_{a_4} \\ L_{a_5} & L_{a_6} \end{bmatrix} \middle| L_{a_i} \in \{L_0, L_1, L_2, ..., L_m, L_{m+1}\}; 1 \le i \le 6 \right\}$$

be a special DSm semivector space of ordinary labels over the semifield $S = \{L_0, L_1, ..., L_m, L_{m+1}\}$.

$$M = \left\{ \begin{bmatrix} L_{a_1} & 0 \\ 0 & L_{a_2} \\ L_{a_3} & 0 \end{bmatrix} \middle| L_{a_i} \in \{L_0, L_1, L_2, ..., L_m, L_{m+1}\}; 1 \le i \le 3 \right\} \subseteq W;$$

M is a special DSm semivector subspace of ordinary labels over the semifield S.

*Example 4.44:* Let

$$V = \left\{ \begin{bmatrix} L_{a_1} \\ L_{a_2} \\ \vdots \\ L_{a_{12}} \end{bmatrix} \middle| L_{a_i} \in \{L_0, L_1, L_2, ..., L_m, L_{m+1}\}; 1 \le i \le 12 \right\}$$

be a DSm special semivector space of ordinary labels over the semifield $S = \{L_0, L_1, ..., L_m, L_{m+1}\}$.

$$M = \left\{ \begin{bmatrix} L_{a_1} \\ 0 \\ L_{a_2} \\ 0 \\ \vdots \\ L_{a_{10}} \\ 0 \\ L_{a_{12}} \end{bmatrix} \middle| L_{a_i} \in \{L_0, L_1, L_2, ..., L_m, L_{m+1}\}; 1 \le i \le 12 \right\} \subseteq V,$$

M is a DSm special semivector subspace of ordinary labels over the semifield S of V.

We can define the notion of direct sum and pseudo direct sum as in case of other vector spaces.

We will illustrate this situation by some examples.

*Example 4.45:* Let

$$V = \left\{ \begin{bmatrix} L_{a_1} & L_{a_2} \\ L_{a_3} & L_{a_4} \\ \vdots & \vdots \\ L_{a_{15}} & L_{a_{16}} \end{bmatrix} \middle| L_{a_i} \in \{L_0, L_1, L_2, \ldots, L_m, L_{m+1}\}; 1 \le i \le 16 \right\}$$

be a DSm semivector space of ordinary labels over the semifield $S = \{0, 1\}$.

Let

$$W_1 = \left\{ \begin{bmatrix} L_{a_1} & L_{a_2} \\ 0 & 0 \\ L_{a_3} & L_{a_4} \\ 0 & 0 \\ \vdots & \vdots \\ 0 & 0 \end{bmatrix} \middle| L_{a_i} \in \{L_0, L_1, L_2, \ldots, L_m, L_{m+1}\}; 1 \le i \le 4 \right\} \subseteq V,$$

$$W_2 = \left\{ \begin{bmatrix} 0 & 0 \\ L_{a_1} & L_{a_2} \\ 0 & 0 \\ L_{a_3} & L_{a_4} \\ \vdots & \vdots \\ 0 & 0 \end{bmatrix} \middle| L_{a_i} \in \{L_0, L_1, L_2, \ldots, L_m, L_{m+1}\}; 1 \le i \le 4 \right\} \subseteq V,$$

$$W_3 = \left\{ \begin{bmatrix} 0 & 0 \\ 0 & 0 \\ 0 & 0 \\ 0 & 0 \\ L_{a_1} & L_{a_2} \\ L_{a_3} & L_{a_4} \\ 0 & 0 \\ 0 & 0 \end{bmatrix} \middle| L_{a_i} \in \{L_0, L_1, L_2, \ldots, L_m, L_{m+1}\}; 1 \le i \le 4 \right\} \subseteq V$$

and

$$W_4 = \left\{ \begin{bmatrix} 0 & 0 \\ 0 & 0 \\ 0 & 0 \\ 0 & 0 \\ 0 & 0 \\ 0 & 0 \\ L_{a_1} & L_{a_2} \\ L_{a_3} & L_{a_4} \end{bmatrix} \middle| L_{a_i} \in \{L_0, L_1, L_2, ..., L_m, L_{m+1}\}; 1 \le i \le 4 \right\} \subseteq V,$$

be four DSm semivector space of ordinary labels over the semifield S = {0, 1} of V.

Clearly $V = \bigcup_{i=1}^{4} W_i$ where

$$W_i \cap W_j = \begin{bmatrix} 0 & 0 \\ 0 & 0 \\ 0 & 0 \\ 0 & 0 \\ \vdots & \vdots \\ 0 & 0 \end{bmatrix} \text{ if } i \ne j, 1 \le i, j \le 4.$$

Thus V is a direct sum of DSm semivector subspace of ordinary labels of V over the semifield S = {0, 1}.

*Example 4.46:* Let

$$V = \left\{ \begin{bmatrix} L_{a_1} & L_{a_3} & L_{a_5} & L_{a_7} & L_{a_9} & L_{a_{11}} & L_{a_{13}} & L_{a_{15}} & L_{a_{17}} \\ L_{a_2} & L_{a_4} & L_{a_6} & L_{a_8} & L_{a_{10}} & L_{a_{12}} & L_{a_{14}} & L_{a_{16}} & L_{a_{18}} \end{bmatrix} \right.$$

$L_{a_i} \in \{L_0, L_1, L_2, ..., L_m, L_{m+1}\}; 1 \le i \le 18\}$ be a DSm semivector space of ordinary labels of V over the semifield S = {0, 1}.

$$W_1 = \left\{ \begin{bmatrix} L_{a_1} & 0 & L_{a_3} & 0 & 0 & ... & 0 \\ 0 & L_{a_2} & 0 & L_{a_4} & 0 & ... & 0 \end{bmatrix} \middle| L_{a_i} \in \{L_0, L_1, ..., L_{m+1}\} \right\} \subseteq V;$$

$$W_2 = \left\{ \begin{bmatrix} 0 & L_{a_2} & 0 & L_{a_4} & L_{a_5} & 0 & 0 & \dots & 0 \\ L_{a_1} & 0 & L_{a_3} & 0 & L_{a_6} & 0 & 0 & \dots & 0 \end{bmatrix} \right.$$

$$\left. L_{a_i} \in \{L_0, L_1, L_2, \dots, L_m, L_{m+1}\}; 1 \le i \le 6 \right\} \subseteq V;$$

$$W_3 = \left\{ \begin{bmatrix} 0 & 0 & 0 & 0 & 0 & L_{a_1} & L_{a_2} & 0 & 0 \\ 0 & 0 & 0 & 0 & 0 & L_{a_3} & L_{a_4} & 0 & 0 \end{bmatrix} \middle| \begin{array}{c} L_{a_i} \in \{L_0, L_1, \dots, L_{m+1}\}; \\ 1 \le i \le 4 \end{array} \right\} \subseteq V$$

and

$$W_4 = \left\{ \begin{bmatrix} 0 & 0 & 0 & 0 & 0 & 0 & L_{a_1} & L_{a_3} \\ 0 & 0 & 0 & 0 & 0 & 0 & L_{a_2} & L_{a_4} \end{bmatrix} \middle| \begin{array}{c} L_{a_i} \in \{L_0, L_1, \dots, L_{m+1}\}; \\ 1 \le i \le 4 \end{array} \right\}$$

$$\subseteq V$$

be four DSm semivector subspaces of V of ordinary labels over the semifield S = {0, 1}.

Clearly $V = \bigcup_{i=1}^{4} W_i$ where

$$W_i \cap W_j = \begin{bmatrix} 0 & 0 & 0 & \dots & 0 \\ 0 & 0 & 0 & \dots & 0 \end{bmatrix} \text{ if } i \ne j.$$

Thus V is a direct sum of DSm semivector subspaces of ordinary labels of V over the semifield S = {0, 1}.

***Example 4.47:*** Let

$$V = \left\{ \begin{bmatrix} L_{a_1} & L_{a_2} & L_{a_3} \\ L_{a_4} & L_{a_5} & L_{a_6} \\ L_{a_7} & L_{a_8} & L_{a_9} \\ L_{a_{10}} & L_{a_{11}} & L_{a_{12}} \\ L_{a_{13}} & L_{a_{14}} & L_{a_{15}} \\ L_{a_{16}} & L_{a_{17}} & L_{a_{18}} \end{bmatrix} \middle| L_{a_i} \in \{L_0, L_1, L_2, \dots, L_m, L_{m+1}\}; 1 \le i \le 18 \right\}$$

be a DSm semivector space of ordinary labels over the semifield S = {0,1}. Consider the semivector subspaces $W_1$, $W_2$, $W_3$, $W_4$, $W_5$.

$$W_1 = \left\{ \begin{bmatrix} L_{a_1} & L_{a_2} & L_{a_3} \\ 0 & 0 & 0 \\ L_{a_4} & L_{a_5} & L_{a_6} \\ 0 & 0 & 0 \\ \vdots & \vdots & \vdots \\ 0 & 0 & 0 \end{bmatrix} \middle| L_{a_i} \in \{L_0, L_1, ..., L_{m+1}\}; 1 \leq i \leq 6 \right\} \subseteq V,$$

$$W_2 = \left\{ \begin{bmatrix} 0 & 0 & 0 \\ L_{a_1} & L_{a_2} & L_{a_3} \\ . & . & . \\ . & . & . \\ . & . & . \end{bmatrix} \middle| L_{a_i} \in \{L_0, L_1, ..., L_{m+1}\}; 1 \leq i \leq 3 \right\} \subseteq V,$$

$$W_3 = \left\{ \begin{bmatrix} 0 & 0 & 0 \\ 0 & 0 & 0 \\ 0 & 0 & 0 \\ L_{a_1} & L_{a_2} & L_{a_3} \\ 0 & 0 & 0 \\ 0 & 0 & 0 \end{bmatrix} \middle| L_{a_i} \in \{L_0, L_1, ..., L_{m+1}\}; 1 \leq i \leq 3 \right\} \subseteq V,$$

$$W_4 = \left\{ \begin{bmatrix} 0 & 0 & 0 \\ \vdots & \vdots & \vdots \\ 0 & 0 & 0 \\ L_{a_1} & L_{a_2} & L_{a_3} \\ 0 & 0 & 0 \end{bmatrix} \middle| L_{a_i} \in \{L_0, L_1, ..., L_{m+1}\}; 1 \leq i \leq 3 \right\} \subseteq V,$$

and

$$W_5 = \left\{ \begin{bmatrix} 0 & 0 & 0 \\ 0 & 0 & 0 \\ \vdots & \vdots & \vdots \\ 0 & 0 & 0 \\ L_{a_1} & L_{a_2} & L_{a_3} \end{bmatrix} \middle| L_{a_i} \in \{L_0, L_1, ..., L_{m+1}\}; 1 \leq i \leq 3 \right\} \subseteq V$$

be a DSm semivector subspaces of ordinary labels over the semifield S = {0, 1} of V.

We see $V = \bigcup_{i=1}^{5} W_i$ where

$$W_i \cap W_j = \begin{bmatrix} 0 & 0 & 0 \\ 0 & 0 & 0 \\ 0 & 0 & 0 \\ 0 & 0 & 0 \\ 0 & 0 & 0 \\ 0 & 0 & 0 \end{bmatrix} \text{ if } i \neq j; \ 1 \leq i, j \leq 5.$$

So V is a direct sum of DSm subsemivector spaces of V ordinary labels over the semifield S = {0, 1}.

***Example 4.48:*** Let

$$V = \left\{ \begin{bmatrix} L_{a_1} & L_{a_2} & L_{a_3} & L_{a_4} & L_{a_5} & L_{a_6} & L_{a_7} \\ L_{a_8} & L_{a_9} & L_{a_{10}} & L_{a_{11}} & L_{a_{12}} & L_{a_{13}} & L_{a_{14}} \\ L_{a_{15}} & L_{a_{16}} & L_{a_{17}} & L_{a_{18}} & L_{a_{19}} & L_{a_{20}} & L_{a_{21}} \end{bmatrix} \middle| \begin{array}{c} L_{a_i} \in \{L_0, L_1, ..., L_{m+1}\}; \\ 1 \leq i \leq 21 \end{array} \right\}$$

be a DSm semivector space of ordinary labels over the semifield S = {0, 1}. Consider $W_1$, $W_2$, $W_3$, $W_4$, $W_5$ and $W_6$ DSm semivector subspaces of ordinary labels of V over the semifield S = {0, 1}, where

$$W_1 = \left\{ \begin{bmatrix} 0 & L_{a_1} & 0 & L_{a_4} & 0 & 0 & L_{a_7} \\ 0 & L_{a_2} & 0 & L_{a_5} & 0 & 0 & 0 \\ 0 & L_{a_3} & 0 & L_{a_6} & 0 & 0 & 0 \end{bmatrix} \middle| \begin{array}{c} L_{a_i} \in \{L_0, L_1, ..., L_{m+1}\}; \\ 1 \leq i \leq 7 \end{array} \right\} \subseteq V$$

$$W_2 = \left\{ \begin{bmatrix} L_{a_1} & 0 & L_{a_5} & L_{a_6} & L_{a_7} & 0 & 0 \\ L_{a_2} & L_{a_4} & 0 & 0 & L_{a_8} & 0 & 0 \\ L_{a_3} & 0 & 0 & 0 & 0 & 0 & 0 \end{bmatrix} \middle| \begin{array}{c} L_{a_i} \in \{L_0, L_1, ..., L_{m+1}\}; \\ 1 \leq i \leq 8 \end{array} \right\} \subseteq V,$$

$$W_3 = \left\{ \begin{bmatrix} L_{a_1} & L_{a_2} & L_{a_5} & 0 & 0 & 0 & L_{a_9} \\ 0 & 0 & L_{a_6} & 0 & 0 & 0 & 0 \\ 0 & 0 & L_{a_7} & L_{a_8} & 0 & 0 & L_{a_{10}} \end{bmatrix} \middle| \begin{array}{l} L_{a_i} \in \{L_0, L_1, ..., L_{m+1}\}; \\ 1 \leq i \leq 10 \end{array} \right\} \subseteq V,$$

$$W_4 = \left\{ \begin{bmatrix} L_{a_1} & 0 & 0 & 0 & 0 & L_{a_6} & 0 \\ 0 & L_{a_4} & 0 & L_{a_5} & 0 & L_{a_7} & 0 \\ L_{a_2} & 0 & 0 & 0 & 0 & L_{a_8} & 0 \end{bmatrix} \middle| \begin{array}{l} L_{a_i} \in \{L_0, L_1, ..., L_{m+1}\}; \\ 1 \leq i \leq 8 \end{array} \right\} \subseteq V,$$

$$W_5 = \left\{ \begin{bmatrix} L_{a_1} & 0 & L_{a_2} & 0 & L_{a_7} & L_{a_8} & L_{a_4} \\ 0 & 0 & 0 & 0 & 0 & 0 & L_{a_6} \\ 0 & 0 & L_{a_3} & 0 & 0 & 0 & L_{a_5} \end{bmatrix} \middle| \begin{array}{l} L_{a_i} \in \{L_0, L_1, ..., L_{m+1}\}; \\ 1 \leq i \leq 8 \end{array} \right\} \subseteq V$$

and

$$W_6 = \left\{ \begin{bmatrix} L_{a_1} & L_{a_2} & 0 & 0 & L_{a_5} & 0 & L_{a_8} \\ 0 & 0 & L_{a_3} & 0 & L_{a_6} & 0 & 0 \\ 0 & 0 & 0 & L_{a_4} & L_{a_7} & 0 & L_{a_9} \end{bmatrix} \middle| \begin{array}{l} L_{a_i} \in \{L_0, L_1, ..., L_{m+1}\}; \\ 1 \leq i \leq 9 \end{array} \right\} \subseteq V$$

are DSm semivector subspaces of ordinary labels over the semifield $S = \{0, 1\}$.

Clearly $W_1 \cup W_2 \cup \ldots \cup W_6 = V$ that is $V = \bigcup_{i=1}^{6} W_i$ but

$$W_i \cap W_j \neq \begin{bmatrix} 0 & 0 & 0 & 0 & 0 & 0 & 0 \\ 0 & 0 & 0 & 0 & 0 & 0 & 0 \\ 0 & 0 & 0 & 0 & 0 & 0 & 0 \end{bmatrix}$$

if $i \neq j$; $1 \leq i, j \leq 6$.

Thus V is only a pseudo direct union of semivector subspaces $W_1$, $W_2$, $W_3$, $W_4$, $W_5$ and $W_6$ of V over the semifield $S = \{0, 1\}$.

Now we will work with the DSm special semivector spaces of ordinary labels over the semifield $S = \{L_0, L_1, \ldots, L_m, L_{m+1}\}$.

*Example 4.49:* Let

$$V = \left\{ \begin{bmatrix} L_{a_1} \\ L_{a_2} \\ L_{a_3} \\ L_{a_4} \end{bmatrix} \middle| L_{b_i} \in \{0 = L_0, L_1, L_2, ..., L_m, L_{m+1}\}; 1 \le i \le 4 \right\}$$

be a DSm special semivector space of ordinary labels over the semifield (semi lattice) $S = \{L_0, L_1, ..., L_m, L_{m+1}\}$. Consider

$$W_1 = \left\{ \begin{bmatrix} L_{a_1} \\ 0 \\ 0 \\ 0 \end{bmatrix} \middle| L_{b_1} \in S \right\} \subseteq V$$

the DSm special semivector subspace of ordinary labels over the semifield S. Now take

$$W_2 = \left\{ \begin{bmatrix} 0 \\ L_{a_1} \\ 0 \\ 0 \end{bmatrix} \middle| L_{b_1} \in S \right\} \subseteq V$$

the DSm special semivector subspace of ordinary labels over the semifield S.

$$W_3 = \left\{ \begin{bmatrix} 0 \\ 0 \\ L_{a_1} \\ 0 \end{bmatrix} \middle| L_{b_1} \in S \right\} \subseteq V$$

the DSm special semivector subspace of ordinary labels over the semifield S of V.

$$W_4 = \left\{ \begin{bmatrix} 0 \\ 0 \\ 0 \\ L_{a_1} \end{bmatrix} \middle| L_{b_1} \in S \right\} \subseteq V$$

the DSm special semivector subspace of V of ordinary labels over S.

We see $V = \bigcup_{i=1}^{4} W_i$ but

$$W_i \cap W_j \neq \begin{bmatrix} 0 \\ 0 \\ 0 \\ 0 \end{bmatrix} \text{ if } i \neq j;\ 1 \leq i, j \leq 4.$$

Thus V is the direct sum of the DSm special semivector subspaces $W_1$, $W_2$, $W_3$ and $W_4$.

*Example 4.50:* Let

$$V = \left\{ \begin{bmatrix} L_{a_1} & L_{a_3} & L_{a_5} & L_{a_7} & L_{a_9} & L_{a_{11}} \\ L_{a_2} & L_{a_4} & L_{a_6} & L_{a_8} & L_{a_{10}} & L_{a_{12}} \end{bmatrix} \middle| \begin{array}{l} L_{a_i} \in \{L_0, L_1, ..., L_{m+1}\}; \\ 1 \leq i \leq 12 \end{array} \right\}$$

be a DSm special semivector space of ordinary labels over the semifield $S = \{L_0, L_1, \ldots, L_m, L_{m+1}\}$. Consider the DSm special semivector subspaces $W_1$, $W_2$, $W_3$, $W_4$ and $W_5$ of V over S the semifield.

Here

$$W_1 = \left\{ \begin{bmatrix} L_{a_1} & L_{a_3} & 0 & 0 & 0 & 0 \\ L_{a_2} & 0 & 0 & 0 & 0 & 0 \end{bmatrix} \middle| \begin{array}{l} L_{a_i} \in \{L_0, L_1, ..., L_{m+1}\}; \\ 1 \leq i \leq 3 \end{array} \right\} \subseteq V,$$

is a DSm special semivector subspace of V over S.

$$W_2 = \left\{ \begin{bmatrix} 0 & 0 & L_{a_2} & 0 & 0 & 0 \\ 0 & L_{a_1} & L_{a_3} & 0 & 0 & 0 \end{bmatrix} \middle| \begin{array}{l} L_{a_i} \in \{L_0, L_1, ..., L_{m+1}\}; \\ 1 \leq i \leq 3 \end{array} \right\} \subseteq V$$

is again a DSm special semivector subspace of V over S.

$$W_3 = \left\{ \begin{bmatrix} 0 & 0 & 0 & L_{a_1} & 0 & 0 \\ 0 & 0 & 0 & L_{a_2} & 0 & 0 \end{bmatrix} \middle| \begin{array}{l} L_{a_i} \in \{L_0, L_1, ..., L_{m+1}\}; \\ 1 \leq i \leq 2 \end{array} \right\} \subseteq V$$

is a DSm special semivector subspace of V over S.

Consider

$$W_4 = \left\{ \begin{bmatrix} 0 & 0 & 0 & 0 & L_{a_1} & 0 \\ 0 & 0 & 0 & 0 & 0 & L_{a_2} \end{bmatrix} \middle| \begin{array}{l} L_{a_i} \in \{L_0, L_1, ..., L_{m+1}\}; \\ 1 \leq i \leq 2 \end{array} \right\} \subseteq V$$

is a DSm special semivector subspace of V over S.

Finally

$$W_5 = \left\{ \begin{bmatrix} 0 & 0 & 0 & 0 & 0 & L_{a_2} \\ 0 & 0 & 0 & 0 & L_{a_1} & 0 \end{bmatrix} \middle| \begin{array}{c} L_{a_i} \in \{L_0, L_1, ..., L_{m+1}\}; \\ 1 \le i \le 2 \end{array} \right\} \subseteq V$$

is a DSm special semivector subspace of V over S. Thus

$$V = \bigcup_{i=1}^{5} W_i$$

but

$$W_i \cap W_j \ne \begin{bmatrix} 0 & 0 & 0 & 0 & 0 & 0 \\ 0 & 0 & 0 & 0 & 0 & 0 \end{bmatrix}$$

if $i \ne j$; $1 \le i, j \le 5$. Thus V is a direct union of the DSm special semivector subspaces $W_1, W_2, ..., W_5$ of V over S.

*Example 4.51:* Let

$$V = \left\{ \begin{bmatrix} L_{a_1} & L_{a_2} & L_{a_3} \\ L_{a_4} & L_{a_5} & L_{a_6} \\ L_{a_7} & L_{a_8} & L_{a_9} \end{bmatrix} \middle| L_{a_i} \in \{L_0, L_1, L_2, ..., L_m, L_{m+1}\}; 1 \le i \le 9 \right\}$$

be a DSm special semivector space of ordinary labels over the semifield $S = \{L_0, L_1, ..., L_m, L_{m+1}\}$.

Consider

$$P_1 = \left\{ \begin{bmatrix} L_{a_1} & 0 & 0 \\ L_{a_2} & 0 & 0 \\ 0 & 0 & 0 \end{bmatrix} \middle| L_{a_i} \in \{L_0, L_1, L_2, ..., L_m, L_{m+1}\}; 1 \le i \le 2 \right\} \subseteq V$$

is a DSm special semi vector subspace of V over the semifield S.

Take

$$P_2 = \left\{ \begin{bmatrix} L_{a_1} & 0 & L_{a_2} \\ 0 & L_{a_3} & 0 \\ 0 & 0 & 0 \end{bmatrix} \middle| L_{a_i} \in \{L_0, L_1, L_2, ..., L_m, L_{m+1}\}; 1 \le i \le 3 \right\} \subseteq V$$

is a DSm special semivector subspace of V over the semifield S.

$$P_3 = \left\{ \begin{bmatrix} 0 & L_{a_1} & L_{a_2} \\ L_{a_3} & 0 & 0 \\ L_{a_4} & 0 & 0 \end{bmatrix} \middle| L_{a_i} \in \{L_0, L_1, L_2, ..., L_m, L_{m+1}\}; 1 \leq i \leq 4 \right\} \subseteq V$$

and

$$P_4 = \left\{ \begin{bmatrix} L_{a_1} & 0 & L_{a_2} \\ 0 & 0 & L_{a_3} \\ 0 & L_{a_4} & L_{a_5} \end{bmatrix} \middle| L_{a_i} \in \{L_0, L_1, L_2, ..., L_m, L_{m+1}\}; 1 \leq i \leq 5 \right\} \subseteq V$$

be special semivector subspace of V over the semifield S.

Consider

$$V = \bigcup_{i=1}^{4} P_i = P_1 \cup P_2 \cup P_3 \cup P_4$$

and

$$P_i \cap P_j \neq \begin{bmatrix} 0 & 0 & 0 \\ 0 & 0 & 0 \\ 0 & 0 & 0 \end{bmatrix} \text{ if } i \neq j; 1 \leq i, j \leq 4.$$

V is a pseudo direct sum of DSm special semivector subspaces of ordinary labels over the semifield $S = \{L_0, L_1, ..., L_m, L_{m+1}\}$.

Now having seen properties related with DSm semivector space of ordinary labels and DSm special semivector space of ordinary labels we now proceed onto give some of their applications.

**Chapter Five**

# APPLICATIONS OF DSM SEMIVECTOR SPACES OF ORDINARY LABELS AND REFINED LABELS

Study of DSm fields and Linear algebra of refined labels have been carried out by [34-5]. They have given applications while dealing qualitative information and other applications [7, 34-5].

We have introduced DSm semivector spaces of ordinary and refined labels and DSm set vector spaces of refined labels defined over a set. DSm group vector space of refined labels defined over a group and so on.

These structures have been transformed to matrices of refined labels and polynomials with refined label coefficients. These structures will find applications in fuzzy models, and mathematical model which uses matrices and in eigen value problems respectively.

Also this study can be used in web designing.

Further all types of social problems can be solved by using the partially ordered ordinary labels or unordered ordinary labels.

Further applications of these structures is to be invented in due course of time when these algebraic structures using these ordinary and refined labels are made more familiar with researchers.

**Chapter Six**

# SUGGESTED PROBLEMS

In this chapter we introduce 130 problem some of them are simple exercise and some difficult and a few of them are open research problems.

1. Let $V = \left\{ \begin{bmatrix} L_a & L_b \\ L_c & L_d \end{bmatrix} \middle| L_a, L_b, L_c, L_d \in L_R \right\}$ be a DSm linear algebra of refined labels over R (R reals).
   a. Find DSm linear subalgebras of V.
   b. Write V as a direct sum of subspaces.
   c. Define DSm linear operator on V which is non invertible.

2. Obtain some interesting properties about DSm vector spaces of refined labels over R (R reals).

3. Let $V = \left\{ \begin{bmatrix} L_{a_1} & L_{a_2} \\ L_{a_3} & L_{a_4} \\ L_{a_5} & L_{a_6} \end{bmatrix} \middle| L_{a_i} \in L_R ; 1 \leq i \leq 6 \right\}$ be a DSm vector space of refined labels over the field R.
   a. Find a basis of V.
   b. Find dimension of V
   c. Find atleast 3 vector subspaces of V.
   d. Find L(V, V)
   e. What is the algebraic structure enjoyed by L(V, V)?

4. Let $W = \left\{ \begin{bmatrix} L_{a_1} & L_{a_2} & L_{a_3} & L_{a_4} \\ L_{a_5} & L_{a_6} & L_{a_7} & L_{a_8} \end{bmatrix} \middle| L_{a_i} \in L_R ; 1 \leq i \leq 8 \right\}$ be a DSm vector space of refined labels over R.

   $V = \left\{ \begin{bmatrix} L_{a_1} & L_{a_2} & L_{a_3} \\ L_{a_4} & L_{a_5} & L_{a_6} \\ L_{a_7} & L_{a_8} & L_{a_9} \end{bmatrix} \middle| L_{a_i} \in L_R ; 1 \leq i \leq 9 \right\}$ be a DSm vector space of refined labels over R.
   a. Find a linear transformation T from V to W so that T is non invertible.
   b. Find a linear transformation T from W to V so that T is invertible.
   c. Find the algebraic structure enjoyed by L (V, W) and L (W, V).
   d. Find L (V, V).
   e. Find L (W, W).
   f. Define a projection on V.
   g. Write V as a direct sum.

5. Let $V = \left\{ \sum_{i=0}^{28} L_{a_i} x^i \middle| L_{a_i} \in L_R ; 0 \leq i \leq 28 \right\}$ be a DSm vector space of refined labels (polynomials with refined labels coefficients) over R.
   a. What is dimension of V?

b. Find a proper subset P of V which is linearly dependent.
c. Find a generating subset of V.
d. Find L (V, V)
e. Find subspaces of V.

6. Let $V = \left\{ \begin{bmatrix} L_{a_1} & L_{a_2} & L_{a_3} \\ L_{a_4} & L_{a_5} & L_{a_6} \\ L_{a_7} & L_{a_8} & L_{a_9} \end{bmatrix} \middle| L_{a_i} \in L_R ; 1 \le i \le 9 \right\}$ be a DSm linear algebra of refined labels over R.

a. What is the dimension of R?
b. Write V as a pseudo direct union of sublinear algebras.
c. Write V as a direct union of DSm sublinear algebras.
d. Is $W = \left\{ \begin{bmatrix} L_{a_1} & L_{a_2} & L_{a_3} \\ 0 & L_{a_4} & L_{a_5} \\ 0 & 0 & L_{a_6} \end{bmatrix} \middle| L_{a_i} \in L_R ; 1 \le i \le 6 \right\} \subseteq V$ a DSm vector subspace of V? What is dimension of W?
e. Let $W = \left\{ \begin{bmatrix} L_{a_1} & 0 & 0 \\ 0 & L_{a_2} & 0 \\ 0 & 0 & L_{a_3} \end{bmatrix} \middle| L_{a_i} \in L_R ; 1 \le i \le 3 \right\} \subseteq V$, be a DSm vector subspace of V over R. Define $\theta : V \to V$ so that W is invariant order θ.
f. Define a projection E on V. Is $E^2 = E$? Justify your answer.
g. Give a linear operator T on V which has non trivial nullity.

7. Let $V = \left\{ \begin{bmatrix} L_{a_1} \\ L_{a_2} \\ L_{a_3} \\ L_{a_4} \end{bmatrix}, \begin{pmatrix} L_{a_1} & L_{a_2} & L_{a_3} \end{pmatrix} \middle| L_{a_i} \in L_R; 1 \leq i \leq 4 \right\}$ be a set vector space of refined labels over the set $S = 3Z^+ \cup 4Z^+ \cup \{0\}$.
   a. Find a linear operator on V.
   b. Let $L(V, V) = \{T : V \to V\}$. What is the algebraic structure enjoyed by V.
   c. Write V as a direct sum.
   d. Write V as a pseudo direct sum.

8. Let $V = \left\{ \sum_{i=0}^{8} L_{a_i} x^i, \begin{bmatrix} L_{a_1} \\ L_{a_2} \\ \vdots \\ L_{a_9} \end{bmatrix}, \begin{pmatrix} L_{a_1} & \cdots & L_{a_{10}} \end{pmatrix}, \begin{bmatrix} L_{a_1} & L_{a_2} \\ L_{a_3} & L_{a_4} \end{bmatrix} \middle| L_{a_i} \in L_R; 1 \leq i \leq 10 \right\}$ be a set vector space of refined labels over the set $S = 8Z^+ \cup \{0\}$.
   a. Find a subset vector subspace of V over $T \subseteq S$ (T a subset of S).
   b. Can V have an invertible linear operator?
   c. Can V be written as a direct sum of set vector subspaces of refined labels over V.

9. Let $X = \left\{ \begin{bmatrix} L_{a_1} & L_{a_2} \\ L_{a_3} & L_{a_4} \end{bmatrix} \middle| L_{a_i} \in L_R; 1 \leq i \leq 4 \right\}$ be a DSm linear algebra of refined labels over the field R.
   a. Find a basis for X.
   b. Find dimension of X.
   c. Write X as a direct union of sublinear algebras.
   d. Can X be written as a pseudo direct sum of DSm linear subalgberas of refined labels of X over R?

e. Can X have subfield linear subalgebras of refined labels over R?
  f. Find a non invertible linear operator on X.
  g. If the field R is replaced by the field Q will X be a finite dimensional linear algebra.
  h. Find a projection $\eta$ on X. Is $\eta^2 = \eta$?

10. Let $W = \left\{ \begin{bmatrix} L_{a_1} & L_{a_2} & L_{a_3} \\ L_{a_4} & L_{a_5} & L_{a_6} \\ L_{a_7} & L_{a_8} & L_{a_9} \end{bmatrix} \middle| L_{a_i} \in L_R ; 1 \le i \le 9 \right\}$ be a DSm semigroup linear algebra of refined labels over the semigroup $S = Z^+ \cup \{0\}$.
  a. What is the dimension of W?
  b. Find a subsemigroup linear subalgebra of refined labels over the subsemigroup $T = 3Z^+ \cup \{0\}$.
  c. Write W as a direct sum of semigroup linear subalgebra of refined labels over S.
  d. Prove W can also be written as a pseudo direct sum of linear subalgebras of refined labels over S.
  e. Let $X = \left\{ \begin{bmatrix} L_{a_1} & L_{a_2} & L_{a_3} \\ 0 & 0 & 0 \\ 0 & 0 & 0 \end{bmatrix} \middle| L_{a_i} \in L_R ; 1 \le i \le 3 \right\} \subseteq W$ is a semigroup linear subalgebra of W of refined labels over the semigroup S. Define a projection $\eta$ : W → X so that $\eta(X) \subseteq X$ and $\eta^2 = \eta$.

11. Let $V = \left\{ \begin{bmatrix} L_{a_1} & L_{a_2} & L_{a_3} \\ L_{a_4} & L_{a_5} & L_{a_6} \\ L_{a_7} & L_{a_8} & L_{a_9} \end{bmatrix} \middle| L_{a_i} \in L_R ; 1 \le i \le 9 \right\}$ be a group linear algebra over the group G = R.
  a. Find a basis of V.
  b. What is the dimension of V over R?
  c. If R is replaced by Q or Z what will the dimensions?
  d. Find the algebraic structure enjoyed by $L_R(V, V)$.

12. Let $V = \left\{ \begin{bmatrix} L_{a_1} \\ L_{a_2} \\ \vdots \\ L_{a_4} \end{bmatrix}, \left( L_{a_1} \quad L_{a_2} \quad L_{a_3} \right), \right.$

$\left. \begin{pmatrix} L_{a_1} & L_{a_2} & L_{a_3} & L_{a_4} & L_{a_5} \\ L_{a_6} & L_{a_7} & L_{a_8} & L_{a_9} & L_{a_{10}} \\ L_{a_{11}} & L_{a_{12}} & L_{a_{13}} & L_{a_{14}} & L_{a_{15}} \\ L_{a_{16}} & L_{a_{17}} & L_{a_{18}} & L_{a_{19}} & L_{a_{20}} \end{pmatrix} \middle| L_{a_i} \in L_R ; 1 \le i \le 20 \right\}$ be a

set vector space of refined labels over the set $S = Z^+ \cup \{0\}$.
   a. Find the dimension of V over $Z^+ \cup \{0\} = S$.
   b. Find set vector subspaces of refined labels of V over $S = Z^+ \cup \{0\}$.
   c. Write V as a direct union of set vector subspaces of V over S.
   d. Can V be written as a pseudo direct union of set vector subspaces of V?

13. Find some interesting properties enjoyed by set linear algebras of refined labels over a set S.

14. Let $W = \left\{ \left( L_{a_1} \quad L_{a_2} \quad \ldots \quad L_{a_{10}} \right) \middle| L_{a_i} \in L_R ; 1 \le i \le 10 \right\}$ be a DSm semigroup linear algebra of refined labels over the semigroup $S = R$ (reals).
   a. Find a basis for W.
   b. What is the dimension of W over R?
   c. Write W as a direct union of semigroup linear subalgebras of refined labels over the semigroup $S = R$.
   d. If $S = R$ is replaced by the semigroup $T = R^+ \cup \{0\}$ study (i) and (ii)
   e. Study (i) and (ii) if R is replaced by $3Z^+ \cup \{0\} = P$.

15. Let $M = \left\{ \begin{bmatrix} L_{a_1} & L_{a_2} & L_{a_3} \\ L_{a_4} & L_{a_5} & L_{a_6} \\ L_{a_7} & L_{a_8} & L_{a_9} \end{bmatrix} \middle| L_{a_i} \in L_R ; 1 \leq i \leq 9 \right\}$ be a group linear algebra over the group G = R (R reals).
   a. Find dimension of M over R.
   b. Find a basis of M over R.
   c. Does M have a linearly independent set over R?
   d. Find a linearly dependent set of order less than the order of the basis over R.
   e. Study (i) to (v) if R is replaced by Z.

16. Let $V = \left\{ \begin{bmatrix} L_{a_1} & L_{a_2} & L_{a_3} \\ L_{a_4} & L_{a_5} & L_{a_6} \\ \vdots & \vdots & \vdots \\ L_{a_{28}} & L_{a_{29}} & L_{a_{30}} \end{bmatrix} \middle| L_{a_i} \in L_R ; 1 \leq i \leq 30 \right\}$ be a DSm group linear algebra refined labels over the group G = Z.
   a. What is dimension of V?
   b. Is V finite dimensional?
   c. Can V be written as a direct union of DSm group linear subalgebras of refined labels over R?
   d. Study $L_z (V, V)$.

17. Obtain some interesting properties enjoyed by semigroup vector spaces V of refined labels over the semigroup S = $5Z^+ \cup \{0\}$. If S is replaced by Z what is the special features enjoyed by V. Study the problem if S is replaced by R (the reals).

18. Let P = $\left\{ \begin{bmatrix} L_{a_1} & L_{a_2} \\ L_{a_3} & L_{a_4} \\ L_{a_5} & L_{a_6} \\ L_{a_7} & L_{a_8} \\ L_{a_9} & L_{a_{10}} \end{bmatrix}, \begin{bmatrix} L_{a_1} & L_{a_2} \\ L_{a_3} & L_{a_4} \end{bmatrix}, (L_{a_1}\ L_{a_2}\ \ldots\ L_{a_{12}}) \middle| \begin{array}{l} L_{a_i} \in L_R; \\ 1 \le i \le 12 \end{array} \right\}$

be a group vector space of refined labels over the group G = Z.
   a. Find a basis for P.
   b. What is the dimension of P over Z?
   c. What is the dimension of V if Z is replaced by Q?
   d. What is the dimension of V if Z is replaced by R?
   e. Find T : P → P such that T is an idempotent linear operator on V.

19. Let V = {All 5 × 5 matrices with entries from $L_R$}. Is V a DSm linear algebra over R?
    Is V a finite dimensional linear algebra over Q?
    Study the properties and compare the structure of V defined on R or Q.

20. Obtain some interesting properties related with group linear algebra of refined labels defined over the group G = Z.

21. Let W = $\left\{ \begin{bmatrix} L_{a_1} & L_{a_2} \\ L_{a_3} & L_{a_4} \\ \vdots & \vdots \\ L_{a_{21}} & L_{a_{22}} \end{bmatrix} \middle| L_{a_i} \in L_R; 1 \le i \le 22 \right\}$ be a DSm

group linear algebra of refined labels over the group G = Z.
   a. Find the dimension of W over Z.
   b. Write W as a direct union of DSm group linear subalgebras of refined labels over the group Z.
   c. Write W as a pseudo direct sum of DSm group linear subalgebras.

d. Find T : W → W which has a non trivial kernel.

22. Let $V = \left\{ \sum_{i=0}^{\infty} L_{a_i} x^i \mid L_{a_i} \in L_R; 0 \le i \le \infty \right\}$ be the DSm linear algebra of polynomials with refined label over the field R.
   a. Find DSm linear subalgebra of V over R.
   b. Find a basis for V.
   c. Prove V is infinite dimensional.
   d. Can we write V as a direct union of DSm linear subalgebras?

23. Let $V = \left\{ \begin{bmatrix} L_{a_1} \\ L_{a_2} \\ L_{a_3} \\ \vdots \\ L_{a_8} \end{bmatrix} \middle| L_{a_i} \in L_R; 1 \le i \le 8 \right\}$ be a DSm linear algebra of refined labels over the field R of reals.
   a. Find the basis of V.
   b. Write V as a direct sum.
   c. Write V as a pseudo direct sum of DSm sublinear algebras of refined labels.
   d. Find L (V, V).

24. Find some nice applications of DSm group linear algebras over the group G.

25. Obtain some interesting result about DSm semigroup vector spaces over the semigroup 53Z.

26. Let
$$V = \left\{ \begin{bmatrix} L_{a_1} & L_{a_2} \\ L_{a_3} & L_{a_4} \\ L_{a_5} & L_{a_6} \\ L_{a_7} & L_{a_8} \end{bmatrix}, (L_{a_1}\ L_{a_2}\ L_{a_3}), \begin{bmatrix} L_{a_1} & L_{a_2} & L_{a_3} & L_{a_4} & L_{a_5} \\ L_{a_6} & L_{a_7} & L_{a_8} & L_{a_9} & L_{a_{10}} \\ L_{a_{11}} & L_{a_{12}} & L_{a_{13}} & L_{a_{14}} & L_{a_{15}} \end{bmatrix} \middle| \begin{array}{l} L_{a_i} \in L_R; \\ 1 \le i \le 15 \end{array} \right\}$$
be a DSm set vector space of refined labels over the set S = $3Z^+ \cup 5Z^+ \cup 7Z^+ \cup \{0\}$ and
$$W = \left\{ [L_{a_1}\ L_{a_2}\ ...\ L_{a_8}], \begin{pmatrix} L_{a_1} & L_{a_2} \\ 0 & L_{a_3} \end{pmatrix}, \begin{bmatrix} L_{a_1} & L_{a_2} & L_{a_3} \\ L_{a_4} & L_{a_5} & L_{a_6} \\ L_{a_7} & L_{a_8} & L_{a_9} \\ L_{a_{10}} & L_{a_{11}} & L_{a_{12}} \\ L_{a_{13}} & L_{a_{14}} & L_{a_{15}} \end{bmatrix} \middle| \begin{array}{l} L_{a_i} \in L_R; \\ 1 \le i \le 15 \end{array} \right\}$$ be
a DSm set vector space of refined labels over the set S = $3Z^+ \cup 5Z^+ \cup 7Z^+ \cup \{0\}$.
   a. Let $L_S$ (V, W) be the collection of linear transformation from V into W. Study the structure of $L_S$ (V, W).
   b. Find direct sum of V and pseudo direct sum of W.
   Define T : V → W such that T preserves set subvector spaces.

27. Give an example of pseudo DSm semigroup linear subalgebra of refined labels of a DSm group linear algebra.

28. Does there exist a DSm group linear algebra of refined labels which has no pseudo DSm semigroup linear subalgebra of refined labels? Justify your claim.

29. Does there exist a DSm group linear algebra of refine labels which is a simple DSm group linear algebra? Justify your claim.

30. What are the essential difference between a DSm set linear algebra of refined labels and DSm group linear algebra of refined labels?

31. Let $W = \left\{ \begin{bmatrix} L_{a_1} & L_{a_2} & L_{a_3} & L_{a_4} \\ L_{a_5} & L_{a_6} & L_{a_7} & L_{a_8} \\ L_{a_9} & L_{a_{10}} & L_{a_{11}} & L_{a_{12}} \end{bmatrix} \middle| L_{a_i} \in L_R ; 1 \leq i \leq 12 \right\}$

be a DSm group linear algebra of refined labels over the group $Z = G$.

a. Is $P = \left\{ \begin{bmatrix} L_{a_1} & L_{a_2} & 0 & 0 \\ L_{a_3} & L_{a_4} & 0 & 0 \\ L_{a_5} & L_{a_6} & 0 & 0 \end{bmatrix} \middle| L_{a_i} \in L_R ; 1 \leq i \leq 6 \right\} \subseteq W$;

a group linear subalgebra of refined labels of W over $G = Z$.

b. Define $T : W \to W$ so that $T(P) \subseteq P$.
c. Define $T : W \subseteq W$ so that $T(P) \not\subseteq P$
d. Let $M = \left\{ \begin{bmatrix} L_{a_1} & 0 & L_{a_2} & 0 \\ 0 & L_{a_3} & 0 & L_{a_4} \\ L_{a_5} & 0 & L_{a_6} & 0 \end{bmatrix} \middle| L_{a_i} \in L_R ; 1 \leq i \leq 6 \right\} \subseteq W$. Is

M a group linear subalgebra of refined labels over $G = Z$ of W?

e. Find a T such that $T(M) \subseteq M$.
f. Does there exists a $T : W \subseteq W$ so that $T(M) \subseteq M$ and $T(P) \subseteq P$?
g. Find a $T : W \to M$ so that $T^2 = T$.
h. Will T in (vii) be such that $T(P) \subseteq P$?

32. Let $M = \left\{ \begin{bmatrix} L_{a_1} & L_{a_2} & L_{a_3} \\ L_{a_4} & L_{a_5} & L_{a_6} \\ L_{a_7} & L_{a_8} & L_{a_9} \end{bmatrix} \middle| L_{a_i} \in L_R; 1 \leq i \leq 9 \right\}$ be a DSm semigroup linear algebra of refined labels over the semigroup $S = R^+ \cup \{0\}$.

  a. If S is replaced by $Z^+ \cup \{0\}$ then what are the properties enjoyed by M using $R^+ \cup \{0\}$.

  b. $N = \left\{ \begin{bmatrix} 0 & L_{a_1} & 0 \\ L_{a_2} & 0 & 0 \\ L_{a_3} & L_{a_4} & 0 \end{bmatrix} \middle| L_{a_i} \in L_R; 1 \leq i \leq 4 \right\} \subseteq M$,

  Find $T : M \to M$ such that $T(N) \subseteq N$.

  c. Find a $T : M \to M$ so that $T^2 = T$.

  d. Let $D = \left\{ \begin{bmatrix} 0 & 0 & 0 \\ L_{a_1} & 0 & L_{a_2} \\ 0 & L_{a_3} & 0 \end{bmatrix} \middle| L_{a_i} \in L_R; 1 \leq i \leq 3 \right\} \subseteq M$,

  find a $T : M \to M$ so that $T(D) \not\subseteq D$.

33. Let

  $P = \left\{ \begin{bmatrix} L_{a_1} & L_{a_2} \\ L_{a_3} & L_{a_4} \\ L_{a_5} & L_{a_6} \end{bmatrix}, \begin{bmatrix} L_{a_1} \\ L_{a_2} \\ \vdots \\ L_{a_9} \end{bmatrix}, \begin{pmatrix} L_{a_1} & L_{a_2} & \ldots & L_{a_{20}} \end{pmatrix} \middle| \begin{array}{l} L_{a_i} \in L_R; \\ 1 \leq i \leq 20 \end{array} \right\}$

  be a set vector space of refined labels over the set $S = 5Z^+ \cup 3Z^+ \cup \{0\}$.

  a. Define a $T : P \to P$ such that $T^2 = T$.
  b. Write P as a direct sum.
  c. Write P as a pseudo direct sum of set vector subspaces of refined labels.

d. Let $K = \left\{ \begin{bmatrix} L_{a_1} & L_{a_2} \\ L_{a_3} & L_{a_4} \\ L_{a_5} & L_{a_6} \end{bmatrix}, \begin{bmatrix} L_{a_1} \\ L_{a_2} \\ \vdots \\ L_{a_9} \end{bmatrix} \middle| L_{a_i} \in L_R; 1 \le i \le 9 \right\} \subseteq P$.

   Is K a set vector subspace of V over S?
   e. Find a $T : P \to P$ so that $T(K) \subseteq K$.
   f. Find a $T : P \to P$ so that $T(K) \not\subseteq K$.

34. Let $V = \{ (L_{a_1}, 0, L_{a_2}, L_{a_3}, L_{a_4}), (0, L_{a_1}, L_{a_2}, 0, 0),$

    $(L_{a_1}, L_{a_2}, 0, 0, L_{a_3}), (0, 0, 0, L_{a_3}, L_{a_4}) \middle| L_{a_i} \in L_R; 1 \le i \le 4 \}$

    be a set vector space of refined labels over the set $S = 3Z^+ \cup \{0\}$.
    a. Find set vector subspaces of V over S.
    b. Prove V is not a set linear algebra of refined labels over R.
    c. Define $T : V \to V$ so that $T^2 = T$.
    d. Define $T : V \to V$ so that $T^{-1}$ exists.

35. Let $V = \{$all $8 \times 8$ matrices with entries from $L_R\}$ be a group linear algebra of refined labels over the group $G = Z$.
    a. Define $T : V \to V$ so that $T^2 = T$.
    b. Write V as a direct sum of group linear subalgebra of refined labels over G.
    c. Is $M = \{$All $8 \times 8$ upper triangular matrices with entries from $L_R\} \subseteq V$ a group linear subalgebra of refined labels over G?
    d. Let $K = $
       $\subseteq V$, Is K a group linear subalgebra of refined labels over G?
    e. Define $T : V \to V$ such that $T(K) \subseteq K$.
    f. Define $T : V \to V$ such that $T^{-1} : V \to V$ exists.

36. Let $K = \left\{ \begin{bmatrix} L_{a_1} \\ L_{a_2} \\ \vdots \\ L_{a_{12}} \end{bmatrix}, (L_{a_1}, L_{a_2}, L_{a_3}), \begin{bmatrix} L_{a_1} & L_{a_2} \\ L_{a_3} & L_{a_4} \\ L_{a_5} & L_{a_6} \end{bmatrix} \middle| \begin{array}{l} L_{a_i} \in L_R; \\ 1 \le i \le 12 \end{array} \right\}$ be

   a group vector space of refined labels over the group G = R.
   a. Prove K is not a group linear algebra of refined labels over the group G.
   b. Let $M = \left\{ (L_{a_1}, L_{a_2}, L_{a_3}), \begin{bmatrix} L_{a_1} \\ L_{a_2} \\ \vdots \\ L_{a_{12}} \end{bmatrix} \middle| L_{a_i} \in L_R; 1 \le i \le 12 \right\}$

   ⊆ K, be the group vector subspace of K of refined labels over R.
   Define T : K → M so that T (M) ⊆ M
   c. Find T : K → K so that $T^2 = T$.
   d. Find a T so that $T^{-1}$ does not exist.

37. Let
   $T = \left\{ \begin{bmatrix} L_{a_1} \\ L_{a_2} \\ L_{a_3} \\ L_{a_4} \end{bmatrix}, (L_{a_1}, L_{a_2}), \begin{bmatrix} L_{a_1} & L_{a_2} \\ L_{a_3} & L_{a_4} \\ L_{a_5} & L_{a_6} \end{bmatrix}, \begin{pmatrix} L_{a_1} & \cdots & L_{a_{12}} \\ L_{a_{13}} & \cdots & L_{a_{24}} \\ L_{a_{25}} & \cdots & L_{a_{36}} \end{pmatrix} \middle| \begin{array}{l} L_{a_i} \in L_R; \\ 1 \le i \le 36 \end{array} \right\}$ be

   a semigroup vector space of refined labels over the semigroup $S = Q^+ \cup \{0\}$.
   a. Write T as a direct union of semigroup vector subspaces over S.
   b. Find $L_S$ (T, T).

38. If T be as in problem (37) and if $W = \left\{ \begin{bmatrix} L_{a_1} \\ L_{a_2} \\ L_{a_3} \\ L_{a_4} \end{bmatrix} \middle| \begin{array}{l} L_{a_i} \in L_R; \\ 1 \le i \le 4 \end{array} \right\}$

$\subseteq$ T. Define $\eta : T \to W$ so that $\eta$ is a projection.
a. Will $\eta^2 = \eta$?
b. Is $\eta(W) \subseteq W$?
c. If P is any other subspace find $\eta(P)$.

39. Let

$$M = \left\{ \begin{bmatrix} L_{a_1} & L_{a_2} & L_{a_3} \\ L_{a_4} & L_{a_5} & L_{a_6} \\ L_{a_7} & L_{a_8} & L_{a_9} \end{bmatrix}, \begin{bmatrix} L_{a_1} \\ L_{a_2} \\ \vdots \\ L_{a_9} \end{bmatrix}, \begin{pmatrix} L_{a_1} & L_{a_2} & \cdots & L_{a_{10}} \\ L_{a_{11}} & L_{a_{12}} & \cdots & L_{a_{20}} \end{pmatrix} \middle| \begin{array}{l} L_{a_i} \in L_R; \\ 1 \le i \le 20 \end{array} \right\}$$

and $V = \left\{ (L_{a_1}, \ldots, L_{a_8}), \begin{bmatrix} L_{a_1} & L_{a_2} & L_{a_3} \\ L_{a_4} & L_{a_5} & L_{a_6} \\ L_{a_7} & L_{a_8} & L_{a_9} \\ L_{a_{10}} & L_{a_{11}} & L_{a_{12}} \\ L_{a_{13}} & L_{a_{14}} & L_{a_{15}} \end{bmatrix}, \begin{bmatrix} L_{a_1} & L_{a_2} \\ L_{a_3} & L_{a_4} \end{bmatrix} \middle| \begin{array}{l} L_{a_i} \in L_R; \\ 1 \le i \le 15 \end{array} \right\}$

be any two group vector spaces of refined labels over the group G = Z.
a. Find a basis of M.
b. Find a basis of V.
c. Find $T : M \to V$ so that ker $T \ne \phi$.
d. Find $S : V \to M$ with ker $S \ne \phi$.
e. Is $L_Z(V, M) \cong L_Z(M, V)$?
f. Write $M = \cup M_i$, $M_i$ group vector subspaces of M over Z.
g. Find $\theta : V \to M$ so that $\theta^{-1}$ exist.
h. Find $h : M \to V$ so that $h^{-1}$ is not defined.

40. Let $P = (L_R[x], +, \times)$ the polynomial ring with refined label coefficients from R.

a. Is $L_R[x] \cong R[x]$?
b. Can $L_R[x]$ have ideals?
c. Is $L_R[x]$ a PID?
d. Can $L_R[x]$ be a S-ring?

41. Let $L_R[x_1, x_2]$ be the polynomial ring in the variables $x_1$ and $x_2$ with refined label coefficients.
    a. Prove $L_R[x_1, x_2] \cong R[x_1, x_2]$.
    b. Is $L_R[x_1, x_2]$ a vector space over the reals R?
    c. Is $L_R[x_1, x_2]$ a linear algebra over the reals R?
    d. Is $L_R[x_1, x_2]$ a principal ideal ring?

42. Is $L_R[x]$ a Euclidean ring?

43. Is $L_R[x]$ a principal ideal ring?

44. Can any polynomial $L_R[x]$ be written in a unique manner as a product of irreducible polynomials in $L_R[x]$?

45. Let $A = (p(x))$ in $L_R[x]$ be an ideal. Is it true A is a maximal ideal if and only if $p(x)$ is irreducible over $L_R$?

46. Is $L_R$ an integral domain?

47. Will $L_R[x_1, x_2, x_3]$ be an integral domain?

48. Can we derive the division algorithm in case of $L_R[x]$?

49. Let $L_R[x]$ be the polynomial ring. The ideal $A = (a_0)$ is a maximal ideal of the ring $L_R[x]$ if and only if $a_0$ is a prime element of $L_R$ ($L_R \cong R$).

50. Can we say $L_R[x]$ is a unique factorization domain?

51. Let $L_R[x]$ be the polynomial ring in the variable x. For $f(x) \in L_R[x]$ define $f'(x)$ the derivative of the polynomial. Prove if $f(x) \in L_R[x]$ where $L_R$ is the field of refined labels then $f(x)$ is divisible by the square of a polynomial if and only if f (x) and f′ (x) have a greatest common divisor d (x) of positive degree.

52. Let V be a vector space of refined labels over the reals R. W be a vector subspace of refined labels of V over R.
    a. Define quotient space.
    b. If T : V → U is a vector space of refined labels over R. (where U and V are vector space of refined labels over R) with kernel T = W then V is isomorphic to U/W and prove there exists a homomorphism of U onto U/W.

53. Let V be a vector space of refined labels over R (or $L_R$). Let $W_1, W_2, \ldots, W_t$ be vector subspaces of refined labels over R (or $L_R$). Let V be direct union or direct sum of $W_1, W_2, \ldots, W_t$. Can we develop the notion of external direct sum using $W_1, \ldots, W_t$.

54. Let S ⊆ V, V a vector space of refined labels over R (or $L_R$) L (S) be the linear span of S.
    a. Prove L (S) is a subspace of refined labels of V.
    b. If T and S are subsets of V then
        i) T ⊆ S implies L (T) ⊆ L (S).
        ii) L (S ∪ T) = L (S) + L (T).
        iii) L (L (S)) = L (S).

55. Let V be a finite dimensional vector space of refined labels over R and if W is a subspace of refined labels of V over R then
    a. W is finite dimensional.
    b. dim V ≥ dim W.
    c. dim V/W = dim V – dim W.

56. If A and B are finite dimensional vector sub spaces of refined labels of V over R then;
    a. A + B is finite dimensional.
    b. dim (A + B) = dim A + dim B – dim (A ∩ B).

57. Let V and W be any two vector spaces of refined labels over the field R.
    Is $Hom_R (V, W)$ a vector space over R? (justify your answer)

58. If V and W are of dimensions m and n over R (V and W vector spaces of refined labels over R) then will $Hom (V, W) = L_R (V, W)$ be of dimension mn over R?

59. Let $V = \{(L_{a_1}, L_{a_2}, ..., L_{a_8}) \mid L_{a_i} \in L_R; 1 \leq i \leq 8\}$ be a vector space of refined labels over R.
    a. Prove V has orthogonal vector u, v such that $u.v = (L_{a_1}, L_{a_2}, ..., L_{a_8}) \cdot (L_{b_1}, L_{b_2}, ..., L_{b_8}) = (0, 0, ..., 0)$.
    b. If W is a subspace of refined labels of V, find $W^\perp$. Will $(W^\perp)^\perp = W$?

60. Let $V = \left\{ \begin{pmatrix} L_{a_1} & L_{a_2} & ... & L_{a_9} \\ L_{a_{10}} & L_{a_{11}} & ... & L_{a_{18}} \end{pmatrix} \middle| L_{a_i} \in L_R; 1 \leq i \leq 18 \right\}$ be a DSm vector space of defined labels over the field R. How many hypersubspaces (hyperspace) in V exists?

61. Let V be a vector space of refined labels over the field R. f is a linear functional on V. Study properties about V.

62. Derive Taylors formula for $L_R$.

63. Define minimal polynomial for any linear operator on a finite dimensional DSm vector space of refined labels over the field F.

64. Can Cayley-Hamilton theorem be derived in case of linear operator on vector space of refined labels over R?

65. Define invariant direct sums in case of a vector space of refined labels over R.

66. Let T be a linear operator on a finite dimensional vector space V of refined labels over R. Find conditions on T so that T is diagonalizable.

67. Can Primary – decomposition Theorem be derived in case of linear operator T on V?

68. Let $V = \left\{ \begin{pmatrix} L_{a_1} & L_{a_2} & L_{a_3} & L_{a_4} \\ L_{a_5} & L_{a_6} & L_{a_7} & L_{a_8} \\ L_{a_9} & L_{a_{10}} & L_{a_{11}} & L_{a_{12}} \\ L_{a_{13}} & L_{a_{14}} & L_{a_{15}} & L_{a_{16}} \end{pmatrix} \middle| L_{a_i} \in L_R; 1 \leq i \leq 16 \right\}$

be a DSm linear algebra of refined labels over the field R.
   a. Write W as a direct sum.
   b. Find two subspaces of V which are not disjoint.
   c. Find dimension of V.
   d. Find a basis for V.
   e. Find $T : V \to V$ such that ker $T = \phi$.

69. Show that the space V of refined labels over R where $V = \left\{ \begin{pmatrix} L_{a_1} & L_{a_2} & L_{a_3} \\ L_{a_4} & L_{a_5} & L_{a_6} \\ L_{a_7} & L_{a_8} & L_{a_9} \end{pmatrix} \middle| L_{a_i} \in L_R; 1 \leq i \leq 9 \right\}$ is a linear algebra over R.

70. Let
$M = \left\{ \begin{bmatrix} L_{a_1} \\ L_{a_2} \\ L_{a_3} \\ L_{a_4} \end{bmatrix}, (L_{a_1}, L_{a_2}, ..., L_{a_{10}}), \begin{bmatrix} L_{a_1} & L_{a_2} & L_{a_3} \\ L_{a_4} & L_{a_5} & L_{a_6} \\ L_{a_7} & L_{a_8} & L_{a_9} \end{bmatrix} \middle| \begin{array}{l} L_{a_i} \in L_R; \\ 1 \leq i \leq 9 \end{array} \right\}$

be the DSm vector space of refined labels over R.
   a. Prove M is not a DSm linear algebra.
   b. Find a basis for M.
   c. Write $M = \cup M_i$ as a direct sum $W_i$'s vector subspaces of M.

d. Write $M = \cup W_i$ as pseudo direct sum.
  e. Define $T : M \to W$ where $W = \left\{ \begin{bmatrix} L_{a_1} \\ L_{a_2} \\ L_{a_3} \\ L_{a_4} \end{bmatrix} \middle| L_{a_i} \in L_R; 1 \leq i \leq 4 \right\} \subseteq$

  M is a subspace such that $T(W) \subseteq W$ and $T^2 = T$.

71. Let $X = \left\{ \begin{bmatrix} L_{a_1} & L_{a_2} & L_{a_3} \\ L_{a_4} & L_{a_5} & L_{a_6} \\ L_{a_7} & L_{a_8} & L_{a_9} \end{bmatrix} \middle| L_{a_i} \in L_R; 1 \leq i \leq 9 \right\}$ be a DSm linear algebra of refined labels over R.
    a. Find $T : X \to X$ such that T is non invertible.
    b. Find $T : X \to X$ such that $T^2 = T$.
    c. Let $W = \left\{ \begin{bmatrix} L_{a_1} & 0 & 0 \\ 0 & L_{a_2} & 0 \\ 0 & 0 & L_{a_3} \end{bmatrix} \middle| L_{a_i} \in L_R; 1 \leq i \leq 3 \right\} \subseteq X$

    be a DSm linear subalgebra of X. Find a $T : X \to X$ such that $T(W) \subseteq W$.

72. Show the set of unordered or partially ordered refined labels forms a lattice.

73. Give nice results about refined label lattice.

74. Give some important / interesting properties related with DSm semivector space of refined labels.

75. What is the major difference between the DSm vector space of refined labels and DSm semivector space of refined labels.

76. Let $V = \left\{ \begin{bmatrix} L_{a_1} \\ L_{a_2} \\ L_{a_3} \\ L_{a_4} \\ L_{a_5} \end{bmatrix}, (L_{a_1}, L_{a_2}, L_{a_3}), \begin{bmatrix} L_{a_1} & L_{a_2} \\ 0 & L_{a_3} \end{bmatrix} \middle| \begin{array}{l} L_{a_i} \in L_R; \\ 1 \leq i \leq 5 \end{array} \right\}$ be

a DSm vector space of refined labels over the field R.
   a. Is V finite dimensional over R?
   b. Find a basis of V.
   c. Is V a linear algebra of refined labels?
   d. Find subspaces $W_i$ of V so that $V = \bigcup_i W_i$.
   e. Write $V = \cup W_i$ as a pseudo direct sum of subspaces.

77. Let $V = \left\{ \begin{bmatrix} L_{a_1} & L_{a_2} \\ L_{a_3} & L_{a_4} \\ L_{a_5} & L_{a_6} \end{bmatrix}, (L_{a_1}, L_{a_2}, ..., L_{a_7}) \middle| \begin{array}{l} L_{a_i} \in L_R; \\ 1 \leq i \leq 7 \end{array} \right\}$ be a

DSm vector space of refined labels over R.
   a. Let $W_1 =$
   $\left\{ \begin{bmatrix} L_{a_1} & 0 \\ 0 & L_{a_2} \\ 0 & L_{a_3} \end{bmatrix}, \begin{bmatrix} 0 & 0 \\ L_{a_1} & L_{a_2} \\ 0 & L_{a_3} \end{bmatrix}, (L_{a_1}, 0, 0, 0, 0, L_{a_2}, L_{a_3}), (0, L_{a_1}, L_{a_2}, L_{a_3}, 0, 0, 0, 0) \right.$

   $\left. L_{a_i} \in L_R; 1 \leq i \leq 3 \right\} \subseteq V$ be a DSm vector subspaces
   refined labels of V over R.
   i) Define $T : V \to V$ so that $T(W_1) \subseteq W_1$.
   ii) $T(W_1) \not\subseteq W_1$
   b. Write $V = \cup W_i$ as pseudo sum of vector subspaces of refined labels.

78. Let $V = \left\{ \sum_{i=0}^{\infty} L_{a_i} x^i \,\middle|\, L_{a_i} \in R \right\}$ be the DSm vector space of refined labels over R. Study all the properties related with V.

79. Let $V = \left\{ \begin{bmatrix} L_{a_1} & L_{a_2} \\ L_{a_3} & L_{a_4} \\ L_{a_5} & L_{a_6} \\ L_{a_7} & L_{a_8} \end{bmatrix} \,\middle|\, L_{a_i} \in L_R; 1 \le i \le 8 \right\}$ be a group linear algebra of refined labels over R.
    a. Find subspaces of refined labels of V.
    b. Is V simple?
    c. Find a basis for V.
    d. What is dimension of V if R is replaced by Q?

80. Let $M = \begin{bmatrix} L_a & L_b & L_c \\ L_d & L_e & L_f \\ L_g & L_h & L_i \end{bmatrix}$ be a refined label matrix find the eigen values and eigen vectors by taking a suitable m + 1 and a, b, c, d, e, f, g, h, i.

81. Construct a model using refined label matrix to study any social problem.

82. When will the labels be analysed using FCM or FRM? Illustrate this situation by some examples.

83. Let $V = \left\{ \begin{bmatrix} L_{a_1} & L_{a_2} \\ L_{a_3} & L_{a_4} \end{bmatrix} \,\middle|\, L_{a_i} \in L_{R^+ \cup \{0\}}; 1 \le i \le n \right\}$ be a semivector space of refined labels over the semifield $S = R^+ \cup \{0\}$.
    a. Find dimension of V over S.
    b. Write V as a direct sum.

84. Obtain some interesting properties enjoyed by semivector space of refined labels over a semifield.

85. Prove $L_{R^+ \cup \{0\}} = \{L_{a_i} \mid a_i \in R^+ \cup \{0\}\}$ is a semiring. Is that a semifield? justify your answer.

86. Obtain some interesting properties about DSm semivector spaces of refined labels defined over the semifield $S = Q^+ \cup \{0\}$.

87. Let $V = \left\{ \begin{bmatrix} L_{a_1} & L_{a_2} & L_{a_3} \\ L_{a_4} & L_{a_5} & L_{a_6} \\ L_{a_7} & L_{a_8} & L_{a_9} \end{bmatrix} \middle| L_{a_i} \in L_{R^+ \cup \{0\}}; 1 \leq i \leq 9 \right\}$ be a DSm semivector space of refined labels over the semifield $S = Z^+ \cup \{0\}$.
   a. Find DSm semivector subspaces of refined labels.
   b. What is dimension of V?
   c. Find a linear operator on V.
   d. Is V a DSm semilinear algebra of refined labels over S?

88. Let $V = \left\{ \sum_{i=0}^{\infty} L_{a_i} x^i \middle| L_{a_i} \in L_{R^+ \cup \{0\}} \right\}$ be a DSm semilinear algebra of refined labels over the semifield $S = Z^+ \cup \{0\}$.
   a. Can we write V as a direct sum of DSm semilinear subalgebras?
   b. Find a basis for V.
   c. Is V a S-semilinear algebra of refined labels over the semifield $S = Z^+ \cup \{0\}$?
   d. Find ideals in V.
   e. Does V contain annihilator ideals?
   f. Study the collection of linear operators from V to V. What is the algebraic structure enjoyed by V?

89. Obtain some interesting properties enjoyed by DSm semilinear algebra of refined labels V built using square matrices with entries from $L_{R^+ \cup \{0\}}$ over $S = R^+ \cup \{0\}$.

90. Let $V = \left\{ \begin{bmatrix} L_{a_1} & L_{a_2} & L_{a_3} & L_{a_4} \\ L_{a_5} & L_{a_6} & L_{a_7} & L_{a_8} \\ L_{a_9} & L_{a_{10}} & L_{a_{11}} & L_{a_{12}} \\ L_{a_{13}} & L_{a_{14}} & L_{a_{15}} & L_{a_{16}} \\ L_{a_{17}} & L_{a_{18}} & L_{a_{19}} & L_{a_{20}} \end{bmatrix} \middle| L_{a_i} \in L_{R^+ \cup \{0\}}; 1 \le i \le 20 \right\}$ be a

DSm semivector space of refined labels over the semifield $S = R^+ \cup \{0\}$. Consider

$M = \left\{ \begin{bmatrix} L_{a_1} & L_{a_2} & L_{a_3} & L_{a_4} \\ 0 & 0 & 0 & L_{a_5} \\ 0 & 0 & L_{a_6} & 0 \\ 0 & L_{a_7} & 0 & 0 \\ L_{a_8} & 0 & 0 & 0 \end{bmatrix} \middle| L_{a_i} \in L_{R^+ \cup \{0\}}; 1 \le i \le 8 \right\} \subseteq V$, M is a

DSm semivector space of refined labels over the semifield S.
 a. Find a T: V → V so that T(W) ⊆ W.
 b. Write $V = \bigcup_i W_i$ as a direct sum.
 c. Write $V = \bigcup_i S_i$ as a pseudo direct sum.
 d. Find a basis of V.
 e. Find a linearly independent subset of V.
 f. Does there exist any relation between the questions (d) and (e)?
 g. Find V/W. Is V/W a DSm semivector space of refined labels over S?
 h. Give a T: V → V such that $T^{-1}$ exists.
 i. Give T : V → V so that $T^{-1}$ does not exist.

91. Can Smarandache DSm vector spaces of refined labels be defined?

92. Can Smarandache DSm semivector space of refined labels be defined?

93. Determine some interesting properties about unordered labels or partially ordered labels of finite order.

94. Prove the set of partially ordered labels is always a lattice ($L_0$ and $L_{m+1}$ adjoined).

95. Does the lattice of labels being distributive imply any special property enjoyed by the label?
    ($L_{m+1}$ the maximal label and $L_0$ the minimal label adjoined)

96. Can one say if the ordinary labels associated with a model is non distributive it enjoys special properties?

97. Study the ordinary lattice of labels which is modular.

98. Prove all totally ordered labels with zero and maximal element is isomorphic with the chain lattice.

99. Prove that such labels or these chain lattices are semifields of finite order.

100. Prove using these semifields we can build semivector spaces of finite order.

101. Give examples of finite semivector spaces of ordinary labels.

102. Find some applications of DSm semiring of ordinary labels.

103. What are the possible applications of DSm semifield of refined labels?

104. Study the applications of DSm semivector space of ordinary labels.

105. Obtain some interesting results related with the DSm semivector spaces of matrices built using ordinary labels over the semifield $S = \{L_0 = 0, L_1, L_2, \ldots, L_m, L_{m+1} = 1\}$

106. Let $V = \left\{ \sum_{i=0} L_{a_i} x^i \mid L_{a_i} \in \{0 = L_0, L_1, L_2, \ldots, L_m, L_{m+1} = 1\} \right\}$
    be the DSm semivector space of ordinary labels over the semifield $S = \{0 = L_0, L_1, \ldots, L_m, L_{m+1} = 1\}$. Determine the important properties enjoyed by these algebraic structure. Can these be applied to eigen value problems?

107. Let $V = \left\{ \begin{bmatrix} L_{a_1} & L_{a_2} & L_{a_3} \\ L_{a_4} & L_{a_5} & L_{a_6} \\ L_{a_7} & L_{a_8} & L_{a_9} \end{bmatrix} \mid L_{a_i} \in \{0 = L_0, L_1, \ldots, L_{m+1} = 1\}; 1 \leq i \leq 9 \right\}$
    be a DSm semivector space of ordinary labels over the semifield $S = \{0, L_1, \ldots, L_m, L_{m+1} = 1\}$.
    a. Is V a DSm semilinear algebra of ordinary labels over S?
    b. Find a basis of V.
    c. Is V finite dimensional?

108. Find some interesting properties enjoyed by $V = \left\{ \sum_{i=0}^{10} L_{a_i} x^i \mid L_{a_i} = L = \{0 = L_0, L_1, \ldots, L_{m+1} = 1\}; 0 \leq i \leq 10 \right\}$ be
    the DSm semivector space of ordinary labels over $L = \{0 = L_0, L_1, \ldots, L_m, L_{m+1} = 1\}$.
    a. Is V finite dimensional over L?
    b. Find some special and important properties enjoyed by V.

109. Let $V = \left\{ \begin{bmatrix} L_{a_1} \\ L_{a_2} \\ \vdots \\ L_{a_9} \end{bmatrix}, \begin{bmatrix} L_{a_1} & L_{a_2} \\ L_{a_3} & L_{a_4} \\ L_{a_5} & L_{a_6} \end{bmatrix}, (L_{a_1}, L_{a_2}, ..., L_{a_6}) \right.$

$\left. L_{a_i} \in \{0 = L_0, L_1, L_2, ..., L_m, L_{m+1} = 1\}; 1 \leq i \leq 9 \right\}$

be a DSm set vector space of ordinary labels over the set $S = \{0, L_1, L_3 \text{ and } L_m, L_{m+1} = 1\}$.
   a. Write V as a direct sum or union.
   b. Write V as a pseudo direct sum.
   c. Find a non invertible operator on V.
   d. Does V contain only a finite number of elements?

110. Let $V = \left\{ \begin{bmatrix} L_{a_1} & L_{a_2} \\ L_{a_3} & L_{a_4} \end{bmatrix}, (L_{a_1}, L_{a_2}, ..., L_{a_{20}}), \begin{bmatrix} L_{a_1} \\ L_{a_2} \\ \vdots \\ L_{a_{12}} \end{bmatrix} \right.$

$\left. L_{a_i} \in L = \{0 = L_0, L_1, ..., L_{m+1} = 1\}; 1 \leq i \leq 20 \right\}$ be a

semigroup DSm vector space of ordinary labels over the semigroup $S = \{0 = L_0, L_1, ..., L_m, L_{m+1} = 1\}$.
   a. Can V have DSm semigroup vector subspaces of ordinary labels?
   b. Can V have infinite number of elements?
   c. Does V have a invertible linear operator?
   d. Write V as a direct union.

111. Is it possible to construct group DSm vector space of ordinary labels? Justify your answer.

112. Can $L = \{0 = L_0, L_1, ..., L_m, L_{m+1} = 1\}$ be given a group structure L an ordinary label?

113. Can the unordered set of ordinary labels be given a group structure? Justify your claim.

114. Will a partially ordered ordinary labels be a group? Justify.

115. Prove the partially ordered ordinary label is a lattice. Give an example of a partially ordered ordinary lattice which is distributive.

116. Let
$$V = \left\{ \begin{bmatrix} L_{a_1} & L_{a_9} \\ L_{a_2} & L_{a_{10}} \\ \vdots & \vdots \\ L_{a_8} & L_{a_{16}} \end{bmatrix}, \begin{bmatrix} L_{a_1} \\ L_{a_2} \\ \vdots \\ L_{a_{20}} \end{bmatrix}, \begin{bmatrix} L_{a_1} & L_{a_2} & L_{a_3} \\ L_{a_4} & L_{a_5} & L_{a_6} \\ L_{a_7} & L_{a_8} & L_{a_9} \end{bmatrix}, \begin{bmatrix} L_{a_1} & L_{a_2} \\ L_{a_3} & L_{a_4} \end{bmatrix} \right\}$$

$L_{a_i} \in \{0 = L_0, L_1, ..., L_{m+1} = 1\}; 1 \le i \le 20\}$ be a DSm set vector space over the set $S = \{L_0 = 0, L_1, ..., L_m, L_{m+1} = 1\}$.
  a. Find the cardinality of V.
  b. Write $V = \cup W_i$ as a direct sum.
  c. Write $V = \bigcup_i W_i$, as a pseudo direct sum.
  d. Find $T : V \to V$ so that T is invertible.

117. Let $W = \left\{ \begin{bmatrix} L_{a_1} & L_{a_2} & L_{a_3} \\ L_{a_4} & L_{a_5} & L_{a_6} \\ L_{a_7} & L_{a_8} & L_{a_9} \end{bmatrix} \middle| \begin{array}{l} L_{a_i} \in \{0 = L_0, L_1, ..., L_{m+1} = 1\}; \\ 1 \le i \le 9 \end{array} \right\}$

be a semigroup DSm vector space of ordinary labels over the semigroup $S = \{0 = L_0, L_1, ..., L_m, L_{m+1} = 1\}$
  a. Is V a semigroup DSm linear algebra of refined labels?
  b. Let $T: V \to V$, find a linear operator on V which is invertible.
  c. Write $W = \cup P_i$ as a direct sum.

118. Let $V = \left\{ \begin{bmatrix} L_{a_1} \\ L_{a_2} \\ L_{a_3} \\ L_{a_4} \end{bmatrix} \middle| L_{a_i} \in \{0 = L_0, L_1, L_2, ..., L_m, L_{m+1} = 1\} \right\}$ be a semigroup of DSm vector space over the semigroup $S = \{0 = L_0, L_1, ..., L_m, L_{m+1} = 1\}$.
   a. What is the dimension of V?
   b. Find number of elements in V.
   c. Write $V = \bigcup_i L_i$ as direct sum.

119. Find some interesting applications of DSm semigroup vector space of ordinary labels over the semigroup $S = \{0 = L_0, L_1, ..., L_m, L_{m+1} = 1\}$.

120. Can unordered ordinary labels be used as fuzzy models / dynamical systems using the matrices built using $L = \{0 = L_0, L_1, ..., L_m, L_{m+1} = 1\}$?

121. Can these semifield of ordinary labels be used in web designing?

122. Can the concept lattices be built using lattices of ordinary labels?

123. Study the ordinary labels as a) lattices b) fields c) semilattices.

124. Show if

$$S = \left\{ \begin{bmatrix} L_{a_1} \\ L_{a_2} \\ L_{a_3} \\ L_{a_4} \end{bmatrix} \middle| L_{a_i} \in \{0 = L_0, L_1, ..., L_{m+1} = 1\}; i = 1, 2, 3, 4 \right\}$$

be a DSm semivector space over the semifield $S = \{0 = L_0, L_1, \ldots, L_m, L_{m+1} = 1\}$.

   a. Is S is a lattice?
   b. Can S be totally ordered?
   c. What is the order of S?

125. Distinguish between the algebraic structures enjoyed by the refined labels and the ordinary labels.

126. Study the possible algebraic structures enjoyed by

$$M = \left\{ \begin{bmatrix} L_{a_{11}} & \cdots & L_{a_{1n}} \\ L_{a_{21}} & \cdots & L_{a_{2n}} \\ \vdots & & \vdots \\ L_{a_{n1}} & \cdots & L_{a_{nn}} \end{bmatrix} \middle| L_{a_i} \in \{L_R, +, \times\} \right\}.$$

   a. Is M a ring?
   b. Can M be a commutative ring?
   c. Can M have zero divisors?

127. Let

$$P = \left\{ \begin{bmatrix} L_{a_{11}} & \cdots & L_{a_{1m}} \\ L_{a_{21}} & \cdots & L_{a_{2m}} \\ \vdots & & \vdots \\ L_{a_{n1}} & \cdots & L_{a_{nm}} \end{bmatrix} \middle| L_{a_{ij}} \in \{L_R, +, \times\} (m \neq n) \right\};$$

   a. What is the possible algebraic structure P can enjoy?
   b. Can P become a ring?
   c. Is P a an additive abelian group?

128. Let

$$L = \left\{ \begin{bmatrix} L_{a_1} & L_{a_2} \\ L_{a_3} & L_{a_4} \\ L_{a_5} & L_{a_6} \\ L_{a_7} & L_{a_8} \end{bmatrix} \middle| \begin{array}{c} L_{a_i} \in \{0 = L_0 < L_1 < ... < L_{m+1} = 1\}; \\ 1 \le i \le 8 \end{array} \right\}$$

a. Can L be a partially ordered set?
(under the ordering if

$$M = \begin{bmatrix} L_{a_1} & L_{a_2} \\ L_{a_3} & L_{a_4} \\ L_{a_5} & L_{a_6} \\ L_{a_7} & L_{a_8} \end{bmatrix}$$

and

$$N = \begin{bmatrix} L_{b_1} & L_{b_2} \\ L_{b_3} & L_{b_4} \\ L_{b_5} & L_{b_6} \\ L_{b_7} & L_{b_8} \end{bmatrix}$$

are in L. M > N if each $L_{a_i} > L_{b_i}$; $1 \le i \le 8$)).

b. Can L be a semilattice under '∪' or '∩'?
c. Can L be a lattice?
d. Can L be a distributive lattice? (justify)

129. Let

$$T = \left\{ \begin{bmatrix} L_{a_1} & L_{a_2} & L_{a_3} \\ L_{a_4} & L_{a_5} & L_{a_6} \\ L_{a_7} & L_{a_8} & L_{a_9} \end{bmatrix} \middle| L_{a_i} \in \{L_R\}; 1 \le i \le 9 \right\}$$

be a DSm linear algebra over $L_R$.
a. Find a basis of T over $L_R$.
b. What is dimension of T?
c. Write $T = \cup W_i$ as a direct sum.

130. Study DSm vector spaces V defined over the DSm field $L_R$. Is V the same if $L_R$ is replaced by R? Justify your claim.

# FURTHER READING

# INDEX

# ABOUT THE AUTHORS

**Dr.W.B.Vasantha Kandasamy** is an Associate Professor in the Department of Mathematics, Indian Institute of Technology Madras, Chennai. In the past decade she has guided 13 Ph.D. scholars in the different fields of non-associative algebras, algebraic coding theory, transportation theory, fuzzy groups, and applications of fuzzy theory of the problems faced in chemical industries and cement industries. She has to her credit 646 research papers. She has guided over 68 M.Sc. and M.Tech. projects. She has worked in collaboration projects with the Indian Space Research Organization and with the Tamil Nadu State AIDS Control Society. She is presently working on a research project funded by the Board of Research in Nuclear Sciences, Government of India. This is her $57^{th}$ book.

On India's 60th Independence Day, Dr.Vasantha was conferred the Kalpana Chawla Award for Courage and Daring Enterprise by the State Government of Tamil Nadu in recognition of her sustained fight for social justice in the Indian Institute of Technology (IIT) Madras and for her contribution to mathematics. The award, instituted in the memory of Indian-American astronaut Kalpana Chawla who died aboard Space Shuttle Columbia, carried a cash prize of five lakh rupees (the highest prize-money for any Indian award) and a gold medal.
She can be contacted at vasanthakandasamy@gmail.com
Web Site: http://mat.iitm.ac.in/home/wbv/public_html/
or http://www.vasantha.in

**Dr. Florentin Smarandache** is a Professor of Mathematics at the University of New Mexico in USA. He published over 75 books and 200 articles and notes in mathematics, physics, philosophy, psychology, rebus, literature.

In mathematics his research is in number theory, non-Euclidean geometry, synthetic geometry, algebraic structures, statistics, neutrosophic logic and set (generalizations of fuzzy logic and set respectively), neutrosophic probability (generalization of classical and imprecise probability). Also, small contributions to nuclear and particle physics, information fusion, neutrosophy (a generalization of dialectics), law of sensations and stimuli, etc. He can be contacted at smarand@unm.edu